\documentclass[12pt,a4paper]{article}
\usepackage{amssymb,amsmath,amsthm,euscript}
\usepackage[matrix,arrow]{xy}
\usepackage{hyperref}

\usepackage[utf8]{inputenc}
\usepackage[english]{babel}

\newtheorem{Th}{Theorem}[section]
\newtheorem{Lem}[Th]{Lemma}
\newtheorem{Prop}[Th]{Proposition}

\newtheorem{Ex}{Example}[section]
\newtheorem{Def}{Definition}[section]
\newtheorem{Rem}{Remark}[section]

\newcommand{\End}{\mathop{\mathrm{End}}\nolimits}
\newcommand{\Res}{\mathop{\mathrm{Res}}\nolimits}
\newcommand{\res}{\mathop{\mathrm{res}}\nolimits}

\newcommand{\Mer}{\mathop{\mathrm{Mer}}\nolimits}

\def\beq#1#2\eeq{%
        \begin{equation}%
        \label{#1}%
            #2%
        \end{equation}%
    }

\newcommand{\p}{\partial}

\renewcommand{\a}{\alpha}

\newcommand{\nad}[2]{\genfrac{}{}{0pt}{}{#1}{#2}}

\numberwithin{equation}{section}

\renewcommand{\le}{\leqslant}
\renewcommand{\ge}{\geqslant}
\newcommand{\wt}{\widetilde}
\newcommand{\wh}{\widehat}

\newcommand{\T}{{\mathcal T}}
\renewcommand{\SS}{{\mathbb S}}
\newcommand{\HH}{\hbox to 15pt {${\cal H}$\hss${\cal H}$}}

\newcommand{\mref}[1]{(\ref{#1})}

\begin{document}

\author{M.\,Feigin\thanks{Misha.Feigin@glasgow.ac.uk}, A.\,Silantyev\thanks{Alexey.Silantyev@glasgow.ac.uk }}
\title{Generalized Macdonald-Ruijsenaars systems}
\date{}

\maketitle

\begin{center}
School of Mathematics and Statistics, University of Glasgow,

University Gardens, Glasgow G12 8QW, UK
\end{center}

\vspace{4mm}

\begin{abstract}
\noindent We consider the polynomial representation of Double Affine Hecke Algebras (DAHAs) and construct its submodules as ideals of functions vanishing on the special collections of affine planes. This generalizes certain results of Kasatani in types $A_n$, $(C_n^\vee,C_n)$.
We obtain commutative algebras of difference operators given by the action of invariant combinations of Cherednik-Dunkl operators in the corresponding quotient modules of the polynomial representation. This gives known and new generalized Macdonald-Ruijsenaars systems.  Thus in the cases of DAHAs of types $A_n$ and $(C_n^\vee,C_n)$ we derive Chalykh-Sergeev-Veselov operators
and a generalization of the Koornwinder operator respectively, together with complete sets of quantum integrals in the explicit form.

\end{abstract}

\begin{keywords}
 Double affine Hecke algebra; polynomial representation; wheel condition; Macdonald operator;  integrable system
\end{keywords}

\tableofcontents

\section{Introduction}

The main goal of this paper is to obtain commutative algebras of difference operators containing generalizations of Macdonald-Ruijsenaars operators using special representations of double affine Hecke algebras (DAHAs). In the work~\cite{Ruijs} the Macdonald-Ruijsenaars operator related to the root system $R=A_{N-1}$ appeared as the relativistic version of quantum Calogero-Moser system, this operator has the form
\begin{align} \label{MR}
 H=\sum_{i=1}^N\prod_{\nad{j=1}{j\ne i}}^N \frac{e^{x_i}- e^{x_j+\eta}}{e^{x_i}-e^{x_j}}\T^\hbar_{x_i},
\end{align}
where $\T^\hbar_{x_i}= e^{\hbar \p_{x_i}}$ is the shift operator and $\eta, \hbar$ are parameters. Macdonald derived and considered operator~\eqref{MR} and its generalizations for any root system and (quasi-)minuscule coweight in the context of multi-variable orthogonal polynomials~\cite{Mac2}.
Commutative algebras of difference operators containing the operators constructed by Macdonald were realized by Cherednik with the help of his DAHA and its polynomial representation~\cite{Cher92DAHA,Cher95}. These algebras are centralizers of the Macdonald-Ruijsenaars operators in the algebras of Weyl-invariant difference operators~\cite{LS}. The whole area of DAHA is a very rich mathematical subject with deep connections and interactions with combinatorics, theory of (quantum) symmetric spaces, representation theory, noncommutative geometry~\cite{CherDAHA}.

The first generalization of the Macdonald-Ruijsenaars operator which is not directly associated to a root system was found by Chalykh~\cite{Chalykhbisp}. The operator is a difference version of the Hamiltonian which describes interaction of the system of Calogero-Moser particles of two possible masses~\cite{CFVjmt}. The approach of~\cite{Chalykhbisp} (see also~\cite{ChalykhAdv}) is based on the multidimensional Baker-Akhiezer functions, in this case all but one particles have equal masses. The generalization of Chalykh's operator to the system of $N_1+N_2$ particles of two types was introduced by Sergeev and Veselov~\cite{SV}, the corresponding operator has the form
\begin{multline}\label{MRdeform}
H_{N_1, N_2}= \sum_{i=1}^{N_1}
\prod_{j=1\atop j\ne i}^{N_1}
\frac{e^{x_i}- e^{x_j+\eta}}{e^{x_i}-e^{x_j}}
\prod_{\ell=1}^{N_2}
\frac{e^{x_i}- e^{y_\ell+\hbar}}{e^{x_i}-e^{y_\ell}}
\;\T^\hbar_{x_i}\,+ \\
+\frac{e^\hbar-1}{e^\eta-1}
\sum_{\ell=1}^{N_2}
  \prod_{i=1}^{N_1}
\frac{e^{y_\ell}- e^{x_i+\eta}}{e^{y_\ell}-e^{x_i}}
  \prod_{m=1\atop m\ne\ell}^{N_2}
\frac{e^{y_\ell}- e^{y_m+\hbar}}{e^{y_\ell}-e^{y_m}}
\;\T^{\eta}_{y_\ell}.
\end{multline}
In the work \cite{SV} the operator~\eqref{MRdeform} was obtained from the action of the Macdonald-Ruijsenaars operator~\eqref{MR} extended to the space of symmetric functions in the infinite number of variables with subsequent restriction of the operator~\eqref{MR} to the special discriminant subvariety. Further rational generalized Macdonald-Ruijsenaars operators were derived in~\cite{Fbisp}, \cite{SVBC} starting from the generalized Calogero-Moser systems and using the bispectrality.

In the present work we develop a uniform approach to the generalized Macdonald-Ruijsenaars systems based on the special representations of DAHAs.
DAHA can be associated to any affine root system~\cite{Mac2}. The main cases are DAHA of type $\wh{R}$ and of type $\wh{R}^\vee$ where
$R$ is a reduced irreducible root system, $\wh{R}$ is its affinization and $\wh{R}^\vee$ denotes the dual affine root system (see Cherednik's papers~\cite{Cher93} and \cite{Cher92DAHA} respectively), and DAHA of type $(C_n^\vee, C_n)$ (see~\cite{Noumi}, \cite{Sahi}).
The original Cherednik's construction \cite{Cher92DAHA,Cher93} allowed to obtain Macdonald-Ruijsenaars systems for any root system $R$ from the polynomial representation of the corresponding DAHA.  The difference operators from the commutative algebra are obtained by the averaging of the Cherednik-Dunkl operators over the Weyl group with subsequent restriction of the result to the space of invariant polynomials. In~\cite{Noumi} Noumi generalized the construction for the $(C_n^\vee,C_n)$ case which allowed to obtain commutative algebra containing the Koornwinder operator~\cite{Koornwinder}.
In order to obtain generalized Macdonald-Ruijsenaars systems we first specialize DAHA-invariant ideals in the polynomial representation, these ideals consist of functions vanishing on certain collections of planes. Then the action of the invariant combinations of Cherednik-Dunkl operators is defined on the functions on these planes. Upon restriction to Weyl group invariants we obtain commutative algebras of difference operators. In this way we recover operator $\eqref{MRdeform}$ (under assumption that $\hbar/\eta \in \mathbb{Z}$) for the case of DAHA of type $A$ and derive further known and new generalized Macdonald-Ruijsenaars operators starting from other root systems. Similar approach to the generalized Calogero-Moser systems was developed by one of the authors in~\cite{F}. It is based on finding special ideals invariant under the rational Cherednik algebras with subsequent restriction of invariant combinations of Dunkl operators thus generalizing Heckman's construction of the usual Calogero-Moser systems~\cite{Heckman}.

The idea to describe submodules of the polynomial representations of DAHAs by vanishing conditions was very successfully used by Kasatani in~\cite{KasA}, \cite{KasCC} (see also discussion in the introduction of~\cite{CherNonsemisimple} where the polynomial representation is studied from another perspective). Thus in the work~\cite{KasA} a filtration of the polynomial representation in type $A$ by the ideals given by vanishing (``wheel") conditions was introduced so that it provides composition factors of the polynomial module (see also~\cite{Enomoto}). The paper~\cite{KasCC} deals with the polynomial representation of DAHA associated to the affine root system $(C_n^\vee,C_n)$. Ideals in the ring of symmetric polynomials invariant under the Macdonald-Ruijsenaars operator~(\ref{MR}) were earlier defined by wheel conditions in \cite{FJMM}.

The structure of the paper is as follows. In Section~\ref{sec2} we recall Cherednik's construction of the commutative algebra of difference operators using DAHA of type $\wh{R}$. In Section~\ref{sec3}  we introduce arrangements of planes $\cal D$ corresponding to subsystems $S_0$ in the root system $R$ and we define ideals of functions $I_{\cal D}$ vanishing on these planes. Our first main result (Theorem~\ref{ThInv}) states that for the special geometry of the subsystem $S_0$ and for the related choice of parameters for DAHA the ideals $I_{\cal D}$ are DAHA-invariant. The arising arrangements $\cal D$ are shifted versions of the linear arrangements appearing in the study of invariant ideals for the rational Cherednik algebras~\cite{F}. In Subsection~\ref{further-ideals} we construct further invariant ideals given by vanishing conditions where we allow parallel planes. This generalizes construction of such ideals for DAHA for classical root systems from~\cite{KasA,KasCC}.

In Section~\ref{sec4} we present our next main result (Theorem~\ref{teorogr}) which explains how commutative algebras of difference operators arise in the quotients of the polynomial representations by the above ideals $I_{\cal D}$. We show in Subsection~\ref{sec42} that these algebras contain complete sets of algebraically independent operators. We verify in Subsection~\ref{sec43} that in type $A$ we obtain operators~\eqref{MRdeform}. We also construct explicitly a complete set of commuting operators for~\eqref{MRdeform} (Theorem~\ref{PropHOMOcomm}). We consider further examples illustrating the approach in Subsections~\ref{sec44},~\ref{sec46}. Thus we demonstrate that starting from the Macdonald operator of type $B$ we obtain trigonometric version of the operators known from~\cite{SVBC}, \cite{Fbisp}. Subsystems in the exceptional root systems provide new families of commutative algebras, we demonstrate this by getting a restricted operator starting from the Macdonald operator for the root system $R=F_4$.

In Section~\ref{sec5} we extend our constructions for the $(C_n^\vee,C_n)$ DAHA. We specify invariant ideals in the polynomial representation by imposing vanishing condition on suitable arrangements of planes (Theorem~\ref{ThInvCC}). These submodules correspond to the intersections of the submodules investigated in~\cite{KasCC}. We introduce generalized Koornwinder operator and derive it as a restriction of the Koornwinder operator in Propositions~\ref{Koornwinder-generalized-1},~\ref{Koornwinder-generalized-2}. This operator depends on six parameters like the original Koornwinder operator, and also it reflects that there are interacting particles of two types.  Although the restriction construction imposes some integrality relation(s) on the parameters we show that the operator can be embedded into a large commutative algebra of difference operators for any values of the parameters. In Theorem \ref{PropHOMOcommCC} we find a complete set of quantum integrals in the explicit form. This generalizes van Diejen's quantum integrals for the Koornwinder operator \cite{vD96}.

In Section~\ref{sec62} we present our constructions for  DAHA of type $\wh{R}^\vee$ \cite{Cher92DAHA}. The invariant ideals are specified in Theorem~\ref{ThInvvee}. \\


{\bf Acknowledgements.} We are grateful to M. Noumi, J. Stokman and A.P. Veselov for stimulating discussions. We are very grateful to J. Shiraishi for interesting discussions which led us to Theorem~\ref{PropHOMOcomm}.
The work of both authors was supported by the EPSRC grant EP/F032889/1. M.F. also acknowledges support of the British Council (PMI2 Research Cooperation award). The work of A.S. was partially supported by the RFBR grant 09-01-00239-a. A.S. also thanks LAREMA and Department of Mathematics of the University of Angers for the hospitality during the visit.

\section{Macdonald-Ruijsenaars systems from DAHA}
\label{sec2}

In this section we recall the construction of the Macdonald-Ruijsenaars systems from the Cherednik-Dunkl operators. We consider here DAHAs of type $\wh R$ (other cases will be considered in Sections~\ref{sec5},~\ref{sec62} and Appendix~\ref{sec6}). Following~\cite{Cher93,Cher95,CherDAHA,Kir} we define DAHA for the $\wh R$ case, then we define the Cherednik-Dunkl operators as elements of DAHA acting in the polynomial representation and describe the Macdonald-Ruijsenaars systems in these terms.

\subsection{Double Affine Hecke Algebra}
\label{sec21}

Let $R\subset\mathbb R^n\subset V=\mathbb C^n$ be an irreducible reduced root system  of rank $n$. Let $(\cdot,\cdot)$ be the standard bilinear symmetric non-degenerated form in $V$ and let $\beta^\vee=2\beta/(\beta,\beta)$. We have  $(\alpha,\beta^\vee)\in\mathbb Z$ for all $\alpha,\beta\in R$. Consider the orthogonal reflections $s_\alpha$ corresponding to the roots $\alpha\in R$, they act on $x\in V$ as $s_\alpha(x)=x-(\alpha,x)\alpha^\vee$. The group $W$ generated by the orthogonal reflections $s_\alpha$, $\alpha\in R$, is a Weyl group.

Let $R_+$ be a set of positive roots and $\Delta=\{\alpha_1,\ldots,\alpha_n\}\subset R_+$ be the set of simple roots. The reflections $s_{\alpha_1},\ldots,s_{\alpha_n}$ generate the Weyl group $W$. Let us introduce lattices of coroots, of coweights and of weights:
\begin{align}
 Q^\vee&=\bigoplus\nolimits_{i=1}^n\mathbb Z\,\alpha_i^\vee, &
 P^\vee&=\bigoplus\nolimits_{i=1}^n\mathbb Z\,b_i, &
 P&=\bigoplus\nolimits_{i=1}^n\mathbb Z\,\omega_i,
\end{align}
where $b_i$ and $\omega_i$ are fundamental coweights and fundamental weights:
\begin{align}
 (\alpha_i,b_j)&=\delta_{ij}, & (\omega_i,\alpha_j^\vee)&=\delta_{ij}.
\end{align}
Introduce also the subsets of dominant and anti-dominant coweights:
\begin{align*}
P^\vee_+=\{\lambda\in P^\vee\mid(\alpha,\lambda)\ge0, \forall\alpha\in R_+\},&&& P^\vee_-=\{\lambda\in P^\vee\mid(\alpha,\lambda)\le0, \forall\alpha\in R_+\}.
\end{align*}

Let $\wh V$ be the space of affine functionals on $V$. Each functional $\alpha\in\wh V$ acts on $x\in V$ as $\alpha(x)=(\alpha^{(1)},x)+\alpha^{(0)}$, where $\alpha^{(1)}\in V$ and $\alpha^{(0)}\in\mathbb C$, we write $\alpha=[\alpha^{(1)},\alpha^{(0)}]$ in this case. The roots $\gamma\in R$ can be considered as elements of $\wh V$ acting as $\gamma(x)=(\gamma,x)$. Define the bilinear form on the space $\wh V$ by the formula $(\alpha,\beta)=(\alpha^{(1)},\beta^{(1)})$, where $\alpha=[\alpha^{(1)},\alpha^{(0)}]$, $\beta=[\beta^{(1)},\beta^{(0)}]$, and $\alpha^\vee=\frac{2\alpha}{(\alpha,\alpha)}=[\frac{2\alpha^{(1)}}{(\alpha^{(1)},\alpha^{(1)})},\frac{2\alpha^{(0)}}{(\alpha^{(1)},\alpha^{(1)})}]$ for any $\alpha=[\alpha^{(1)},\alpha^{(0)}]\in\wh V$.

The reflection $s_\alpha$ corresponding to $\alpha=[\alpha^{(1)},\alpha^{(0)}]\in \wh V$ is the following affine map: $s_\alpha(x)=x-\alpha(x)(\alpha^{(1)})^\vee$. It is the orthogonal reflection with respect to the hyperplane $H(\alpha)=\{x\in V\mid \alpha(x)=0\}$. The action $s_\alpha$ on $V$ induces the action on $\wh V$ given by the formula $s_\alpha(\beta)=\beta-(\beta,\alpha^\vee)\alpha$, where $\alpha,\beta\in \wh V$.

Let $\hbar\in\mathbb C$ be a non-zero complex parameter. Introduce the functional $\alpha_0=[-\theta,\hbar]$, where $\theta\in R$ is the maximal root for $R$ with respect to the basis $\Delta$. For $\lambda\in V$ denote by $\tau(\lambda)$ the shift operator $\tau(\lambda)x=x+\hbar\lambda$, $\forall x\in V$.

{\it Affine Weyl group} $\wt W$ is generated by the reflections $s_\alpha$, with $\alpha=[\alpha^{(1)},\alpha^{(0)}]$ such that $\alpha^{(1)}\in R$ and $\alpha^{(0)}\in\hbar\mathbb Z$. The set $\wh R=\{[\alpha^{(1)},\hbar k]\mid \alpha^{(1)}\in R, k\in\mathbb Z\}$ is the affine root system corresponding to the root system $R$~\cite{Mac1}. The group $\tau(Q^\vee)$ is a normal subgroup in $\wt W$ and $\wt W=W\ltimes \tau(Q^\vee)$. The group $\wt W$ is generated by the reflections $s_0,s_1,\ldots,s_n$, where $s_i=s_{\alpha_i}$ for $i=0,1,\ldots,n$, (see~\cite{H}). The subset $\wt\Delta=\{\alpha_0,\alpha_1,\ldots,\alpha_n\}\subset\wh R$ is a basis of $\wh V$ called the set of affine simple roots. For $0\le i, j\le n$, $i\ne j$ we denote by $m_{ij}=m_{ji}$ the order (finite or infinite) of the element $s_is_j\in\wt W$.

{\it Extended affine Weyl group} is defined as $\wh W=W\ltimes\tau(P^\vee)$ with the commutation relations $w\tau(\lambda)=\tau(w\lambda)w$, where $w\in W$, $\lambda\in P^\vee$. The action of the group $\wh W$ on $V$ induces the following linear action on the space $\wh V$:
\begin{align}
 w\tau(\lambda)(\alpha)&=[w(\alpha^{(1)}),\alpha^{(0)}-\hbar(\alpha^{(1)},\lambda)],
\end{align}
where $w\in W$, $\lambda\in P^\vee$, $\alpha=[\alpha^{(1)},\alpha^{(0)}]$.

Let $\Pi=\{\pi\in\wh W\mid \pi\wt\Delta=\wt\Delta\}\subset\wh W$ be a subgroup of elements preserving the set of simple affine roots. Then $\wh W=\Pi\ltimes\wt W$. The group $\Pi$ is finite: $\Pi=\{\pi_r\mid r\in{\cal O}\}$, where ${\cal O}\subset\{0,1,\ldots,n\}$ and each element $\pi_r$ is uniquely determined by the condition $\pi_r(\alpha_0)=\alpha_r$. They have the form $\pi_0=1$, $\pi_r=\tau(b_r)u_r^{-1}$ for $r\in{\cal O}'$, where ${\cal O}'={\cal O}\backslash\{0\}$ and $u_r\in W$ are transformations preserving the set $\{-\theta,\alpha_1,\ldots,\alpha_n\}$ and satisfying $u_r(\alpha_r)=-\theta$. The corresponding fundamental coweights $b_r$ form the entire set of {\it minuscule coweights}, that is $\{b_r\mid r\in{\cal O}'\}=\{\lambda\in P^\vee\mid0\le(\alpha,\lambda)\le1, \forall\alpha\in R_+\}$.

Let $\wh P=\{\mu=[\mu^{(1)},\mu^{(0)}]\mid\mu^{(1)}\in P,\,\mu^{(0)}\in\frac\hbar{\epsilon_0}\mathbb Z\}\subset \wh V$, where $\epsilon_0$ is an integer such that $(P,P^\vee)\subset\epsilon_0^{-1}\mathbb Z$. We consider the weight lattice $P$ as a sublattice of $\wh P$. Now we define DAHA corresponding to the affine system $\wh R$. Let $t\colon\wh R\to\mathbb C\backslash\{0\}$ be a $\wh W$-invariant function on $\wh R$, that is $t_\beta=t_{\hat w(\beta)}$, for all $\beta\in\wh R$, $\hat w\in\wh W$, where we denote $t_\beta=t(\beta)$. We also denote $t_i=t_{\alpha_i}$ for $i=0,1,\ldots,n$, these values uniquely determine the function $t$. Let $q=e^{\hbar/2}$. We will assume that $\hbar\notin\pi i\mathbb Q$, so $q$ is not a root of unity.

\begin{Def} {\normalfont(\cite{Cher93}, see also~\cite{Kir}).} 
 Double Affine Hecke Algebra (for $\wh{R}$ case) is defined as the algebra $\HH_t=\HH_{\hbar,t}^{\wh R}$ generated by the elements $T_0,T_1,\ldots,T_n$, by the set of pairwise commuting elements $\{X^\mu\}_{\mu\in\hat P}$, and by the group $\Pi=\{\pi_r\mid r\in{\cal O}\}$ with the relations{\normalfont
\begin{align}
 X^\mu X^{\mu'}&=X^{\mu+\mu'},                & &\mu,\mu'\in\wh P;  & & \label{DAHAxx} \\
 X^\nu&=q^{2s}=e^{\hbar s},                & &\nu=[0,\hbar s]\in\wh P, & &s\in\epsilon_0^{-1}\mathbb Z; \label{DAHAxq} \\
 (T_i-t_i)&(T_i+t_i^{-1})=0; & & & & \label{DAHAtt} \\
 \underbrace{T_iT_jT_i\ldots}_{\text{$m_{ij}$ factors}}&=\underbrace{T_jT_iT_j\ldots}_{\text{$m_{ij}$ factors}},   & &\text{if $i\ne j$} & &\text{and $m_{ij}<\infty$}; \label{DAHAttt} \\
 \pi_r T_i\pi_r^{-1}&=T_j, & &\text{if $\pi_r(\alpha_i)=\alpha_j$},& &r\in{\cal O}; \label{DAHApt} \\
 T_i^{-1}X^\mu T_i^{-1}&=X^{\mu-\alpha_i}, & &\text{if $(\mu,\alpha_i^\vee)=1$},& &\mu\in\wh P; \label{DAHAtx1} \\
 T_iX^\mu&=X^\mu T_i, & &\text{if $(\mu,\alpha_i^\vee)=0$},& &\mu\in\wh P; \label{DAHAtx0} \\
 \pi_r X^\mu\pi_r^{-1}&=X^{\pi_r(\mu)}, & &r\in{\cal O},& &\mu\in\wh P; \label{DAHApx}
\end{align}}
where $i,j=0,1,\ldots,n$.
\end{Def}

Note that the generators $T_0,T_1,\ldots,T_n$, $X^{\mu}$, $\pi_r$ of the algebra $\HH_t$ are invertible elements. We define an element $T_{\hat w}\in\HH_t$ for each element $\hat w\in\wh W$ according to its reduced expression:
\begin{align} \label{Twdef}
 T_{\hat w}&=\pi_rT_{i_1}\ldots T_{i_l}, & &\text{for $\hat w=\pi_r s_{i_1}\cdots s_{i_l}$},
\end{align}
where $l=l(\hat w)$ is the length of $\hat w\in\wh W$: $l(\hat w)=|\wh R_-\cap\hat w\wh R_+|$, where $\wh R_+=-\wh R_-=R_+\cup\{[\alpha,\hbar k]\mid\alpha\in R, k\in\mathbb Z_{>0}\}$. In particular, $T_{\pi_r}=\pi_r$ and $T_{s_i}=T_i$. For $\hat w,\hat v\in\wh W$ one has the relation
\begin{align}
 T_{\hat v\hat w}&=T_{\hat v}T_{\hat w}, & &\text{if $l(\hat v\hat w)=l(\hat v)+l(\hat w)$}.
\end{align}
Since for the fundamental coweights we have $l\big(\tau(b_i+b_j)\big)=l\big(\tau(b_i)\big)+l\big(\tau(b_j)\big)$ it leads to $T_{\tau(b_i)}T_{\tau(b_j)}=T_{\tau(b_j)}T_{\tau(b_i)}$ for all $i,j=1,\ldots,n$. Let us define the commutative family of elements $Y^\lambda$ for $\lambda\in P^\vee$:
\begin{align} \label{Ydef}
 Y^\lambda&=\prod_{i=1}^n (T_{\tau(b_i)})^{m_i}, &\text{if $\lambda=\sum\limits_{i=1}^n m_i b_i$,\quad $m_i\in\mathbb Z$.}
\end{align}
They span a subalgebra $\mathbb C[Y]$ of $\HH_t$ isomorphic to $\mathbb C[P^\vee]$ and they satisfy the relations (see \cite{Kir})
\begin{align}
 Y^\lambda Y^{\lambda'}&=Y^{\lambda+\lambda'},                & &\lambda,\lambda'\in P^\vee;  \\
 T_iY^\lambda-Y^{s_i\lambda}T_i&=(t_i-t_i^{-1})\frac{Y^\lambda-Y^{s_i\lambda}}{1-Y^{-\alpha_i^\vee}}, & &\lambda\in P^\vee, \label{TiY}
\end{align}
where $i=1,\ldots,n$.  These formulae imply that the generators $T_1,\ldots,T_n$ commute with the elements of the subalgebra $\mathbb C[Y]^W\subset\mathbb C[Y]$ consisting of the polynomials invariant under the usual action $w\colon Y^\lambda\mapsto Y^{w\lambda}$ of the Weyl group $W$.

Let ${\cal F}=C^\infty(V)^{2\pi i Q^\vee}$ be the space of smooth complex-valued periodic functions $f$ on $V$ with the lattice of periods $2\pi i Q^\vee$, that is $f(x+2\pi i\alpha_j^\vee)=f(x)$ for $j=1,\ldots,n$. The group $\wh W$ acts on ${\cal F}$ as
\begin{align}
 \hat w f(x)&=f(\hat w^{-1}x), & & \hat w\in\wh W,\quad f\in{\cal F},\quad x\in V.
\end{align}

The elements of the algebra $\HH_t$ can be represented as operators acting in ${\cal F}$:
\begin{align}
 X^\mu&=e^{\mu(x)}, & \pi_r&=\pi_r, & T_i=\frac{t_i-t_i^{-1}}{1-e^{-\alpha_i(x)}}+\frac{t_i^{-1}-t_ie^{-\alpha_i(x)}}{1-e^{-\alpha_i(x)}}\,s_i, \label{ReprDAHA}
\end{align}
where $\mu\in\wh P$, $r\in{\cal O}$, $i=0,1,\ldots,n$. Rearranging the last formula into the form
\begin{align}
  T_i=t_i+\frac{t_i^{-1}-t_ie^{-\alpha_i(x)}}{1-e^{-\alpha_i(x)}}\,(s_i-1), \label{ReprDAHA_T}
\end{align}
we see that $T_i f\in{\cal F}$ for all $f\in{\cal F}$. Indeed, for any $i=0,1,\ldots,n$ the periodic function $g(x)=f(x)-f(s_ix)=f(x)-f\big(x-\alpha_i(x)(\alpha_i^{(1)})^\vee\big)$ vanishes whenever $\alpha_i(x)\in2\pi i\mathbb Z$, hence the function $\tilde g(x)=\dfrac{g(x)}{1-e^{-\alpha_i(x)}}$ is periodic and smooth: $\tilde g\in{\cal F}$. Therefore, the function
\begin{align} \label{Tifx}
 T_i f(x)=t_i f(x)+\big(t_i^{-1}-t_ie^{-\alpha_i(x)}\big)\tilde g(x)
\end{align}
also belongs to ${\cal F}$.

The formulae~\eqref{ReprDAHA} define a faithful representation of DAHA $\HH_t$ on the space $\mathcal F$ called {\it polynomial representation} (extended to smooth functions). The elements $Y^{\lambda}$ in this representation are called {\it Cherednik-Dunkl operators}.

\subsection{Macdonald-Ruijsenaars systems}
\label{sec22}

Denote by ${\cal F}^W$ the subspace of invariant functions $f\in{\cal F}$, that is $f(w x)=f(x)$ for all $w\in W$ and $x\in V$. Consider the subalgebra $\mathbb C[Y]^W\subset\mathbb C[Y]\subset\HH_t$. It is spanned as a vector space by the elements
\begin{align}
m_\lambda(Y)\stackrel{}{=}\sum_{\nu\in W\lambda}Y^{\nu}=\frac1{|W_\lambda|}\sum_{w\in W}Y^{w\lambda}, \label{flambdadef}
\end{align}
where $W\lambda=\{w\lambda\mid w\in W\}$ is the \href{http://en.wikipedia.org/wiki/Group_action#Orbits_and_stabilizers}{orbit} of the coweight $\lambda\in P^\vee$ and $W_\lambda=\{w\in W\mid w\lambda=\lambda\}$ is the \href{http://en.wikipedia.org/wiki/Group_action#Orbits_and_stabilizers}{stabilizer} of $\lambda$. We will normally be choosing $\lambda\in P^\vee_-$ for the operators $m_\lambda(Y)$.

In the representation~\eqref{ReprDAHA} the elements $m_\lambda(Y)$ have the form
\begin{align}
 m_\lambda(Y)=\sum_{\nu\in P^\vee}\sum_{w\in W}g^\lambda_{\nu,w}(x)\tau(\nu)w, \label{VIAg}
\end{align}
where $g^\lambda_{\nu,w}(x)$ are some meromorphic functions and $\tau(\nu)$ are the shift operators acting on functions by $\tau(\nu)f(x)=f(x-\hbar\nu)$.

The subspace ${\cal F}^W\subset{\cal F}$ can be characterized as a space annihilated by all the operators $T_1-t_1$, \ldots, $T_n-t_n$. Since the elements of $\mathbb C[Y]^W$ commute with $T_1,\ldots,T_n$ they preserve this subspace: $m_\lambda(Y)\big({\cal F}^W\big)\subset{\cal F}^W$. Thus we can restrict the elements of $\mathbb C[Y]^W$ to ${\cal F}^W$ \cite{Cher92DAHA,Cher92}. The restriction of~\eqref{VIAg} is
\begin{align}
 M_\lambda=\Res\big(m_\lambda(Y)\big)=m_\lambda(Y)\Big|_{{\cal F}^W}=\sum_{\nu\in P^\vee}\sum_{w\in W}g^\lambda_{\nu,w}(x)\tau(\nu).\label{Flamndadef}
\end{align}

Let $\Mer(V)$ be the space of meromorphic functions $\varphi(x)$ on $V$ such that $\varphi(x+2\pi i\alpha_j^\vee)=\varphi(x)$ for all $j=1,\ldots,n$. Define the subspace ${\cal M}\subset\Mer(V)$ consisting of functions $\varphi\in\Mer(V)$ holomorphic out of the union of hyperplanes $\bigcup_{ \beta\in{\wh R}}\{x\in V\mid \beta(x)\in2\pi i\mathbb Z \}$. Consider the \href{http://en.wikipedia.org/wiki/Group_Hopf_algebra#Hopf_module_algebras_and_the_Hopf_smash_product}{smash-product algebra} $\mathbb D={\cal M}\mathop{\#}\mathbb C[\tau(P^\vee)]$ defined by the commutation relation $\tau(\lambda)\varphi(x)=\varphi(x-\hbar\lambda)\tau(\lambda)$. This is the algebra of difference operators acting on the space ${\cal M}$. Any element $a\in\mathbb D$ has the form $a=\sum_{\lambda\in P^\vee}\sum_{w\in W}g_{\lambda,w}(x)\tau(\lambda)w$, where $g_{\lambda,w}\in\mathcal M$ and the sum over $\lambda$ is finite.
Note that $M_\lambda \in \mathbb D$ for any $\lambda \in P^\vee$.

Let $r\in{\cal O}'$, so $b_r\in P^\vee$ is a minuscule coweight. Then the restriction of $m_{-b_r}(Y)$ to ${\cal F}^W$ can be written explicitly~\cite{Cher92DAHA,Cher95,Mac2}:
\begin{align}
 M_{-b_r}=\sum_{\lambda\in W(-b_r)}\Big(\prod_{\alpha\in R \atop (\alpha,\lambda)=1}\frac{t_\alpha e^{(\alpha,x)}-t_\alpha^{-1}}{e^{(\alpha,x)}-1}\Big)\tau(\lambda)= \label{RMop_br}
\end{align}
\begin{align}
\sum_{\lambda\in W(-b_r)}\Big(\prod_{\alpha\in R \atop (\alpha,\lambda)=1}\frac{t_\alpha e^{(\alpha,x)}-t_\alpha^{-1}}{e^{(\alpha,x)}-1}\Big)\Big(\tau(\lambda)-1\Big)+m_{-b_r,t}\,, \label{RMop_br1}
\end{align}
where $m_{\lambda,t}=\sum_{\nu\in W\lambda}\prod_{\a\in R_+} t_\a^{(\a,\nu)}$.
Note that in the cases $R=E_8$, $F_4$, $G_2$ the minuscule coweights do not exist. For any $R$ there exists quasi-minuscule coweight
 $\lambda=-\theta^\vee$. The corresponding operator has the form~\cite{Cher95,Mac2}:
\begin{align}
 M_{-\theta^\vee}=\sum_{\beta\in W\theta}\Big(\prod_{\alpha\in R\atop (\alpha,\beta^\vee)>0}\frac{t_\alpha e^{(\alpha,x)}-t_\alpha^{-1}}{e^{(\alpha,x)}-1}\Big)\;\frac{q^{-1}t_0e^{(\beta,x)}-qt_0^{-1}}{q^{-1}e^{(\beta,x)}-q}\;\Big(\tau(\beta^\vee)-1\Big)+m_{-\theta^\vee,t}\,. \label{RMop_th}
\end{align}
The operators~\eqref{RMop_br}, \eqref{RMop_th} are called {\it Macdonald operators} (see e.g.~\cite{Mac88}).

Consider the subspace $\mathfrak R\subset\mathbb D$ spanned by the elements $M_\lambda$ defined by~\eqref{Flamndadef}. Since the restriction operation $\Res\colon m_\lambda(Y)\mapsto M_\lambda$ is an algebra homomorphism the subspace $\mathfrak R=\Res\big(\mathbb C[Y]^W\big)$ is a commutative subalgebra in $\mathbb D$. We will recall that $\mathfrak R$ contains $n$ algebraically independent elements.
 First we introduce some notations. For each $\lambda\in P^\vee$ there exists a unique dominant coweight $\lambda_+\in P^\vee_+$ such that $\lambda_+\in W\lambda$.
Let $\nu,\lambda\in P^\vee$ be two coweights. We write $\nu>\lambda$ if $\nu\ne\lambda$ and $\nu-\lambda=\sum_{j=1}^n k_j\alpha^\vee_j$ for some $k_j\in\mathbb Z_{\ge0}$. We also introduce another order on $P^\vee$: we write $\nu\prec\lambda$ if $\nu_+<\lambda_+$ or ($\nu_+=\lambda_+$ and $\nu>\lambda$), (see~\cite{Kir}).

As it is explained in~\cite{Kir} the Cherednik-Dunkl operators have the form
\begin{align}
 Y^\lambda=\sum_{w\in W}\tilde g^\lambda_w(x)\tau(\lambda)w+\sum_{\nu\prec\lambda}\sum_{w\in W}\tilde g^\lambda_{\nu,w}(x)\tau(\nu)w,
\end{align}
for some functions $\tilde g^\lambda_{\nu,w}(x)$, $\tilde g^\lambda_w(x)$.
For an anti-dominant coweight $\lambda\in P^\vee_-$ one has a more exact formula~\cite{Kir}
\begin{align}
 Y^\lambda=g^\lambda(x)\tau(\lambda)+\sum_{\nu\prec\lambda}\sum_{w\in W}\tilde g^\lambda_{\nu,w}(x)\tau(\nu)w.
\end{align}
Here $g^\lambda(x)=\prod_{\alpha\in\wh R_\lambda}\frac{t_\alpha e^{\alpha(x)}-t_\alpha^{-1}}{e^{\alpha(x)}-1}$, where $\wh R_\lambda=\wh R_-\cap(\tau(\lambda)\wh R_+)$ is a subset in $\wh R$.
Thus, for any $\lambda\in P^\vee_-$ we obtain
\begin{align}
 m_\lambda(Y)=g^\lambda(x)\tau(\lambda)+\sum_{\nu\prec\lambda}\sum_{w\in W} g^\lambda_{\nu,w}(x)\tau(\nu)w \label{mlambaY}
\end{align}
for some functions $\tilde g^\lambda_{\nu,w}(x)$.
In accordance with the formula~\eqref{Flamndadef} the restriction of~\eqref{mlambaY} to the subspace ${\cal F}^W$ gives the following difference operator
\begin{align} \label{MlambdaY}
 M_\lambda=g^\lambda(x)\tau(\lambda)+\sum_{\nu\prec\lambda}g^\lambda_{\nu}(x)\tau(\nu),
\end{align}
where $g^\lambda_{\nu}(x)=\sum_{w\in W} g^\lambda_{\nu,w}(x)$. Since in the case $\lambda\in P^\vee_-$ the condition $\nu\prec\lambda$ implies $\nu>\lambda$ we can replace the order $\prec$ with the order $>$ in the sum in the formula~\eqref{MlambdaY}.
It follows now that for any $m_1,\ldots,m_n\in\mathbb Z_{\ge0}$ one has
\begin{align}
 \prod_{j=1}^n\big(M_{-b_j}\big)^{m_j}=g^{m_1,\ldots,m_n}(x)\tau(\lambda)+\sum_{\nu>\lambda}g^{m_1,\ldots,m_n}_{\nu}(x)\tau(\nu), \label{Fbmonom}
\end{align}
where $g^{m_1,\ldots,m_n}(x)\not\equiv0$ and $\lambda=-\sum_{i=1}^n m_jb_j\in P^\vee_-$. Since the right hand sides of~\eqref{Fbmonom} are linearly independent the operators $M_{-b_1},\ldots,M_{-b_n}$ are algebraically independent.

\begin{Rem} \normalfont One can show that the operator~\eqref{Fbmonom} has the form $c^\lambda_\lambda M_\lambda+\sum_{\nu\in P^\vee_-, \nu>\lambda}c^\lambda_\nu M_\nu$ such that $c^\lambda_\lambda\ne0$. Then \cite[VI, \S 3, 4, Lemma 4]{Bourb} implies that the operators~\eqref{Fbmonom} form a basis of the space $\mathfrak R$. It follows the commutative algebra $\mathfrak R$ is isomorphic to the algebra of polynomials $\mathbb C[z_1,\ldots,z_n]$, and the isomorphism can be given by $M_{-b_j}\mapsto z_j$. In particular the elements $M_{-b_1},\ldots,M_{-b_n}$ generate the algebra $\mathfrak R$.
\end{Rem}

The subalgebra $\mathfrak R\subset\mathbb D$ is called the {\it Macdonald-Ruijsenaars system}. The Macdonald-Ruijsenaars operator~\eqref{RMop_br} or the operator~\eqref{RMop_th} can be considered as a Hamiltonian of this system.

\section{Invariant shifted parabolic ideals}
\label{sec3}

The goal of this section is to construct some $\HH_t$-submodules of ${\cal F}$, which will be used in the next section. Note that any ideal of the algebra ${\cal F}$ is invariant under the action of the generators $X^\mu\in\HH_t$. We define some ideals as the sets of functions vanishing on certain planes. We show that these ideals are invariant under the actions of other generators of DAHA if parameters satisfy special conditions.

Let $\eta\colon R\to\mathbb C$ be a $W$-invariant function on the irreducible reduced root system $R$, that is $\eta_{w(\alpha)}=\eta_\alpha=\eta(\alpha)$, $\forall\alpha\in R$, $\forall w\in W$. We will also denote $\eta_i=\eta_{\alpha_i}$ for the simple roots $\alpha_i$. To each subset $S\subset R$ we associate the plane
\begin{align} \label{HetaS}
 H_\eta(S)=\{x\in V \mid (\alpha,x)=\eta_\alpha,\forall\alpha\in S\} =\bigcap_{\alpha\in S}\{x\in V\mid (\alpha,x)=\eta_\alpha\}.
\end{align}
Note that $H_\eta(wS)=wH_\eta(S)$ for any $w\in W$.

Let $S_0$ be a subset of the set of simple roots $\Delta$ and let $W_0\subset W$ be the subgroup generated by $s_\alpha, \alpha\in S_0$. We associate to $S_0$ the following subset
\begin{align}
 W_+&=\{w\in W\mid wS_0\subset R_+\}\subset W. \label{DefWplus}
\end{align}
Let $V_0=\langle S_0\rangle$ be the linear subspace of $V$ spanned by the set $S_0$. Then $R_0=V_0\cap R$ is the root system with the set of simple roots $S_0$. Note that the transformations $w\in W_+$ restricted to $R_0$ respect the standard order. More exactly, if $\beta\le\gamma$, $\beta,\gamma\in R_0$ then $w\beta\le w\gamma$. (Here $\delta_1\le\delta_2$ means $\delta_2-\delta_1=\sum_{i=1}^n m_i\alpha_i$, where $m_i\in\mathbb Z_{\ge0}$).

\begin{Rem} \label{Rem31} \normalfont
 The elements of $W_+$ are the shortest coset representatives of the left cosets of $W_0$ in $W$ (see~\cite{H}).
\end{Rem}

Let us consider the plane
\begin{align} \label{HetaS0}
 \mathcal D_0&=H_\eta(S_0)=\{x\in V \mid (\alpha,x)=\eta_\alpha,\forall\alpha\in S_0\}
\end{align}
and the union of planes
\begin{align}
 {\cal D}&=\bigcup_{w\in W_+}H_\eta(wS_0)=W_+\mathcal D_0. \label{DandIdefD}
\end{align}

The arrangement $\mathcal D$ is obtained by shifting the linear arrangement $WH_0(S_0)$, where $H_0(S_0)$ is defined by the formula~\eqref{HetaS0} with the function $\eta\equiv0$. More exactly, the following two statements take place.

\begin{Prop} \label{AppRemLem1}
 For any $w\in W$ there exists $w'\in W_+$ such that $H_0(wS_0)=H_0(w'S_0)$. Thus each plane $wH_0(S_0)$ is parallel to some plane $w'{\cal D}_0 \subset \cal D$.
\end{Prop}

\noindent{\bf Proof.} The element $w\in W$ can be decomposed as $w=w'w_0$ where $w'\in W_+$, $w_0\in W_0$ (see~\cite{H}). Since $H_0(S_0)=H_0(R_0)$ and $w_0R_0=R_0$ we have $w_0H_0(S_0)=H_0(S_0)$. Hence $H_0(wS_0)=w'w_0H_0(S_0)=w'H_0(S_0)=H_0(w'S_0)$. \qed

\begin{Prop} \label{AppRemLem2} If the planes $w\mathcal D_0$ and $w'\mathcal D_0$ are parallel for some $w,w'\in W_+$ then they coincide: $w\mathcal D_0=w'\mathcal D_0$.
\end{Prop}

\noindent{\bf Proof.} Let $w,w'\in W_+$. If the planes $w\mathcal D_0$ and $w'\mathcal D_0$ are parallel then $wV_0=w'V_0$. In particular, this implies
\begin{align} \label{walpha}
 w\alpha&=\sum_{\beta\in S_0}\gamma_{\alpha\beta}w'\beta,  &&\text{for all $\alpha\in S_0$,}
\end{align}
where $\gamma_{\alpha\beta}\in\mathbb C$. Applying $(w')^{-1}$ to~\eqref{walpha} we conclude that $\gamma_{\alpha\beta}\in\mathbb Z_{\ge0}$ $\forall\beta\in S_0$ or $\gamma_{\alpha\beta}\in\mathbb Z_{\le0}$ $\forall\beta\in S_0$. Taking into account $w\alpha\in R_+$ and $w'\beta\in R_+$ we deduce from~\eqref{walpha} that $\gamma_{\alpha\beta}\in\mathbb Z_{\ge0}$ $\forall\alpha,\beta\in S_0$.
Analogously we have
\begin{align} \label{wprimebeta}
 w'\beta&=\sum_{\alpha\in S_0}\wt\gamma_{\beta\alpha}w\alpha,  &&\text{for all $\beta\in S_0$,}
\end{align}
where $\wt\gamma_{\beta\alpha}\in\mathbb Z_{\ge0}$. The formulae~\eqref{walpha}, \eqref{wprimebeta} lead to the following system of equations
\begin{align} \label{gammasystem}
\sum_{\beta\in S_0}\gamma_{\alpha\beta}\wt\gamma_{\beta\wt\alpha}=\delta_{\alpha,\wt\alpha}, &&&\alpha,\wt\alpha\in S_0,
\end{align}
where $\delta_{\alpha,\beta}$ is the delta-symbol: $\delta_{\alpha,\alpha}=1$, $\delta_{\alpha,\beta}=0$ for $\alpha\ne\beta$. Under conditions $\gamma_{\alpha\beta},\wt\gamma_{\beta\alpha}\in\mathbb Z_{\ge0}$ the general solution of the system~\eqref{gammasystem} is
\begin{align}
 \gamma_{\alpha\beta}=\wt\gamma_{\beta\alpha}=\delta_{\sigma(\alpha),\beta},
\end{align}
where $\sigma\colon S_0\to S_0$ is a bijection. Hence
\begin{align} \label{walpha_}
 w\alpha&=\sum_{\beta\in S_0}\delta_{\sigma(\alpha),\beta}w'\beta=w'\sigma(\alpha),  &&\text{for all $\alpha\in S_0$.}
\end{align}
Using this one yields
\begin{multline}
 w\mathcal D_0=\{x\in V\mid (w\alpha,x)=\eta_\alpha,\forall\alpha\in S_0\}=\{x\in V\mid (w'\sigma(\alpha),x)=\eta_\alpha,\forall\alpha\in S_0\}= \\
 =\{x\in V\mid (w'\sigma(\alpha),x)=\eta_{\sigma(\alpha)},\forall\alpha\in S_0\}=\{x\in V\mid (w'\alpha,x)=\eta_\alpha,\forall\alpha\in S_0\}=w'\mathcal D_0, \notag
\end{multline}
where we used $\eta_\alpha=\eta_{\sigma(\alpha)}$ obtained from~\eqref{walpha_}. \qed \\

Let $R_1,\ldots,R_m$ be irreducible components of $R_0=V_0\cap R$, $R_0=\bigsqcup_{\ell=1}^m R_\ell$. In each root subsystem $R_\ell$ we have the corresponding set of simple roots $S_\ell=S_0\cap R_\ell$, so that $S_0=\bigsqcup_{\ell=1}^m S_\ell$, where $S_\ell\cap S_{\ell'}=\emptyset$ for $\ell\ne\ell'$. The reflections $s_\alpha$, $\alpha\in R_\ell$, generate the Weyl group $W_\ell\subset W$ corresponding to the root system $R_\ell$. For each root system $R_\ell$ there exists the maximal root $\theta_\ell\in R_\ell$ with respect to $S_\ell$. Let  $k_\ell=|S_\ell|$ and $S_\ell=\{\alpha^\ell_1,\alpha^\ell_2,\ldots,\alpha^\ell_{k_\ell}\}$, then $\theta_\ell=\sum_{j=1}^{k_\ell}n^\ell_j\alpha^\ell_j$, where $n^\ell_j\ge1$ are integers.

\begin{Def} \label{defgCn}
Let $\theta=\sum_{i=1}^n n_i\alpha_i$ be the maximal root of $R$ with respect to the set of simple roots $\Delta$. The number $$h_\eta=\eta_\theta+\sum_{i=1}^n n_i\eta_{\alpha_i}$$ is called generalized Coxeter number for the root system $R$ and the invariant function $\eta$.
\end{Def}

We note that the generalised Coxeter number for the root system $R$ does not depend on the choice of $R_+$ (or, equivalently, on the choice of $\Delta$). For a constant function $\eta_\alpha\equiv\eta$ the generalized Coxeter number $h_\eta=\eta\cdot h$, where $h=1+\sum_{i=1}^n n_i$ is the usual Coxeter number (see~\cite{Bourb}).

\begin{Rem} \normalfont
 The generalized Coxeter number can be equivalently defined as the coefficient of proportionality of the following $W$-invariant inner products (c.f.~\cite{F}):
\begin{align}
 \sum_{\alpha\in R}\frac{\eta_\alpha(\alpha,u)(\alpha,v)}{(\alpha,\alpha)}=h_\eta\cdot(u,v), &&&\forall u,v\in V.
\end{align}
\end{Rem}

Note that the restriction of the function $\eta\colon R\to\mathbb C$ to the subsystems $R_\ell\subset R$ is a $W_\ell$-invariant function on $R_\ell$. Since $R_\ell$ is irreducible one can define the generalized Coxeter numbers for this subsystem:
\begin{align}
 h^\ell_\eta=\eta_{\theta_\ell}+\sum_{j=1}^{k_\ell} n^\ell_j\eta^\ell_j, \label{helleta}
\end{align}
where $\eta^\ell_j=\eta_{\alpha^\ell_j}$.

Define the ideal of functions vanishing on the union of planes~\eqref{DandIdefD}:
\begin{align}
 I_{\cal D}&=\{f\in{\cal F}\mid f(x)=0, \forall x\in{\cal D}\}. \label{DandIdef}
\end{align}
The main result of this section is the following theorem.

\begin{Th} \label{ThInv}
Suppose that the subsystem $S_0\subset\Delta$ and the function $\eta$ are such that
\begin{align}
 h^\ell_\eta&=\hbar & &\text{for all\quad $\ell=1,\ldots,m$}. \label{condInvGCN}
\end{align}
Let the parameters of DAHA $\HH_t$ satisfy $t_\alpha^2=e^{\eta_\alpha}$ for all $\alpha\in R$. Then the ideal $I_{\cal D}\subset{\cal F}$ is invariant under the action of DAHA: $\HH_t I_{\cal D}\subset I_{\cal D}$.
\end{Th}

We will prove this theorem in the next tree subsections. We establish the invariance under $T_i$, $1\le i\le n$ in Subsection~\ref{sec31}, the invariance under $\pi_r$ is established in Subsection~\ref{sec32} and the invariance under $T_0$ is established in Subsection~\ref{sec33}.

\begin{Rem} \normalfont
The degenerations of the ideals $I_{\cal D}$ at $\eta\to0$ give {\it invariant parabolic ideals} for the rational Cherednik algebra $H_c(W)$ with $c=\eta/\hbar$. Any radical ideal invariant under $H_c(W)$ can be obtained this way (see~\cite{F}).
\end{Rem}

\begin{Rem} \normalfont
 The invariant ideals $I_{\cal D}$ for the $A$ series are the ideals obtained by Kasatani in~\cite{KasA}, namely these are the ideals $I_m^{(k,r)}$ with $r=2$ and the same $m$ as we use. (The ideals in \cite{KasA} are considered in the space of polynomials of multiplicative variables).
\end{Rem}

Let us record the introduced planes in terms of affine functionals. For a root $\alpha\in R_+$ we denote $\alpha'=[\alpha,-\eta_\alpha]\in \wh V$. Note that $w(\alpha')=[w\alpha,-\eta_\alpha]=[w\alpha,-\eta_{w\alpha}]=(w\alpha)'$ for all $w\in W_+$, $\alpha\in R_+$. For a subset $\mathcal S\subset \wh V$ let
\begin{align}
 H(\mathcal S)=\{x\in V\mid \beta(x)=0,\forall\beta\in\mathcal S\}.
\end{align}
We will write $H\big(\beta_1,\beta_2,\ldots\big)$ instead of $H\big(\{\beta_1,\beta_2,\ldots\}\big)$. In these notations the definitions~\eqref{HetaS} and~\eqref{DandIdefD} take the form
\begin{align}
 H_\eta(S)&=H(S')=\{x\in V\mid \alpha'(x)=0,\forall\alpha\in S\}, &
 {\cal D}&=\bigcup_{w\in W_+}H(wS_0'),
\end{align}
where $S'=\{\alpha'\mid \alpha\in S\}$, $S'_0=\{\alpha'\mid \alpha\in S_0\}$. Note also that $wH(S')=H(wS')$.

\subsection{Invariance under $T_1$,\ldots, $T_n$}
\label{sec31}

\begin{Lem} \label{LemSumPhi}
 Let $S=wS_0$ where $S_0\subset\Delta$, and $w\in W_+$. Let $\beta\in\Delta$. If $\beta\not\in S$ then $\beta$ cannot be expressed as a linear combination of the elements of $S$.
\end{Lem}

\noindent{\bfseries Proof.} Suppose that there exist numbers $\gamma_\alpha\in\mathbb C$ for $\alpha\in S_0$ such that $\beta=\sum_{\alpha\in S_0}\gamma_\alpha w\alpha$. Then, $w^{-1}(\beta)=\sum_{\alpha\in S_0}\gamma_\alpha\alpha$. Since $S_0\subset\Delta$ this equation implies that the numbers $\gamma_\alpha$ are integers with the same sign: $\gamma'_\alpha\ge0$ ~$\forall\alpha\in S_0$ or $\gamma'_\alpha\le0$ ~$\forall\alpha\in S_0$. Taking into account that $w \a \in R_+$ for any $\a \in S_0$ we see that the relation $\beta=\sum_{\alpha\in S_0}\gamma_\alpha w\alpha$ holds only if $\beta\in S$, but this contradicts the condition of the lemma. \qed

\begin{Prop} \label{propInvTi}
  Let parameters of DAHA $\HH_t$ satisfy $t_\alpha^2=e^{\eta_\alpha}$ for all $\alpha\in R$. Then the ideal $I_{\cal D}$ given by~\eqref{DandIdef}, \eqref{DandIdefD}, \eqref{HetaS0} is invariant under the action of the operators $T_1,\ldots,T_n$, moreover $T_i I_{\cal D}=I_{\cal D}$ for all $i=1,\ldots,n$.
\end{Prop}

\noindent{\bfseries Proof.} Let $f\in I_{\cal D}$. Consider the subset $S=wS_0\subset R_+$, where $w\in W_+$. We need to prove that the function $T_if(x)$ vanishes on the plane $H(S')$.

Consider first the case $\alpha_i\in S$. Consider the expression~\eqref{Tifx} for the function $T_if(x)$. Let $x\in H(S')$, then the first term in this expression vanishes since the function $f(x)$ vanishes on the plane $H(S')$ and the second term vanishes since $t_i^{-1}-t_ie^{-\alpha_i(x)}=t_i^{-1}-t_ie^{-\eta_i}=0$.

Now consider the case $\alpha_i\not\in S$. Then we have $s_i(S)\subset R_+$. This means that $s_i w\in W_+$, $s_iH(S')=H(s_iwS_0')\subset{\cal D}$ and $f(s_i x)=0$ whenever $x\in H(S')$. Thus, if $x\in H(S')$ and $e^{-\alpha_i(x)}\ne1$ then both terms in
\begin{align} \label{Talphafx}
T_if(x)=\frac{t_i-t_i^{-1}}{1-e^{-\alpha_i(x)}}f(x)+\frac{t_i^{-1}-t_ie^{-\alpha_i(x)}}{1-e^{-\alpha_i(x)}}\,f(s_ix)
\end{align}
vanish: $T_i f(x)=0$. Note that for any $p\in\mathbb Z$ the intersection of the planes $H(S')$ and $H([\alpha_i,2\pi ip])$ has positive codimension in $H(S')$. Indeed, otherwise $H(S')\subset H([\alpha_i,2\pi ip])$ and this implies, in particular, that $\alpha_i$ can be obtained as a linear combination of the elements of $S$, but this is impossible by Lemma~\ref{LemSumPhi}. Thus the subset $\{x\in H(S')\mid e^{-\alpha_i(x)}\ne1\}$ is everywhere dense in the plane $H(S')$. Since the function $T_i f(x)$ is continuous it vanishes on the whole plane: $T_i f(x)=0$, $\forall x\in H(S')$.

We have proved $T_iI_{\cal D}\subset I_{\cal D}$. The relation $T_i^{-1}=T_i+(t_i^{-1}-t_i)$ implies $T_i^{-1}I_{\cal D}\subset I_{\cal D}$ and then $T_iI_{\cal D}=I_{\cal D}$. \qed

\begin{Rem} \normalfont If the $W$-invariant function $\eta\colon R\to\mathbb C$ is fixed then one can choose a $\wh W$-invariant function $t\colon\wh R\to\mathbb C$ satisfying  the conditions of Proposition~\ref{propInvTi} by the formula $t_{[\beta,\hbar k]}=e^{\eta_{\beta}/2}$ for all $\beta\in R$ and $k\in\mathbb Z$.
\end{Rem}

\subsection{Invariance under the group $\Pi$}
\label{sec32}

\begin{Lem} \label{lemOn_r} {\normalfont(c.f.~\cite{Bourb}).}
Let $\theta=\sum_{i=1}^n n_i\alpha_i$ be the maximal root of $R$ with respect to the set of simple roots $\Delta$. Then $n_r=1$ for some $r$ if and only if there exists $u_r\in W$ such that $u_r(\{-\theta\}\cup\Delta)=\{-\theta\}\cup\Delta$ and $u_r(\alpha_r)=-\theta$. Furthermore, if $u_r(\alpha_i)=\alpha_j$ then $n_i=n_j$.
\end{Lem}

\noindent{\bf Proof.} In the notations of Subsection~\ref{sec21} the statement means ${\cal O}'=\{r\mid n_r=1\}$ and it is proved in~\cite{Bourb}. Now let $r\in{\cal O}'$. Define the number $r^*\in{\cal O}'$ such that $\pi_{r^*}=\pi_r^{-1}$, then $u_{r^*}=u_r^{-1}$, which leads to $u_r(-\theta)=\alpha_{r^*}$. Taking into account that  $u_r(\alpha_r)=-\theta$, $u_r(\alpha_i)=\alpha_j$, $n_r=1$ and $u_r(\Delta\backslash\{\alpha_r\})=\Delta\backslash\{\alpha_{r^*}\}$, we rearrange $u_r(-\sum_{i=1}^n n_i\alpha_i)=\alpha_{r^*}$ as
 $$n_i\alpha_j+\sum_{k\ne j,r^*}^n\tilde n_k\alpha_k-\sum_{k=1}^n n_k\alpha_k+\alpha_{r^*}=0,$$
where $\tilde n_k\in\mathbb Z$. We obtain $n_i=n_j$ since $j \ne r^*$ by extracting the coefficient at $\alpha_j$. \qed

\begin{Prop} \label{propInvPi}
Suppose that the subsystem $S_0\subset\Delta$ and the function $\eta$ are such that $h^\ell_\eta=\hbar$ for all $\ell=1,\ldots,m$. Then the arrangement~\eqref{DandIdefD} and the ideal~\eqref{DandIdef} satisfy $\pi_r{\cal D}={\cal D}$ and $\pi_r I_{\cal D}=I_{\cal D}$ for any $r\in{\cal O}$.
\end{Prop}

\noindent{\bf Proof.} It is sufficient to establish that $\pi_r^{-1}{\cal D}\subset{\cal D}$ for any $r\in{\cal O}'$. Consider the subset $wS_0\subset R_+$, where $w\in W_+$. We will show that there exists $\tilde w\in W_+$ such that $\pi^{-1}_r H(wS'_0)=H(\tilde wS'_0)$.

First we note that
\begin{align} \label{alphabr0}
 \text{if $\alpha\in R_+$ and $(\alpha,b_r)=0$}&&&\text{then $u_r\alpha\in R_+$}.
\end{align}
Indeed since the transformation $\pi_r^{-1}$ preserves the set $\wt\Delta$ it also preserves the set $\wh R_+$. Therefore $u_r\alpha=\pi_r^{-1}\alpha\in\wh R_+\cap R=R_+$, where we used $u_r=\pi_r^{-1}\tau(b_r)$ and $u_r\alpha\in R$.

Lemma~\ref{lemOn_r} and $r\in{\cal O}'$ imply $n_r=1$. Since the roots $w\alpha_j^\ell$ and $w\theta_\ell=\sum_{j=1}^{k_\ell}n_j^\ell w\alpha_j^\ell$ are positive and coweight $b_r$ is minuscule the products $(w\theta_\ell,b_r)$ and $(w\alpha_j^\ell,b_r)$ are equal to $0$ or~$1$. So for each $\ell$ one of the following two possibilities holds:
\begin{itemize}
\item[(i)] $(w\theta_\ell,b_r)=1$, hence $(w\alpha_j^\ell,b_r)=\delta_{jp_\ell}$ and $n^\ell_{p_\ell}=1$ for some $p_\ell$ such that $1\le p_\ell\le k_\ell$;
\item[(ii)] $(w\theta_\ell,b_r)=0$, hence $(w\alpha_j^\ell,b_r)=0$ for all $j=1,\ldots,k_\ell$.
\end{itemize}
We can fix $p_\ell=1$ by renumbering the roots $\alpha^\ell_1,\ldots,\alpha^\ell_{k_\ell}$ if necessary.

Consider first the case (i): $\gamma^\ell_{rj}=\delta_{j1}$ and $n^\ell_1=1$.
Let us first rearrange the plane $H(wS'_\ell)$ as follows. We can replace in
\begin{align}
 H(wS'_\ell)&=H\Big([w\alpha^\ell_1,-\eta^\ell_1],[w\alpha^\ell_2,-\eta^\ell_2],\ldots,[w\alpha^\ell_{k_\ell},-\eta^\ell_{k_\ell}]\Big) \notag
\end{align}
the functional $[w\alpha^\ell_1,-\eta^\ell_1]$ by the linear combination $\sum_{j=1}^{k_\ell}n^\ell_j[-w\alpha^\ell_j,\eta_j^\ell]$; so taking into account $w\theta_\ell=\sum_{j=1}^{k_\ell}n^\ell_jw\alpha^\ell_j$ and the formula~\eqref{helleta} we obtain
\begin{align}
 H(wS'_\ell)
 &=H\Big([-w\theta_\ell,h^\ell_\eta-\eta_{\theta_\ell}],[w\alpha^\ell_2,-\eta^\ell_2],\ldots,[w\alpha^\ell_{k_\ell},-\eta^\ell_{k_\ell}]\Big). \label{HwSvarth}
\end{align}
This gives
\begin{align} \label{pi_rHwSl}
 \pi_r^{-1}H(wS'_\ell)=H\Big([-u_rw\theta_\ell,h^\ell_\eta-\eta_{\theta_\ell}-\hbar],[u_rw\alpha^\ell_2,-\eta^\ell_2],\ldots,[u_rw\alpha^\ell_{k_\ell},-\eta^\ell_{k_\ell}]\Big)
\end{align}
using $\pi_r^{-1}=u_r\tau(-b_r)$ and $u_r\tau(-b_r)\alpha=[u_r\alpha^{(1)},\alpha^{(0)}+\hbar(\alpha^{(1)},b_r)]$, where $\alpha=[\alpha^{(1)},\alpha^{(0)}]\in\wh V$,
It follows from~\eqref{alphabr0} that $u_rw\alpha^\ell_j\in R_+$ for $j=2,\ldots,k_\ell$.

We note that $-u_rw\theta_\ell\in R_+$ since $(u_rw\theta_\ell,u_rb_r)=1$ and $u_rb_r\in P^\vee_-$. The latter property is proved as follows. Since $\Pi$ is a group there exists $r^*\in{\cal O}'$ such that $\pi_r^{-1}=\pi_{r^*}$. Rearranging $\pi_r=\pi_{r^*}^{-1}$ as $\tau(b_r)u_r^{-1}=u_{r^*}\tau(-b_{r^*})$ we obtain $u_rb_r=-b_{r^*}\in P^\vee_-$, where we used the relation $u_r\tau(b_r)u_r^{-1}=\tau(u_rb_r)$.

Since $n^\ell_1=1$ it follows by Lemma~\ref{lemOn_r} applied to the root system $R_\ell$ that there exists an element $v_\ell\in W_\ell$ such that $v_\ell(S_\ell)=\{-\theta_\ell,\alpha^\ell_2,\ldots,\alpha^\ell_{k_\ell}\}$. So using $h^\ell_\eta=\hbar$ we obtain
\begin{align} \label{pi_rHwSl_}
 \pi_r^{-1}H(wS'_\ell)&=H(u_rwv_\ell S'_\ell) &&\text{and} && u_rwv_\ell S_\ell\subset R_+\,.
\end{align}

In the case (ii) the formulae~\eqref{pi_rHwSl_} are valid for $v_\ell=1$. Indeed $\pi_r^{-1}H(wS'_\ell)=H(u_rwS'_\ell)$ and it follows from~\eqref{alphabr0} that $u_rw\alpha_j^\ell\in R_+$ for all $j=1,\ldots,k_\ell$.

Define now $\tilde w=u_rw\prod_{\ell=1}^m v_\ell$. Since $v_\ell\in W_\ell$ acts identically on $S_{\ell'}$ for $\ell'\ne\ell$ we derive
\begin{align}
 \tilde wS_0&=\bigcup_{\ell=1}^m u_rwv_\ell S_\ell\subset R_+,
\end{align}
  hence $\tilde w\in W_+$, and
\begin{align} \label{eq325}
 \pi_r^{-1}H(wS'_0)&=\bigcap_{\ell=1}^m\pi_r^{-1}H(wS'_\ell)
  =\bigcap_{\ell=1}^mH(u_rwv_\ell S'_\ell)=H\Big(\bigcup_{\ell=1}^mu_rwv_\ell S'_\ell\Big)=H(\tilde w S'_0).
\end{align}
Thus the invariance of ${\cal D}$ and, as consequence, the invariance of $I_{\cal D}$ are proven. \qed

\subsection{Invariance under $T_0$}
\label{sec33}

If ${\cal O}\ne\{0\}$ then there exists $r\in{\cal O}'$ and the invariance of $I_{\cal D}$ under $T_0$ follows from the formula $T_0=\pi_r^{-1}T_r\pi_r$. Otherwise (${\cal O}$=\{0\} for $R$=$E_8$, $F_4$, $G_2$) we need to prove the invariance under $T_0$ separately. We will prove it in the general case. We start with two preliminary lemmas.

\begin{Lem} \label{LemSumSell}
 Let $w\in W_+$. Suppose the maximal roots for the system $R$ and the subsystem~$R_\ell$ for some $\ell$ satisfy $w\theta_\ell\ne\theta$. Then $\theta\not\in wR_\ell$.
\end{Lem}

\noindent{\bf Proof.} Suppose $\theta\in wR_\ell$ and let $\beta=w^{-1}\theta$. Then $\beta\in R_\ell$ hence $\beta\le\theta_\ell$. Since $w\in W_+$ and $\beta,\theta_\ell\in R_0$ it follows $\theta=w\beta\le w \theta_\ell$ (see the comment after the formula~\eqref{DefWplus}). Since $\theta$ is the maximal root we get $\theta=w\theta_\ell$, which is a contradiction. \qed

\begin{Lem} \label{Lem_inv_s0D}
Suppose the subsystem $S_0\subset\Delta$ and the function $\eta$ are such that $h^\ell_\eta=\hbar$ for all $\ell=1,\ldots,m$. Let $w\in W_+$ and suppose $w\theta_\ell\ne\theta$ for all $\ell=1,\ldots,m$. Then $s_0 H(wS'_0)\subset{\cal D}$.
\end{Lem}

\noindent{\bf Proof.} We will prove the existence of an element $\tilde w\in W_+$ such that $s_0 H(wS'_0)=H(\tilde wS'_0)$. First we note that Lemma~\ref{LemSumSell} implies $\theta\ne w\alpha_j^\ell$ for all $j=1,\ldots,k_\ell$, $\ell=1,\ldots,m$.
Recall that if $\alpha\in R_+$ and $\alpha\ne\theta$ then $(\alpha,\theta^\vee)$ equals $0$ or $1$ (see~\cite{Bourb}). Since $w\theta_\ell=\sum_{j=1}^{k_\ell} n^\ell_jw\alpha^\ell_j\ne\theta$ there are only two possibilities for each~$\ell$:
\begin{itemize}
\item[(i)] $(w\theta_\ell,\theta^\vee)=1$ and hence $(w\alpha^\ell_j,\theta^\vee)=\delta_{jl_\ell}$, $n^\ell_{p_\ell}=1$ for some $p_\ell$ such that $1\le p_\ell\le k_\ell$;
\item[(ii)] $(w\theta_\ell,\theta^\vee)=0$ and hence $(w\alpha^\ell_j,\theta^\vee)=0$ for all $j=1,\ldots,k_\ell$.
\end{itemize}
Here we took into account that $w\theta_\ell,w\alpha^\ell_j\in R_+$.

Consider the case (i). We can assume $p_\ell=1$. Rewriting $H(wS'_\ell)$ as~\eqref{HwSvarth} and acting to this plane by $s_0=s_\theta \tau(-\theta^\vee)$ one yields
\begin{align} \label{s0HwSl}
 s_0H(wS'_\ell)=H\Big([-s_\theta w\theta_\ell,h^\ell_\eta-\eta^\ell_{\theta_\ell}-\hbar],[s_\theta w\alpha^\ell_2,-\eta^\ell_2],\ldots,[s_\theta w\alpha^\ell_{k_\ell},-\eta^\ell_{k_\ell}]\Big).
\end{align}
Taking into account that $s_\theta(\alpha)=\alpha-(\alpha,\theta^\vee)\theta$, $\forall\alpha\in R$, we conclude $s_\theta w\alpha^\ell_j=w\alpha^\ell_j\in R_+$ for $j=2,\ldots,k_\ell$ and $-s_\theta w\theta_\ell=\theta-w\theta_\ell\in R_+$.
Since $n^\ell_1=1$ there exists $v_\ell\in W_\ell$ such that $v_\ell(S_\ell)=\{-\theta_\ell,\alpha^\ell_2,\ldots,\alpha^\ell_{k_\ell}\}$ (due to Lemma~\ref{lemOn_r}). Then the condition $h^\ell_\eta=\hbar$ implies
\begin{align} \label{s0HwSl_}
 s_0H(wS'_\ell)&=H(s_\theta wv_\ell S'_\ell), & s_\theta wv_\ell S_\ell\subset R_+\,.
\end{align}

In the case (ii) the formulae~\eqref{s0HwSl_} are valid for $v_\ell=1$, because in this case $s_0w\alpha^\ell_j=s_\theta w\alpha^\ell_j=w\alpha^\ell_j\in R_+$ for all $j=1,\ldots,k_\ell$.

Thus we have $s_0H(wS'_0)=H(\tilde wS'_0)$, where $\tilde w=s_\theta w\prod_{\ell=1}^m v_\ell\in W_+$. \qed \\

Now the invariance of the ideal $I_{\cal D}$ under the action of the operator $T_0$ can be established.

\begin{Prop} \label{propInvT0}
Suppose the subsystem $S_0\subset\Delta$ and the functions $\eta$, $t$ satisfy $h^\ell_\eta=\hbar$ for all $\ell=1,\ldots,m$ and $t_0^2=e^{\eta_\theta}$. Then $T_0 I_{\cal D}=I_{\cal D}$.
\end{Prop}

\noindent{\bfseries Proof.} We need to prove that for any $f\in I_{\cal D}$ and $w\in W_+$ the function $T_0f(x)$ vanishes on the plane $H(S')$, where $S=wS_0$.

Let us first suppose that $w\theta_\ell=\theta$ for some $\ell$. Let $x\in H(S')$. Since $h^\ell_\eta=\hbar$ and $\eta_{\theta_\ell}=\eta_\theta$ one yields $\alpha_0(x)=-(w\theta_\ell,x)+\hbar=\eta_{\theta_\ell}-h^\ell_\eta+\hbar=\eta_\theta$ and hence $e^{-\alpha_0(x)}=t_0^{-2}$. Therefore the right hand side in~\eqref{Tifx} with $i=0$ vanishes.

Otherwise $w\theta_\ell\ne\theta$ for all $\ell$ and by virtue of Lemma~\ref{Lem_inv_s0D} we have $s_0H(S')\subset{\cal D}$. Hence the function
\begin{align} \label{T0fx}
T_0f(x)=\frac{t_0-t_0^{-1}}{1-e^{-\alpha_0(x)}}f(x)+\frac{t_0^{-1}-t_0e^{-\alpha_0(x)}}{1-e^{-\alpha_0(x)}}\,f(s_0x)
\end{align}
vanishes on the subset $H(S')\cap\{x\in V\mid e^{-\alpha_0(x)}\ne1\}$. Let us prove that this subset is everywhere dense in $H(S')$. Indeed, otherwise $H(S')$ is contained in a plane $H(\alpha_0+[0,2\pi ip])$ for some $p\in\mathbb Z$. This leads to $[-\theta,\hbar+2\pi ip]=\sum_{\alpha\in wS_0}\gamma_\alpha \alpha'$. Therefore $\theta\in wV_0$, where $V_0\subset V$ is the subspace spanned by $S_0$. Since $w^{-1}\theta\in R$ then $w^{-1}\theta\in R_0=V_0\cap R$, and the formula $R_0=\bigsqcup_{\ell=1}^m R_\ell$ implies that $w^{-1}\theta\in R_\ell$ for some $\ell$. But this is impossible due to Lemma~\ref{LemSumSell}. Thus the continuous function $T_0f(x)$ vanishes on all the plane $H(S')$. \qed

Theorem~\ref{ThInv} follows from Propositions~\ref{propInvTi},~\ref{propInvPi},~\ref{propInvT0}.

\subsection{Further invariant ideals}
\label{further-ideals}

In this subsection we generalize previously constructed ideals $I_{\mathcal D}$ by requiring functions to vanish on the families of parallel planes containing the planes $w\mathcal D_0$ with $w\in W_+$. For DAHA associated with the classical root systems such submodules were considered by Kasatani in~\cite{KasA,KasCC}.

Consider a subset $wS_0\subset R$, where $w\in W$, and its decomposition $wS_0=\bigsqcup_{\ell=1}^m wS_\ell$. For a function $\xi\colon w S_0\to\mathbb Z_{\ge0}$ we denote $|\xi|^w_\ell=\sum_{\alpha\in wS_\ell}\xi(\alpha)$, where $\ell=1,\ldots,m$.
Note that if $wS_0=w'S_0$ for some $w'\in W$ then the values $|\xi|^w_1,\ldots,|\xi|^w_m$ and $|\xi|^{w'}_1,\ldots,|\xi|^{w'}_m$ are related by $|\xi|^{w'}_\ell=|\xi|^w_{\sigma(\ell)}$, where a permutation $\sigma\in \mathfrak{S}_m$ is such that $R_{\sigma(\ell)}\cong R_\ell$.

Let $r_1,\ldots,r_m\in\mathbb Z_{\ge0}$ be such that $r_\ell=r_{\ell'}$ if $R_\ell\cong R_{\ell'}$. For an element $w\in W$ we define the following set of functions
\begin{align}
 \Xi_{\vec r}(w)=\{\xi\colon wS_0\to\mathbb Z_{\ge0}\mid\xi^{-1}(0)\subset R_+,\; |\xi|^w_\ell\le r_\ell, \ell=1,\ldots,m\},
\end{align}
where $\vec r=(r_1,\ldots,r_m)$. The condition $\xi^{-1}(0)\subset R_+$  means that $\xi(\alpha)=0$ implies $\alpha\in R_+$.

Let $\eta\colon R\to\mathbb C$ be a fixed $W$-invariant function. To a subset $S\subset R$ and a function $\xi\colon S\to\mathbb Z_{\ge0}$ we associate the plane
\begin{align} \label{HetaSxi}
 H_\eta(S;\xi)=\{x\in V \mid (\alpha,x)=\eta_\alpha-\hbar\xi(\alpha),\forall\alpha\in S\}.
\end{align}
Starting from the subset $S_0$ and $\vec r$ we define the union of planes
\begin{align}
 {\cal D}_{\vec r}&=\bigcup_{w\in W}\bigcup_{\xi\in\Xi_{\vec r}(w)}H_\eta(wS_0;\xi). \label{Drdef}
\end{align}
We will assume that $\vec r$ is such that
\begin{align}
 GCD(&r_\ell+1,k_\ell+1)=1, &&\text{if $R_\ell\cong A_{k_\ell}$,} && \ell=1,\ldots,m, \label{rConditA} \\
 &r_\ell=0, &&\text{if $R_\ell\not\cong A_{k_\ell}$,} && \ell=1,\ldots,m, \label{rCondit}
\end{align}
where $GCD$ stands for the greatest common divisor.
Note that in the case $\vec r=0$ the union $\mathcal D_{\vec r}$ coincides with $\mathcal D=W_+\mathcal D_0$ defined by the formula~\eqref{DandIdefD}.
We also generalize the conditions $h^\ell_\eta=\hbar$ as
\begin{align}
 h^\ell_\eta&=(r_\ell+1)\hbar. \label{gCnr}
\end{align}

\begin{Th} \label{LemFII}
Suppose the parameters of DAHA $\HH_t$ satisfy $t_\alpha^2=e^{\eta_\alpha}$ for all $\alpha\in R$. Suppose the subsystem $S_0\subset\Delta$, the function $\eta$ and the parameters $r_1,\ldots,r_m$ are such that the conditions~\eqref{rConditA}, \eqref{rCondit} and \eqref{gCnr} are satisfied. Then the ideal
\begin{align} \label{IDr}
 I_{\mathcal D_{\vec r}}=\{f\in\mathcal F\mid f(x)=0, \text{ for all $x\in\mathcal D_{\vec r}$}\}
\end{align}
where $\mathcal D_{\vec r}$ is given by~\eqref{Drdef}, is invariant under the action of DAHA $\HH_t$.
\end{Th}

First we establish the invariance under $T_1,\ldots,T_n$.

\begin{Lem} \label{LemFIITi} Under assumptions of Theorem~\ref{LemFII} one has $T_iI_{\mathcal D_{\vec r}}\subset I_{\mathcal D_{\vec r}}$ for $i=1,\ldots,n$.
\end{Lem}

\noindent{\bf Proof.} Similar to the proof of Proposition~\ref{propInvTi} we show that the function $(T_if)(x)$ vanishes on the plane $H_\eta(wS_0;\xi)$, where $f\in I_{\mathcal D_{\vec r}}$, $w\in W$ and $\xi\in\Xi_{\vec r}(w)$. If $\alpha_i\in\xi^{-1}(0)$ then it follows vanishing of the function $t_i^{-1}-t_ie^{-\alpha_i(x)}$ (see the formula~\eqref{Tifx}). Suppose now that $\alpha_i\notin\xi^{-1}(0)$. It follows $s_i(\xi^{-1}(0))=(\xi\circ s_i)^{-1}(0)\subset R_+$ and, therefore, $\xi\circ s_i\in\Xi_{\vec r}(s_iw)$. If $x\in H_\eta(wS_0;\xi)$ then $s_ix\in H_\eta(s_iwS_0;\xi\circ s_i)\subset\mathcal D_{\vec r}$ and the function $f(s_ix)$ vanishes as well as the function $f(x)$. We are left to prove that $H_\eta(wS_0;\xi)\not\subset H([\alpha_i,2\pi ip])$ for all $p\in\mathbb Z$.
 Suppose the opposite, that is $\alpha_i=\sum_{\alpha\in S_0}\gamma_\alpha w\alpha$ and $2\pi i p=\sum_{\alpha\in S_0}\gamma_\alpha(\hbar\xi(w\alpha)-\eta_\alpha)$ for some $\gamma_\alpha \in \mathbb{C}$ and $p\in\mathbb Z$. Since $w^{-1}\alpha_i=\sum_{\alpha\in S_0}\gamma_\alpha \alpha\in R_0$, $\gamma_\alpha$ are integers of the same sign and there exists $\ell$, $1 \le \ell \le m$ such that $w^{-1}\alpha_i\in R_\ell$. Therefore $2\pi i p=\sum_{\alpha\in S_\ell}\gamma_\alpha(\hbar\xi(w\alpha)-\eta_\alpha)$. If $r_\ell=0$ then $\xi\big|_{wS_\ell}\equiv0$ hence $w S_\ell \subset R_+$ and $\a_i \in w S_\ell$ which is not possible  as $\xi(\a_i)\ne 0$.
If $r_\ell\ne0$ then $R_\ell\cong A_{k_\ell}$, hence $\eta\big|_{R_\ell}=const$ and $\theta_\ell=\sum_{\alpha\in S_\ell}\alpha$. Note that $\pm w^{-1}\alpha_i=\pm\sum_{\alpha\in S_\ell}\alpha\le\theta_\ell$ hence $|\gamma_\alpha|\le1$.  The condition~\eqref{gCnr} implies
\begin{align} \label{2piip_}
2\pi i p=\hbar\Big(\sum_{\alpha\in S_\ell}\gamma_\alpha\xi(w\alpha)-\frac{r_\ell+1}{k_\ell+1}\sum_{\alpha\in S_\ell}\gamma_\alpha\Big).
\end{align}
Since $\hbar\notin2\pi i\mathbb Q$ the expression in the brackets vanishes. As $|\sum_{\alpha\in S_\ell}\gamma_\alpha|=\sum_{\alpha\in S_\ell}|\gamma_\alpha|\le k_\ell$ we get a contradiction with the condition~\eqref{rConditA}. \qed

We will need the following technical lemma for the rest of the proof of Theorem~\ref{LemFII}.

\begin{Lem} \label{ulambda}
 Let $u\in W$, $\lambda\in P^\vee_+$ such that $u\lambda\in P^\vee_-$. Suppose that $\hat u R_+\subset\wh R_+$, where $\hat u=u\tau(-\lambda)$. Let $R_\ell\cong A_{k_\ell}$ for some $\ell$ and let $S_\ell=\{\alpha_1^\ell,\ldots,\alpha_{k_\ell}^\ell\}$ be a basis of simple roots. Assume that $w\in W$ and $\xi\in\Xi_{\vec r}(w)$ satisfy the conditions $|\xi|^w_\ell+m_0^\ell\le r_\ell+1$ and $\xi_j^\ell+m_j^\ell\ge0$ for all $j=1,\ldots,k_\ell$, where $m_0^\ell=(w\theta_\ell,\lambda)$, $\xi_j^\ell=\xi(w\alpha_j^\ell)$, $m_j^\ell=(w\alpha_j^\ell,\lambda)$. Denote $\xi_\ell=\xi\big|_{wS_\ell}$. Then $\hat u H_\eta(wS_\ell;\xi_\ell)=H_\eta(uwv_\ell S_\ell;\tilde\xi_\ell)$ for some element $v_\ell\in W_\ell$ and some map $\tilde\xi_\ell\colon uwv_\ell S_\ell\to\mathbb Z_{\ge0}$ such that $\sum_{\alpha\in uwv_\ell S_\ell}\tilde\xi(\alpha)\le r_\ell$ and $\tilde\xi_\ell^{-1}(0)\subset R_+$.
\end{Lem}

\noindent{\bf Proof.} Since $R_\ell\cong A_{k_\ell}$ the maximal root in $R_\ell$ is $\theta_\ell=\sum_{j=1}^{k_\ell}\alpha_j^\ell$ and hence $m_0^\ell=\sum_{j=1}^{k_\ell}m_j^\ell$. Denoting $\zeta_j^\ell=\xi_j^\ell+m_j^\ell$ we have $\zeta_j^\ell\ge0$.  Let us show that
\begin{align}
 &\zeta_j^\ell=0 &&\text{implies} &&uw\alpha_j^\ell\in R_+. \label{zeta0implies}
\end{align}
If $\zeta_j^\ell=0$ there are two cases: $\xi_j^\ell=m_j^\ell=0$ and $\xi_j^\ell=-m_j^\ell>0$. Consider the first case. Note that $m_j^\ell=0$ implies that $uw\alpha_j^\ell=\hat uw\alpha_j^\ell$. Since $\xi_j^\ell=0$ we have $w\alpha_j^\ell\in R_+$ and using $\hat u R_+\subset\wh R_+$ we obtain $uw\alpha_j^\ell\in R\cap\wh R_+=R_+$. In the case when $-m_j^\ell>0$ we have $(uw\alpha_j^\ell,u\lambda)=m_j^\ell<0$. Since $u\lambda\in P^\vee_-$ we also deduce $uw\alpha_j^\ell\in R_+$.

There are two possible cases:
\begin{itemize}
\item[(i)] $|\xi|^w_\ell+m^\ell_0=r_\ell+1$ and hence $m^\ell_0\ge1$;
\item[(ii)] $|\xi|^w_\ell+m^\ell_0\le r_\ell$.
\end{itemize}

Consider first the case (i). Since $m^\ell_0\ge1$ there exists $p_\ell$, $1\le p_\ell\le k_\ell$ such that $m_{p_\ell}^\ell\ge1$. We can assume $p_\ell=1$, that is $m_1^\ell\ge1$. The formula~\eqref{HwSvarth} can be generalized in this case as
\begin{align}
 H_\eta(wS_\ell;\xi_\ell)= H\Big([-w\theta_\ell,-\hbar|\xi|^w_\ell+h^\ell_\eta-\eta_{\theta_\ell}],[w\alpha^\ell_2,\hbar\xi_2^\ell-\eta^\ell_2],\ldots,[w\alpha^\ell_{k_\ell},\hbar\xi_{k_\ell}^\ell-\eta^\ell_{k_\ell}]\Big).
\end{align}
Acting by $\hat u=u\tau(-\lambda)$ and taking into account~\eqref{gCnr} we obtain
\begin{align} \label{uHwSlxi}
 \hat u H_\eta(wS_\ell;\xi_\ell)= H\Big([-uw\theta_\ell,-\eta_{\theta_\ell}],[uw\alpha^\ell_2,\hbar\zeta_2^\ell-\eta^\ell_2],\ldots,[uw\alpha^\ell_{k_\ell},\hbar\zeta_{k_\ell}^\ell-\eta^\ell_{k_\ell}]\Big).
\end{align}
Due to Lemma~\ref{lemOn_r} there exists $v_\ell\in W_\ell$ such that $\{-\theta_\ell,\alpha_2^\ell,\ldots,\alpha_{k_\ell}^\ell\}=v_\ell S_\ell$. Define the function $\tilde\xi_\ell\colon uwv_\ell S_\ell\to\mathbb Z_{\ge0}$ by $\tilde\xi_\ell(-uw\theta_\ell)=0$ and $\tilde\xi_\ell(uw\alpha_j^\ell)=\zeta_j^\ell$ where $j=2,\ldots,k_\ell$. Let us check that the sum of values of the function $\tilde\xi_\ell$ is not greater than $r_\ell$. Indeed, $\sum_{j=2}^{k_\ell}\zeta_j^\ell=|\xi|^w_\ell+m_0^\ell-\xi_1^\ell-m_1^\ell = r_\ell+1-\xi_1^\ell-m_1^\ell\le r_\ell$, where we used $m_1^\ell\ge1$. Since $u\lambda\in P^\vee_-$ and $(uw\theta_\ell,u\lambda)=m_0^\ell\ge1$ we have $-uw\theta_\ell\in R_+$. Taking into account~\eqref{zeta0implies} we derive $\tilde\xi_\ell^{-1}(0)\subset R_+$.

In the case (ii) we set $v_\ell=1$ and $\tilde\xi_\ell(uw\alpha_j^\ell)=\zeta_j^\ell$ where $j=1,\ldots,k_\ell$. The sum of values of the function $\tilde\xi_\ell$ is still not grater than $r_\ell$, indeed $\sum_{j=1}^{k_\ell}\zeta_j^\ell=|\xi|^w_\ell+m_0^\ell\le r_\ell$. From~\eqref{zeta0implies} we derive $\tilde\xi_\ell^{-1}(0)\subset R_+$ for this case as well.

In both cases we arrive at $\hat u H_\eta(wS_\ell;\xi_\ell)=H_\eta(uwv_\ell S_\ell;\tilde\xi_\ell)$. \qed

Now we are ready to establish the invariance of $I_{\mathcal D_{\vec r}}$ under the operators $\pi_r$ and $T_0$.

\begin{Lem} \label{LemFIIpi}
 Let $r\in{\cal O}$. Under assumptions of Theorem~\ref{LemFII} one has $\pi_r I_{\mathcal D_{\vec r}}\subset I_{\mathcal D_{\vec r}}$.
\end{Lem}

\noindent{\bf Proof.} It is sufficient to establish $\pi_r^{-1}{\cal D}_{\vec r}\subset{\cal D}_{\vec r}$ for arbitrary $r\in{\cal O}'$. Consider an element $w\in W$ and a function $\xi\in\Xi_{\vec r}(w)$. We will show that there exists an element $\tilde w\in W$ and a function $\xi\in\Xi_{\vec r}(\tilde w)$ such that $\pi^{-1}_r H_\eta(wS_0;\xi)=H_\eta(\tilde wS_0;\tilde\xi)$. Note that we have
\begin{align}
 \pi_r^{-1}H_\eta(wS_0;\xi)=\bigcap_{\ell=1}^m\pi_r^{-1}H_\eta(wS_\ell;\xi_\ell),
\end{align}
where $\xi_\ell=\xi\big|_{wS_\ell}$. For each $\ell$ we have two possible cases: $R_\ell\cong A_{k_\ell}$ or $R_\ell\not\cong A_{k_\ell}$.

Suppose first that $R_\ell\cong A_{k_\ell}$. For $u=u_r$ and $\lambda=b_r$ we can apply Lemma~\ref{ulambda}. Indeed, it easy to see that $u_rb_r\in P^\vee_-$ (see the proof of Proposition~\ref{propInvPi}). The transformation $\hat u=u_r\tau(-b_r)=\pi_r^{-1}$ preserves the set $\wt\Delta$ and hence $\hat u R_+\subset\wh R_+$. Since $b_r$ is a minuscule coweight $m_i^\ell\in\{-1,0,1\}$ for $i=0,1,\ldots,k_\ell$ (we use the notations of Lemma~\ref{ulambda}) and hence $|\xi|^w_\ell+m_0^\ell\le r_\ell+1$. Further we have $\xi_j^\ell+m_j^\ell\ge0$ because otherwise $\xi_j^\ell=0$ and $m_j^\ell=-1$; since $\xi^{-1}(0)\subset R_+$ the first equality means that $w\alpha_j^\ell\in R_+$, but this contradicts the second equality.

So due to Lemma~\ref{ulambda} there exists an element $v_\ell\in W_\ell$ and a map $\tilde\xi_\ell\colon u_rwv_\ell S_\ell\to\mathbb Z_{\ge0}$ with the sum of values not greater than $r_\ell$ such that $\tilde\xi_\ell^{-1}(0)\subset R_+$ and
\begin{align} \label{pi_rHwSl_xi}
 \pi_r^{-1}H_\eta(wS_\ell;\xi_\ell)&=H_\eta(u_rwv_\ell S_\ell;\tilde\xi_\ell).
\end{align}

In the case $R_\ell\not\cong A_{k_\ell}$ we have $r_\ell=0$ and therefore $\xi_\ell\equiv0$,
$H_\eta(wS_\ell;\xi_\ell)=H(wS'_\ell)$. So the relation~\eqref{pi_rHwSl_} holds for an appropriate $v_\ell\in W_\ell$ (see the proof of Proposition~\ref{propInvPi}). In this case the equality~\eqref{pi_rHwSl_xi} is valid for the function  $\tilde\xi_\ell\colon u_rwv_\ell S_\ell\to\mathbb Z_{\ge0}$ defined as $\tilde\xi_\ell\equiv0$. Note that we have $\tilde\xi_\ell^{-1}(0)=u_rwv_\ell S_\ell\subset R_+$.

Define the element $\tilde w=u_rw\prod_{\ell=1}^m v_\ell$ and the function $\tilde\xi\colon\tilde wS_0\to\mathbb Z_{\ge0}$ by the formula $\tilde\xi\big|_{\tilde wS_\ell}=\tilde\xi_\ell$. We already checked in fact that $|\tilde\xi|^{\tilde w}_\ell\le r_\ell$. Since $\tilde\xi_\ell^{-1}(0)\subset R_+$ for all $\ell=1,\ldots,m$ one has $\tilde\xi^{-1}(0)\subset R_+$. So $\tilde\xi\in\Xi_{\vec r}(\tilde w)$. Thus we obtain $\pi_r^{-1}H_\eta(wS_0;\xi)=H_\eta(\tilde wS_0;\tilde\xi)$ as required (c.f.~\eqref{eq325}). \qed

\begin{Lem} \label{LemFIIT0}
Under assumptions of Theorem~\ref{LemFII} one has
 $T_0I_{\mathcal D_{\vec r}}\subset I_{\mathcal D_{\vec r}}$.
\end{Lem}

\noindent{\bf Proof.} Let $f\in I_{\mathcal D_{\vec r}}$, $w\in W$ and $\xi\in\Xi_{\vec r}(w)$. Then we need to prove that the function $(T_0f)(x)$ vanishes on the plane $H_\eta(wS_0;\xi)$. Let $\xi_\ell$, $\alpha_j^\ell$ and $\xi_j^\ell$ be defined as in Lemma~\ref{ulambda}. Then there are three different possibilities:
\begin{itemize}
\item[(a)] $\theta=w\theta_\ell$ and $|\xi|^w_\ell=r_\ell$ for some $\ell$;
\item[(b)] $w\alpha^\ell_j=-\theta$ and $\xi_j^\ell=1$ for some $j$ and $\ell$;
\item[(c)] the other cases satisfying neither the condition (a) nor (b).
\end{itemize}

In the cases (a) and (b) we have $\alpha_0(x)=\eta_\theta$ for $x\in H_\eta(wS_0;\xi)$. So that the expression~\eqref{Tifx} with $i=0$ for the function $(T_0f)(x)$ vanish on the plane $H_\eta(wS_0;\xi)$.

Now consider the case (c). We will prove in this case $s_0H_\eta(wS_0;\xi)\subset\mathcal D_{\vec r}$ by providing $\tilde w\in W$ and $\tilde\xi\in\Xi_{\vec r}(\tilde w)$ such that $s_0H_\eta(wS_0;\xi)=H_\eta(\tilde wS_0;\tilde\xi)$.

Let $R_\ell\cong A_{k_\ell}$. Put $m_0^\ell=(w\theta_\ell,\theta^\vee)$ and $m_j^\ell=(w\alpha_j^\ell,\theta^\vee)$, then $m_0^\ell=\sum_{j=1}^{k_\ell}m_j^\ell$. We have $|m_i^\ell|\le2$ for all $i=0,1,\ldots,k_\ell$. Firstly note that $|\xi|^w_\ell+m_0^\ell\le r_\ell+1$. Indeed, in the case (c) we have $\theta\ne w\theta_\ell$ or $|\xi|^w_\ell\le r_\ell-1$. If the first condition holds then the inequality we are proving follows from $m_0^\ell\le1$ and $|\xi|^w_\ell\le r_\ell$, otherwise it follows from $m_0^\ell\le2$ and $|\xi|^w_\ell\le r_\ell-1$.

 Secondly $\xi_j^\ell+m_j^\ell\ge0$ for all $j=1,\ldots,k_\ell$. Indeed, if $\xi_j^\ell+m_j^\ell<0$ then $m_j^\ell<0$, which implies $w\alpha_j^\ell\in R_-$ and hence $\xi_j^\ell\ne0$; so that the only possibility is $\xi_j^\ell=1$ and $m_j^\ell=-2$. This implies $w\alpha_j^\ell=-\theta$, which is impossible in the case (c).

Note that $s_0=s_\theta\tau(-\theta^\vee)$
satisfies $s_0 R_+\subset\wh R_+$. Thus using Lemma~\ref{ulambda} for $u=s_\theta$ and $\lambda=\theta^\vee$ we have
\begin{align} \label{s0HwSl_xi}
 s_0H_\eta(wS_\ell;\xi_\ell)&=H_\eta(s_\theta wv_\ell S_\ell;\tilde\xi_\ell)
\end{align}
for some element $v_\ell\in W_\ell$ and some map $\tilde\xi_\ell\colon s_\theta wv_\ell S_\ell\to\mathbb Z_{\ge0}$ with sum of values not greater than $r_\ell$ such that $\tilde\xi_\ell^{-1}(0)\subset R_+$.

If $R_\ell\not\cong A_{k_\ell}$ then $r_\ell=0$ and, therefore, $\xi_\ell\equiv0$,
$H_\eta(wS_\ell;\xi_\ell)=H(wS'_\ell)$. So we have the equality~\eqref{s0HwSl_} for the appropriate $v_\ell\in W_\ell$ (see the proof of Lemma~\ref{Lem_inv_s0D}). Thus in this case the relation~\eqref{s0HwSl_xi} is valid for the map $\tilde\xi_\ell\colon s_\theta wv_\ell S_\ell\to\mathbb Z_{\ge0}$ defined as $\tilde\xi_\ell\equiv0$. Note that we have $\tilde\xi_\ell^{-1}(0)=s_\theta wv_\ell S_\ell\subset R_+$.

Define the element $\tilde w=s_\theta w\prod_{\ell=1}^m v_\ell$ and the function $\tilde\xi\colon\tilde wS_0\to\mathbb Z_{\ge0}$ by the formula $\tilde\xi\big|_{\tilde wS_\ell}=\tilde\xi_\ell$. We already checked that $|\tilde\xi|^{\tilde w}_\ell\le r_\ell$. Since $\tilde\xi_\ell^{-1}(0)\subset R_+$ for all $\ell=1,\ldots,m$ one has $\tilde\xi(0)\in R_+$, so $\tilde\xi\in\Xi_{\vec r}(\tilde w)$. Then we obtain $s_0H_\eta(wS_0;\xi)=H_\eta(\tilde wS_0;\tilde\xi)$.

Thus $f(x)=f(s_0x)=0$ for all $x\in H_\eta(wS_0;\xi)$ and we are left to prove that the denominator $1-e^{-\alpha_0(x)}$ in the formula~\eqref{T0fx} does not vanish identically on the plane $H_\eta(wS_0;\xi)$. Indeed otherwise $\alpha_0+[0,2\pi i p]=\sum_{\alpha\in S_0}\gamma_\alpha[w\alpha,\hbar\xi(w\alpha)-\eta_\alpha]$ for some $p\in\mathbb Z$ and $\gamma_\alpha\in\mathbb Z$. Then there exists $\ell$, $1\le \ell \le m$ such that $\gamma_\alpha\ne0$ only for $\alpha\in S_\ell$. Hence $2\pi i p=\hbar(\sum_{\alpha\in S_\ell}\gamma_\alpha\xi(w\alpha)-1)-\sum_{\a\in S_\ell}\gamma_\a \eta_\alpha$ and $\gamma_\alpha$ are all non-negative or all non-positive. If $r_\ell=0$ then $\xi\big|_{wS_\ell}\equiv0$
and hence $w S_\ell \subset R_+$. Since $\theta\in wR_\ell$ it follows that $\theta=w\theta_\ell$ (c.f.~Lemma~\ref{LemSumSell}) which is a contradiction as we are in case (c) and not in case (a).
If $r_\ell\ne0$ then $R_\ell\cong A_{k_\ell}$, $\eta\big|_{R_\ell}=const$. The condition~\eqref{gCnr} implies
\begin{align} \label{2piip_1}
2\pi i p=\hbar\Big(\sum_{\alpha\in S_\ell}\gamma_\alpha\xi(w\alpha)-1-\frac{r_\ell+1}{k_\ell+1}\sum_{\alpha\in S_\ell}\gamma_\alpha\Big)
\end{align}
with $|\gamma_\alpha|\le1$. Since $\hbar\notin2\pi i\mathbb Q$ the expression in the brackets vanishes and due to $|\sum_{\alpha\in S_\ell}\gamma_\alpha|=\sum_{\alpha\in S_\ell}|\gamma_\alpha|\le k_\ell$ we have a contradiction with the condition~\eqref{rConditA}. \qed \\

\noindent Theorem~\ref{LemFII} follows from Lemmas~\ref{LemFIITi},~\ref{LemFIIpi} and~\ref{LemFIIT0}.

\begin{Ex} \normalfont
 Let $R=A_n$, then $R_\ell\cong A_{k_\ell}$, $m-1+\sum_{\ell=1}^m k_\ell\le n$ and $\eta_\alpha=\eta_1$ for all $\alpha\in R$. In this case the conditions~\eqref{rConditA} and \eqref{gCnr} imply $r_1=\ldots=r_m$ and $k_1=\ldots=k_m$ such that $GCD(r_1+1,k_1+1)=1$ and $(k_1+1)\eta_1=(r_1+1)\hbar$. The corresponding invariant ideals (in the representation of DAHA of type $GL_{n+1}$) were constructed in~\cite{KasA}.
\end{Ex}

\section{Generalized Macdonald-Ruijsenaars systems}
\label{sec4}

In this section we use the invariant ideals $I_{\cal D}$ found in Section~\ref{sec3} to obtain commutative algebras of difference operators acting in the space of functions defined on the affine planes. We work in the assumptions of Theorem~\ref{ThInv}, so that the ideal $I_{\mathcal D}$ is invariant and the quotient module is defined. The symmetrized Cherednik-Dunkl operators $m_\lambda(Y)$ considered in the quotient representation define a commutative algebra of difference operators. We consider explicit examples of operators arising this way.

\subsection{The operation of restriction}
\label{sec41}

In Section~\ref{sec2} we reviewed how the Macdonald-Ruijsenaars systems are obtained in terms of the polynomial representation of DAHAs. In this section we show how the generalized Macdonald-Ruijsenaars systems are obtained in terms of the representation
\begin{align} \label{rho}
 \rho_{\mathcal D}\colon\HH_t\to\End({\cal F}/I_{\mathcal D})
\end{align}
on the quotient module ${\cal F}/I_{\cal D}$. Thus $\rho_{\mathcal D}(a)(g)=a(f)\big|_{\mathcal D}$, where $a\in\HH_t$, $f\in{\cal F}$ and $g=[f]=f\big|_{\mathcal D}$ is a class of the function $f$ identified with its restriction on $\mathcal D$.

Consider the subspace ${\cal F}_0\subset{\cal F}/I_{\mathcal D}$ of functions on $\mathcal D$ which have $W$-invariant extension in $V$, that is ${\cal F}_0=\{[f]\in{\cal F}/I_{\cal D}\mid f\in{\cal F}^W\}$. Since $\mathcal D\subset W\mathcal D_0$ the elements of ${\cal F}_0$ are determined by their values on the plane ${\cal D}_0$ and we can identify the space $\mathcal F_0$ with the space $\{f\big|_{\mathcal D_0}\mid f\in\mathcal F^W\}$ of functions on $\mathcal D_0$.

Since $\mathbb C[Y]^W\mathcal F^W\subset\mathcal F^W$ the operators belonging to the algebra $\rho_{\mathcal D}\big(\mathbb C[Y]^W\big)$ preserve the space ${\cal F}_0$.
For any $a \in \mathbb C[Y]^W$ let
$$
\mbox{res} \rho_{\cal D} (a) = \rho_{\cal D}(a)|_{{\cal F}_0}.
$$
This operator as an operator on ${\cal D}_0$ is given more explicitly by the next theorem. Before its formulation we introduce some notations.

Let $\overline V=\{x\in V\mid(\alpha,x)=0, \forall\alpha\in S_0\}$. Denote by $\bar\lambda$ the orthogonal projection of a vector $\lambda\in V$ onto the subspace $\overline V$ and by $\overline{P^\vee}$ the set of projections of the coweights: $\overline{P^\vee}=\{\bar\lambda\mid\lambda\in P^\vee\}\subset\overline V$. The shift operators $\tau(\bar\lambda)$ preserve the affine plane $\mathcal D_0$ and therefore they act on functions $\psi$ on this plane: $\tau(\bar\lambda)\psi(x)=\psi(x-\hbar\bar\lambda)$, where $x\in\mathcal D_0$.
Let $\overline{\cal M}$ be the space of meromorphic functions on ${\cal D}_0$ which are restriction $\psi|_{{\cal D}_0}$, $\psi \in \cal M$.
Define the algebra $\mathbb D$ of difference operators  on ${\cal D}_0$ by $\overline{\mathbb D} = \overline{\cal M} \mathop{\#} \tau(\overline{P^\vee})$. Note that any element $a\in\HH_t$ in the representation~\eqref{ReprDAHA} has the form
\begin{align}
 a=\sum_{\lambda\in P^\vee}\sum_{w\in W}g_{\lambda,w}\tau(\lambda)w, \label{aPolRepr}
\end{align}
where $g_{\lambda,w}\in\mathcal M$

\begin{Th}\label{teorogr}
Let $a\in\mathbb C[Y]^W$ have the form~\eqref{aPolRepr}. Then
\begin{align}
 \res \rho_{\mathcal D}(a) = \sum_{\lambda'\in\overline{P^\vee}}\Big(\sum_{w\in W}\sum_{\lambda\in P^\vee\atop\bar\lambda=\lambda'}g_{\lambda,w}\Big)\Big|_{\mathcal D_0}\tau(\lambda'). \label{resadef}
\end{align}
Formula~\eqref{resadef} defines the algebra homomorphism $\res\colon\rho_{\mathcal D}\big(\mathbb C[Y]^W\big)\to {\overline{\mathbb D}}$.
\end{Th}

\noindent{\bf Proof.} Let us write $v\approx v'$ for $v,v'\in V$ if $wv-v'\in 2\pi iQ^\vee$ for some $w\in W$. Denote $P^x_\lambda=\{\nu\in P^\vee\mid x-\hbar\nu\approx x-\hbar\lambda\}$. Let $f^{x_0}_\lambda(x)\in\mathcal F^W$ be a function which equals $1$ in the points $x\approx x_0-\hbar\lambda$ and vanishes outside the small neighbourhoods of these points. Note that for any $f\in\mathcal F^W$ the function $(af)(x)$ is smooth and that $f\big|_{\mathcal D_0}=0$ implies $(af)\big|_{\mathcal D_0}=0$ since $aI_{\mathcal D}\subset I_{\mathcal D}$. Substitute $f=f^{x_0}_{\lambda_0}$ for $\lambda_0\in P^\vee$ and $x=x_0$ into the formula
\begin{align}
 (af)(x)=\sum_{\lambda\in P^\vee}\sum_{w\in W}g_{\lambda,w}(x)f(x-\hbar\lambda). \label{aPolReprf}
\end{align}
We get that for any $x_0\in V$ the functions $G^{x_0}_\lambda(x)=\sum_{\nu\in P^{x_0}_\lambda}\sum_{w\in W}g_{\nu,w}(x)$ are defined in the point $x=x_0$ and that if $G^{x_0}_\lambda(x_0)\ne0$ for some $x_0\in\mathcal D_0$, $\lambda\in P^\vee$ then there exists $w\in W$ such that $w(x_0-\hbar\lambda)\in\mathcal D_0+2\pi iQ^\vee$.

Let $x_0$ be a generic point of $\mathcal D_0$. Then it follows from $w(x_0-\hbar\lambda)\in\mathcal D_0+2\pi iQ^\vee$ where $w\in W$ that $w_0(x_0-\hbar\lambda)\in\mathcal D_0+2\pi iQ^\vee$ for some $w_0\in W_0$. Indeed, for generic $x_0\in\mathcal D_0$ the formula $w(x_0-\hbar\lambda)\in\mathcal D_0+2\pi iQ^\vee$ implies $w^{-1}S_0\subset V_0$ and hence $wV_0=V_0$. Decomposing the element $w\in W$ as $w=w'w_0$ with $w'\in W_+$, $w_0\in W_0$ (see Remark~\ref{Rem31} and~\cite{H}) we obtain $w'V_0=V_0$. Then by applying Proposition~\ref{AppRemLem2} we derive $w'\mathcal D_0=D_0$ and hence $w_0(x_0-\hbar\lambda)\in\mathcal D_0+2\pi iQ^\vee$.

Let the $W$-invariant function $\eta$ satisfying the conditions $h^\ell_\eta=\hbar$ be generic. Then one can show that $w_0(x_0-\hbar\lambda)\in\mathcal D_0+2\pi iQ^\vee$ implies $x_0-\hbar\lambda\approx x_0-\hbar\bar\lambda$. Note that the right hand side of~\eqref{aPolReprf} at $x=x_0$ has the form $(af)(x_0)=\sum_{\lambda\in P^\vee}\frac1{|P^{x_0}_\lambda|}G^{x_0}_\lambda(x_0)f(x_0-\hbar\lambda)$. Since each term contributes only if $G^{x_0}_\lambda(x_0)\ne0$ we can rewrite the expression as
\begin{align}
 (af)(x_0)=\sum_{\lambda\in P^\vee}\frac1{|P^{x_0}_\lambda|}G^{x_0}_\lambda(x_0)f(x_0-\hbar\bar\lambda), \label{aPolReprfbar}
\end{align}
where we used the $W$-invariance and $2\pi iQ^\vee$-periodicity of the function $f(x)$. One can prove that if $x_0-\hbar\bar\lambda\approx x_0-\hbar\bar\nu$ for some $\lambda,\nu\in P^\vee$ then $\bar\lambda=\bar\nu$. Thus we get $\bar\nu=\bar\lambda$ if $G^{x_0}_\lambda(x_0)\ne0$ and $\nu\in P^{x_0}_\lambda$.
Hence the formula~\eqref{aPolReprfbar} can be rearranged to
\begin{align}
 (af)(x_0)=\sum_{\lambda'\in\overline{P^\vee}}\Big(\sum_{w\in W}\sum_{\lambda\in P^\vee\atop\bar\lambda=\lambda'}g_{\lambda,w}\Big)(x_0)f(x_0-\hbar\lambda'). \label{aPolReprfbar1}
\end{align}
We established the equality~\eqref{aPolReprfbar1} for generic $x_0\in\mathcal D_0$ and $\eta$. Since the left hand side is defined for all $x_0\in\mathcal D_0$ and $\eta$ the equality is valid for all $x_0\in\mathcal D_0$ and $\eta$. Since $\big(\res\rho_{\mathcal D}(a)f\big|_{\mathcal D_0}\big)(x_0)=(af)(x_0)$ for $f\in\mathcal F^W$ the equality~\eqref{resadef} follows. Note that by putting $f=f^{x_0}_{\lambda'_0}$ where $\lambda'_0\in\overline{P^\vee}$ in~\eqref{aPolReprfbar1} we see that each function $\big(\sum_{w\in W}\sum_{\lambda\in P^\vee\atop\bar\lambda=\lambda'_0}g_{\lambda,w}\big)(x_0)$ is defined for any admissible $\eta$ and generic $x_0\in\mathcal D_0$.

Notice that if $A\in\overline{\mathbb D}$ is such that $A\phi=0$ for any $\phi\in\mathcal F_0$ then $A=0$. Indeed consider the action of $A=\sum_{\lambda'\in\overline{P^\vee}}g_{\lambda'}\tau(\lambda')$ where $g_{\lambda'}\in\overline{\cal M}$ on $\phi=f^{x_0}_{\lambda'_0}\big|_{\mathcal D_0}$. Then one yields $g_{\lambda'_0}(x_0)=0$ for generic $x_0\in\mathcal D_0$ and all $\lambda'_0\in\overline{P^\vee}$. It follows that the linear map $\res\colon\rho_{\mathcal D}\big(\mathbb C[Y]^W\big)\to {\overline{\mathbb D}}$ defined by the formula~\eqref{resadef} is an algebra homomorphism. \qed \\

We are going to derive generalized Macdonald-Ruijsenaars operators using the operation $\mbox{res}  \rho_{\cal D}$ from Theorem~\ref{teorogr}. First it is useful to note that restrictions to the invariant functions and restrictions to the quotient module can be interchanged.
More specifically, let $D\in\mathbb D$ be an operator of the form $D=\sum_{\lambda\in P^\vee}g_\lambda\tau(\lambda)$.
Define the operator
\begin{equation}
\label{dbar}
\overline D=\sum_{\lambda'\in\overline{P^\vee}}\Big(\sum_{\lambda\in P^\vee\atop\bar\lambda=\lambda'}g_\lambda\Big)\Big|_{\mathcal D_0}\tau(\lambda')\in\overline{\mathbb D}
\end{equation}
 in the case when for any $\lambda'\in\overline{P^\vee}$ the restricted function $\Big(\sum_{\lambda\in P^\vee\atop\bar\lambda=\lambda'}g_\lambda\Big)\Big|_{\mathcal D_0}$ is well defined meromorphic function on ${\cal D}_0$.

\begin{Prop} \label{Prop42}
 The subspace $\overline{\mathfrak R}=\res\circ\rho_{\mathcal D}\big(\mathbb C[Y]^W\big)\subset\overline{\mathbb D}$ is a commutative subalgebra in $\overline{\mathbb D}$. The map $\overline{\rho_{\cal D}}\colon D\mapsto\overline D$ given by (\ref{dbar}) is well-defined on the subspace $\mathfrak R=\Res\big(\mathbb C[Y]^W\big)\subset\mathbb D$.
The algebra $\overline{\mathfrak R}$ coincides with the image of the map $\overline{\rho_{D}}\colon\mathfrak R\to\overline{\mathbb D}$, thus it is spanned by the elements $\overline{M_\lambda}$, $\lambda\in P^\vee$.
The following diagram of surjective linear maps is commutative:
 \begin{align} \label{diag1}
 &\xymatrix{
  \mathbb C[Y]^{W}\ar@{->}[d]^{\Res}\ar@{->}[r]^{\rho_{\mathcal D}} & \rho_{\mathcal D}\big(\mathbb C[Y]^W\big)\ar@{->}[d]^{\res} \\
  \mathfrak R\ar@{->}[r]^{\overline{\rho_{D}}} & \overline{\mathfrak R}
 }
\end{align}
and the map $\overline{\rho_{D}}\colon\mathfrak R\to\overline{\mathfrak R}$ is an algebra homomorphism.
\end{Prop}

\noindent{\bf Proof.} The space $\overline{\mathfrak R}$ is a commutative algebra as an image of the homomorphism $\res\circ\rho_{\mathcal D}$. The definition of the map $\res\colon\rho(\mathbb C[Y]^W\big)\to\overline{\mathbb D}$ implies that $\res\circ\rho_{\mathcal D}(a)=\overline{\Res(a)}$ for any $a\in\mathbb C[Y]^W$. Due to this formula and surjectivity of the map $\Res\colon\mathbb C[Y]^W\to\mathfrak R$ the operator $\overline D$ is defined for any $D\in\mathfrak R$. We also conclude that $\overline{\mathfrak R}$ coincides with image of the map $\overline{\rho_{D}}\colon\mathfrak R\to\overline{\mathbb D}$ and that the diagram~\eqref{diag1} commutes. Since the operators $M_\lambda$ form a basis in $\mathfrak R$ their images $\overline{M_\lambda}$ span the space $\overline{\mathfrak R}$. The commutativity of the diagram implies in turn that $\overline{\rho_{D}}\colon\mathfrak R\to\overline{\mathfrak R}$ is a homomorphism: $\overline{a'\cdot b'}=\overline{\Res(a)\cdot\Res(b)}
=\overline{\Res(a\cdot b)}=\res\circ\rho_{\mathcal D}(a\cdot b)
=\res\circ\rho_{\mathcal D}(a)\cdot\res\circ\rho_{\mathcal D}(b)= \overline{a'}\cdot\overline{b'}$, where $a,b\in\mathbb C[Y]^W$ such that $a'=\Res(a)$ and $b'=\Res(b)$. \qed \\

We will say that the image $\overline{D}$ of $D\in\mathfrak R$ under $\overline{\rho_{\cal D}}$ is the {\it restriction} of the operator $D$ onto $\mathcal D_0$. The commutative subalgebra $\overline{\mathfrak R}\subset\overline{\mathbb D}$ is the {\it generalized Macdonald-Ruijsenaars system}. The restrictions $\overline{M_{-b_r}}$ and $\overline{M_{-\theta^\vee}}$ of the Macdonald operators $M_{-b_r}$ and $M_{-\theta^\vee}$ are the {\it generalized Macdonald-Ruijsenaars operators}.

\subsection{Integrability of  generalized Macdonald-Ruijsenaars systems}
\label{sec42}

Let $\bar n=\dim\mathcal D_0=n-|S_0|$, where $|S_0|$ is a number of elements in the subset $S_0\subset\Delta$.
The set of indexes $J=\{j\in\mathbb Z_{>0} \mid \alpha_j\notin S_0\}$ has size $\bar n$, let $J=\{j_1,\ldots,j_{\bar n}\}$.
If $j\in J$ then the corresponding fundamental coweight belongs to the subspace $\overline V$, that is $\overline{b_j}=b_j$, $b_j\in\overline{P^\vee}$ and $\tau(\pm b_j)$ preserve the plane $\mathcal D_0$.

\begin{Th} \label{Th43}
 The operators $\overline{M_{-b_{j_1}}},\ldots,\overline{M_{-b_{j_{\bar n}}}}$ are algebraically independent elements of $\overline{\mathfrak R}$.
\end{Th}

\noindent{\bf Proof.} By the formula~\eqref{MlambdaY} with $\lambda=-b_j$ we have
\begin{align}
 M_{-b_j}=g^{-b_j}(x)\tau({-b_j})+\sum_{\nu\prec-b_j}g^{-b_j}_{\nu}(x)\tau(\nu),
\end{align}
where $g^{-b_j}_{\nu}\in\mathcal M$ and $g^{-b_j}(x)=\prod_{\alpha\in\wh R_{-b_j}}\frac{t_\alpha e^{\alpha(x)}-t_\alpha^{-1}}{e^{\alpha(x)}-1}$.

Let $j\in J$. Firstly note that $\nu\prec-b_j$ implies $\bar\nu\ne-b_j$. Indeed, if a weight $\nu\in P^\vee$ satisfies $\nu\prec-b_j$ then it follows the corresponding inequality between the lengths: $|\nu|\le|-b_j|$. If $\bar\nu=\nu$ then $\bar\nu=\nu\ne-b_j$. If $\bar\nu\ne\nu$ then one has $|\bar\nu|<|\nu|$ hence $|\bar\nu|<|-b_j|$ and therefore $\bar\nu\ne-b_j$. Taking into account that $\overline{b_j}=b_j$ one yields
\begin{align} \label{Mbj1}
 \overline{M_{-b_j}}=g^{-b_j}\big|_{\mathcal D_0}\tau({-b_j})+\sum_{\nu'\in\overline{P^\vee}\atop\nu'\ne-b_j}\Big(\sum_{\nu\prec-b_j\atop\bar\nu=\nu'}g^{-b_j}_{\nu}\Big)\Big|_{\mathcal D_0}\tau(\nu').
\end{align}

Let us prove that the function $g^{-b_j}\big|_{\mathcal D_0}$ does not vanish identically. Indeed, otherwise there exists an affine root $\alpha=[\alpha^{(1)},\alpha^{(0)}]\in\wh R_{-b_j}$ such that $e^{\alpha(x)}\equiv e^{(\alpha^{(1)},x)+\alpha^{(0)}}=t_\alpha^{-2}$ for all $x\in\mathcal D_0$. This implies that $\alpha^{(1)}\in R$ is a linear combination of the elements of $S_0$ and, consequently, $(\alpha^{(1)},b_j)=0$. But on the other hand $$\wh R_{-b_j}=\big(\tau({-b_j})\wh R_+\big)\cap\wh R_-=\{[\beta,\hbar k]\mid\beta\in R_-, k\in\mathbb Z, (\beta,b_j)< k\le0\},$$ which implies $(\alpha^{(1)},b_j)<0$. Hence $g^{-b_j}\big|_{\mathcal D_0}\not\equiv0$.

Let us introduce the order on the set $\overline{P^\vee}$: we will write $\nu'\gtrdot\lambda'$ for $\nu',\lambda'\in\overline{P^\vee}$ if $\nu'\ne\lambda'$ and $\nu'-\lambda'=\sum_{i=1}^n\gamma_i\overline{\alpha_i^\vee}$ for some $\gamma_i\in\mathbb Z_{\ge0}$. One can check that this is a partial order on $\overline{P^\vee}$. The ordering $\nu\prec-b_j$ for $\nu\in P^\vee$ implies $\nu>-b_j$, $\bar\nu\ne-b_j$ and hence $\bar\nu\gtrdot-b_j$. So denoting $\bar g^{-b_j}_{\nu'}=\sum_{\nu\prec-b_j\atop\bar\nu=\nu'}g^{-b_j}_{\nu}\in\overline{\mathcal M}$ for $\nu'\ne-b_j$ we obtain
\begin{align} \label{Mbj2}
 \overline{M_{-b_j}}=g^{-b_j}\big|_{\mathcal D_0}\,\tau({-b_j})+\sum_{\nu'\in\overline{P^\vee}\atop\nu'\gtrdot-b_j}\bar g^{-b_j}_{\nu'}\,\tau(\nu'),
\end{align}
where $g^{-b_j}\big|_{\mathcal D_0}\ne0$.
Let $m_1,\ldots,m_{\bar n}\in\mathbb Z_{\ge0}$ and $\lambda'=-\sum_{l=1}^{\bar n}m_l b_{j_l}\in\overline{P^\vee_-}$, then
\begin{align}
 \prod_{l=1}^{\bar n}\Big(\overline{M_{-b_{j_l}}}\Big)^{m_l}=g^{m_1,\ldots,m_{\bar n}}\,\tau(\lambda')+\sum_{\nu'\in\overline{P^\vee}\atop\nu'\gtrdot\lambda'}g^{m_1,\ldots,m_{\bar n}}_{\nu'}\,\tau(\nu'), \label{Fbmonombar}
\end{align}
where $g^{m_1,\ldots,m_{\bar n}},g^{m_1,\ldots,m_{\bar n}}_{\nu'}\in\overline{\mathcal M}$ and $g^{m_1,\ldots,m_{\bar n}}\ne0$. Since the operators~\eqref{Fbmonombar} are linearly independent the elements $\overline{M_{-b_{j_1}}},\ldots,\overline{M_{-b_{j_{\bar n}}}}$ are algebraically independent. \qed \\


Now we introduce some notations. Let $f$ be a function of several variables depending, in particular, on a complex variable $z$. Then we denote by $\T^a_z$ the operator acting on $f$ by the shift:
\begin{align}
 (\T^a_zf)(\ldots,z,\ldots)=f(\ldots,z+a,\ldots), \label{TazDef}
\end{align}
where $a\in\mathbb C$. Let $x_i=(e_i,x)+c_i$ for $i=1,\ldots,N$, where $\{e_i\}$ is a set of pairwise orthogonal vectors in $\mathbb C^N$ and $c_i\in\mathbb C$. Then we have $\T^{\hbar K_i}_{x_i}=\tau(-e_i)$ on the functions $f(x)=f(x_1,\ldots,x_N)$, where $K_i=(e_i,e_i)$.

\subsection{${A_n}$ type examples}
\label{sec43}

The root system $R=A_n$ is defined as $\{e_i-e_j\mid i,j=1,\ldots n+1, i\ne j\}$ considered in the space $V=\{x=\sum_{i=1}^{n+1}x_ie_i\in\mathbb C^{n+1}\mid\sum_{i=1}^{n+1}x_i=0\}$, where $\{e_i\}$ is the standard orthonormal basis in $\mathbb C^{n+1}$. The simple roots are $\alpha_1=e_1-e_2$, $\alpha_2=e_2-e_3$, \ldots, $\alpha_n=e_n-e_{n+1}$ and the maximal root is $\theta=e_1-e_{n+1}$ All fundamental coweights $b_r$, $r=1,\ldots,n$ are minuscule. Any $\wh W$-invariant function $t\colon\wh R\to\mathbb C\backslash\{0\}$ is constant.

We will consider functions $f$ on $V$ as functions of $n+1$ variables satisfying the condition $\prod_{i=1}^{n+1}\T^a_{x_i}f(x_1,\ldots,x_{n+1})=f(x_1,\ldots,x_{n+1})$ $\forall a\in\mathbb C$. Then we have $\T^\hbar_{x_i}f=\tau(-e_i)f=\tau(-\hat e_i)f$, where $\hat e_i$ is the orthogonal projection of the vector $e_i$ to the subspace $V\subset\mathbb C^{n+1}$.

Taking the coweight $b_1=\hat e_1$ in~\eqref{RMop_br} we obtain the operator
\begin{align}
 M_{-b_1}=\sum_{i=1}^{n+1}\prod_{j=1 \atop j\ne i}^{n+1} \frac{te^{x_j-x_i}-t^{-1}}{e^{x_j-x_i}-1}\;\T^\hbar_{x_i}. \label{M1An}
\end{align}

Any root subsystem $R_0$ has irreducible components $R_\ell$ of type $A$ only in its decomposition $R_0=\bigsqcup_{\ell=1}^m R_\ell$. The invariant function $\eta=const$. The generalised Coxeter number is $h^\ell_\eta=k\eta$ for $R_\ell\cong A_{k-1}$. Thus the conditions~\eqref{condInvGCN} lead to $R_1\cong R_2\cong\ldots\cong R_m\cong A_{k-1}$ and $\hbar=k\eta$ for some $k\in\mathbb Z_{\ge2}$. Let us choose
\begin{align}
 S_\ell&=\{e_{N_1+k(\ell-1)+1}-e_{N_1+k(\ell-1)+2},\ldots,e_{N_1+k\ell-1}-e_{N_1+k\ell}\}, &&\ell=1,\ldots,N_2, \label{SellAn}
\end{align}
where $N_1=n+1-mk$ and $N_2=m$. Introduce new variables
\begin{align}
 y_\ell&=\frac1k\sum_{s=0}^{k-1}(x_{N_1+k\ell-s}-s\eta)=\frac1k\sum_{s=0}^{k-1}x_{N_1+k\ell-s}-\frac{k-1}2\eta, && \ell=1,\ldots,N_2.
\end{align}
Note that $x_{N_1+k\ell-s}\big|_{\mathcal D_0}=y_\ell\big|_{\mathcal D_0}+s\eta$ for all $s=0,\ldots,k-1$. The values of the variables $x_i,y_\ell$ with $i=1,\ldots,N_1$, $\ell=1,\ldots,N_2$ satisfying $\sum_{i=1}^{N_1}x_i+k\sum_{\ell=1}^{N_2}y_\ell=-k(k-1)\eta/2$ uniquely determine a point on the affine plane $\mathcal D_0$. The functions $f\colon\mathcal D_0\to\mathbb C$ are functions $f(x_1,\ldots,x_{N_1},y_1,\ldots,y_{N_2})$ satisfying $\prod_{i=1}^{N_1}\T^a_{x_i}\prod_{\ell=1}^{N_2}\T^{a/k}_{y_\ell}f=f$, $\forall a\in\mathbb C$.

\begin{Prop} \label{PropM1Anbar}
Let $t=e^{\eta/2}$. The restriction of the operator~\eqref{M1An} onto $\mathcal D_0$ takes the form
\begin{multline}
 \overline{M_{-b_1}}= \sum_{i=1}^{N_1}
\prod_{j=1\atop j\ne i}^{N_1}\frac{te^{x_j-x_i}-t^{-1}}{e^{x_j-x_i}-1}\;
\prod_{\ell'=1}^{N_2}\frac{qe^{y_{\ell'}-x_i}-q^{-1}}{e^{y_{\ell'}-x_i}-1}  \;\T^\hbar_{x_i}+ \\
+\frac{q-q^{-1}}{t-t^{-1}}\sum_{\ell=1}^{N_2}
  \prod_{j=1}^{N_1}\frac{te^{x_j-y_\ell}-t^{-1}}{e^{x_j-y_\ell}-1}
 \prod_{\ell'=1\atop\ell'\ne\ell}^{N_2}\frac{qe^{y_{\ell'}-y_\ell}-q^{-1}}{e^{y_{\ell'}-y_\ell}-1} \;\T^{\eta}_{y_\ell}. \label{M1Anbar}
\end{multline}
\end{Prop}

\noindent{\bf Proof.} The coefficient at $\T^\hbar_{x_i}$ in $M_{-b_1}$ has non-zero restriction on $\mathcal D_0$ only if $1\le i\le N_1$ or $i=N_1+k\ell$ for some $\ell=1,\ldots,N_2$. These restrictions can be calculated explicitly using the formulae (valid on $\mathcal D_0$)
\begin{align*}
 \prod_{s=0}^{k-1} g(x_{N_1+k\ell'-s}-a;t,1)&=g(y_{\ell'}-a;t^k,1), \\
\prod_{s=1}^{k-1} g(x_{N_1+k\ell-s}-x_{N_1+k\ell};t,1)&=g(0;t^k,t)=\dfrac{t^k-t^{-k}}{t-t^{-1}},
\end{align*}
where $a\in\mathbb C$ and $g(z;b,c)=\dfrac{be^{z}-b^{-1}}{ce^{z}-c^{-1}}$. Further since $\overline{e_i}=e_i$ for $i=1,\ldots,N_1$ the restriction of $\T^\hbar_{x_i}=\tau(-e_i)$ in this case coincides with $\T^\hbar_{x_i}$. For $i=N_1+k\ell$ the restriction of $\T^\hbar_{x_i}=\tau(-e_i)$ is equal to $\T_{y_\ell}^{\hbar/k}$ due to $\overline{e_i}=\frac1k\sum_{s=0}^{k-1}e_{N_1+k\ell-s}$ and $y_\ell=(\overline{e_i},x)-(k-1)\eta/2$ (see the end of Subsection~\ref{sec42}). Taking into account the condition $\hbar=k\eta$ and its exponential form $q=t^k$ we can write the restriction of~\eqref{M1An} in the form~\eqref{M1Anbar}. \qed \\

Note that the operator~\eqref{M1Anbar} does not effectively depend on $k$, it depends on $N_1$, $N_2$, $\hbar$, $\eta$ only, where it is assumed that $\hbar/\eta\in\mathbb Z_{>0}$.
The operator~\eqref{M1Anbar} was previously obtained and investigated in~\cite{SV} for general parameters $\hbar$ and $\eta$ and earlier in~\cite{Chalykhbisp} in the case $N_2=1$.

The formula~\eqref{RMop_br} gives the explicit form of $n$ Macdonald operators
\begin{align}
M_{-b_r}=\sum_{I\subset\{1,\ldots,n+1\}\atop|I|=r}\prod_{i\in I\atop j\notin I}\frac{te^{x_j-x_i}-t^{-1}}{e^{x_j-x_i}-1}\prod_{i\in I}\T^\hbar_{x_i}, &&&r=1,\ldots,n, \label{Mbr}
\end{align}
where $|I|$ is the number of elements in the set $I\subset\{1,\ldots,n+1\}$. Similar to the proof of Proposition~\ref{PropM1Anbar} we calculate the restrictions of the operators~\eqref{Mbr} onto $\mathcal D_0$. We obtain
\begin{multline}
 \overline{M_{-b_r}}=\sum_{I\subset\{1,\ldots,N_1\}\atop {m_\ell\in\mathbb Z_{\ge0}\atop|I|+\sum_{\ell=1}^{N_2}m_\ell=r}}\prod_{i,j=1\atop {i\in I\atop j\notin I}}^{N_1}\frac{te^{x_j-x_i}-t^{-1}}{e^{x_j-x_i}-1}\;\prod_{\ell,\ell'=1}^{N_2}\prod_{s=0}^{m_\ell-1}\frac{qt^{-s}e^{y_{\ell'}-y_\ell}-q^{-1}t^s}{t^{m_{\ell'}-s}e^{y_{\ell'}-y_\ell}-t^{s-m_{\ell'}}}\; \times \\
\prod_{\ell=1\atop i\in I}^{N_2}
\frac{qe^{y_\ell-x_i}-q^{-1}}{t^{m_\ell}e^{y_\ell-x_i}-t^{-m_\ell}}
\prod_{\ell=1\atop j\notin I}^{N_2}\frac{te^{x_j-y_\ell}-t^{-1}}{t^{1-m_\ell}e^{x_j-y_\ell}-t^{m_\ell-1}}  \;\prod_{i\in I}\T^\hbar_{x_i}\prod_{\ell=1}^{N_2}\T^{m_\ell\eta}_{y_\ell}, \label{MbrG}
\end{multline}
where $q=e^{\hbar/2}$, $t=e^{\eta/2}$. In the case $r=1$ the operator~\eqref{MbrG} takes the form~\eqref{M1Anbar}.

\begin{Th} \label{PropHOMOcomm}
 For any values of the parameters $\hbar$ and $\eta$ the operators $\overline{M_{-b_r}}$ given by the formula~\eqref{MbrG}, $r\in\mathbb Z_{>0}$, pairwise commute. The operators $\overline{M_{-b_r}}$ with $r=1,\ldots,N_1+N_2-1$ are algebraically independent at generic values of $\hbar$ and $\eta$.
\end{Th}

\noindent{\bf Proof.} For any $r,r'\in\mathbb Z_{>0}$ the commutator has the form
\begin{equation}
\Big[\overline{M_{-b_r}},\overline{M_{-b_{r'}}}\Big]=\sum_{k_1,\ldots,k_{N_1}\in\mathbb Z_{\ge0}\atop m_1,\ldots,m_{N_2}\in\mathbb Z_{\ge0}}
 C^{k_1,\ldots,k_{N_1}}_{m_1,\ldots,m_{N_2}}(x,y;q,t) \prod_{i=1}^{N_1}\T^{k_i\hbar}_{x_i}\prod_{\ell=1}^{N_2}\T^{m_\ell\eta}_{y_\ell}, \label{commCkm}
\end{equation}
where the coefficients $C^{k_1,\ldots,k_{N_1}}_{m_1,\ldots,m_{N_2}}(x,y;q,t)$ are rational functions of $q=e^{\hbar/2}$ and $t=e^{\eta/2}$. It follows from Proposition~\ref{Prop42} that these functions equal zero if $\hbar,\eta\in\mathbb C\backslash\pi i\mathbb Q$ such that $\hbar/\eta\in\mathbb Z_{>1}$ and $N_1+\frac\hbar\eta N_2>\max(r,r')$. Thus the coefficients vanish for any parameters $\hbar$ and $\eta$, and the operators commute.

The second part follows from the fact that the operators $\overline{M_{-b_r}}$ at $\hbar=\eta$ are specialisations of the Macdonald operators $M_{-b_r}$ with $n=N_1+N_2-1$, which are algebraically independent. \qed

\subsection{$B_n$ type examples}
\label{sec44}

The root system $R=B_n$ is defined as the set $\{\epsilon e_i+\epsilon' e_j\mid 1\le i<j\le n,\epsilon,\epsilon'=\pm1\}\cup\{\epsilon e_i\mid1\le i\le n,\epsilon=\pm1\}$. The simple roots are $\alpha_1=e_1-e_2$, $\alpha_2=e_2-e_3$, \ldots, $\alpha_{n-1}=e_{n-1}-e_n$, $\alpha_n=e_n$. The maximal root is $\theta=e_1+e_2$. The group $\wh W$ has two orbits on $\wh R$ so the function $t\colon\wh R\to\mathbb C\backslash\{0\}$ is determined by the values $t_1=t_{\pm e_i\pm e_j}$ and $t_n=t_{\pm e_i}$.  There is a unique minuscule coweight $b_1=e_1$. The corresponding Macdonald operator is
\begin{multline}
 M_{-e_1}=\sum_{i=1}^n\Bigg(\prod_{j=1 \atop j\ne i}^n \frac{(t_1e^{x_i-x_j}-t_1^{-1})(t_1e^{x_i+x_j}-t_1^{-1})}{(e^{x_i-x_j}-1)(e^{x_i+x_j}-1)}\frac{t_ne^{x_i}-t_n^{-1}}{e^{x_i}-1}\;\T^{-\hbar}_{x_i}+ \\
  +\prod_{j=1 \atop j\ne i}^n \frac{(t_1e^{x_j-x_i}-t_1^{-1})(t_1e^{-x_i-x_j}-t_1^{-1})}{(e^{x_j-x_i}-1)(e^{-x_i-x_j}-1)}\frac{t_ne^{-x_i}-t_n^{-1}}{e^{-x_i}-1}\;\T^{\hbar}_{x_i}\Bigg), \label{M1Bn}
\end{multline}
where $x_i=(e_i,x)$ are coordinates on $V$.

The invariant function $\eta$ satisfies $\eta_1=\eta_2=\ldots=\eta_{n-1}$ and $\eta_\theta=\eta_1$. The generalised Coxeter number is $h^\ell_\eta=k\eta_1$ for $R_\ell\cong A_{k-1}$ and $h^\ell_\eta=2(k-1)\eta_1+2\eta_n$ if $R_\ell\cong B_k$.

We have two cases $\alpha_n\notin S_0$ and $\alpha_n\in S_0$ depending on the choice of $S_0\subset\Delta$. Consider first the case $\alpha_n\notin S_0$. All the irreducible components $R_\ell$ are of type $A$. Then the conditions $h^\ell_\eta=\hbar$ imply $\hbar=k\eta_1$ for some $k\in\mathbb Z_{>0}$ and $R_\ell\cong A_{k-1}$ for all $\ell$. Hence $q=t_1^k$ and we can choose
\begin{align}
 S_\ell&=\{e_{N_1+k(\ell-1)+1}-e_{N_1+k(\ell-1)+2},\ldots,e_{N_1+k\ell-1}-e_{N_1+k\ell}\}, &&\ell=1,\ldots,N_2, \label{SellBn}
\end{align}
where $N_1=n-mk$ and $N_2=m$. Define new variables
\begin{align}
 y_\ell&=\frac1k\sum_{s=0}^{k-1}x_{N_1+k\ell-s}=\frac{k-1}2\eta_1+\frac1k\sum_{s=0}^{k-1}(x_{N_1+k\ell-s}-s\eta_1), && \ell=1,\ldots,N_2. \label{ylBn}
\end{align}
Note that $x_{N_1+k\ell-s}\big|_{\mathcal D_0}=y_\ell\big|_{\mathcal D_0}-(k-1)\eta_1/2+s\eta_1$ for all $s=0,\ldots,k-1$. The variables $x_i,y_\ell$ with $i=1,\ldots,N_1$, $\ell=1,\ldots,N_2$ form a system of orthogonal coordinates on the affine plane $\mathcal D_0$.

\begin{Prop}
Let $t_1=e^{\eta_1/2}$. In the case $\alpha_n\notin S_0$ the restriction of the operator~\eqref{M1Bn} onto $\mathcal D_0$ takes the form
\begin{align} \label{M1BnG}
 \overline{M_{-e_1}}= \sum_{i=1\atop\epsilon=\pm1}^{N_1}
A_i^\epsilon(x,y)  \;\T^{-\epsilon\hbar}_{x_i}
+\frac{q-q^{-1}}{t_1-t_1^{-1}}\sum_{\ell=1\atop\epsilon=\pm1}^{N_2}B_\ell^\epsilon(x,y)
 \;\T^{-\epsilon\eta_1}_{y_\ell},
\end{align}
where
\begin{multline}
A_i^\epsilon(x,y)=\prod_{j=1\atop j\ne i}^{N_1}\big(g(\epsilon x_i-x_j; t_1,1)g(\epsilon x_i+x_j;t_1,1)\big) \\
\prod_{\ell'=1}^{N_2}\big(g(\epsilon x_i-y_{\ell'};q^{1/2}t_1^{1/2},q^{-1/2}t_1^{1/2})g(\epsilon x_i+y_{\ell'};q^{1/2}t_1^{1/2},q^{-1/2}t_1^{1/2})\big)\;g(\epsilon x_i;t_n,1),
\end{multline}
\begin{multline}
B_\ell^\epsilon(x,y)=\prod_{j=1}^{N_1}\big(g(\epsilon y_\ell-x_j;q^{1/2}t_1^{1/2},q^{1/2}t_1^{-1/2})g(\epsilon y_\ell+x_j;q^{1/2}t_1^{1/2},q^{1/2}t_1^{-1/2})\big) \\
 \prod_{\ell'=1\atop\ell'\ne\ell}^{N_2}\big(g(\epsilon y_\ell-y_{\ell'};q,1)g(\epsilon y_\ell+y_{\ell'};q,1)\big)\;
 g(\epsilon y_\ell;q^{1/2}t_1^{-1/2}t_n,q^{1/2}t_1^{-1/2})
 g(2\epsilon y_\ell;qt_1^{-1},1),
\end{multline}
and
\begin{align} \label{gdef}
 g(z;a,b)=\frac{a e^{z}-a^{-1}}{b e^{z}-b^{-1}}, &&&a,b\in\mathbb C.
\end{align}
\end{Prop}

\noindent{\bf Proof.} The coefficient at $\T^{-\epsilon\hbar}_{x_i}$ in $M_{-e_1}$ has non-zero restriction on $\mathcal D_0$ only if $1\le i\le N_1$ or $i=N_1+k\ell-(1+\epsilon)(k-1)/2$ for some $\ell=1,\ldots,N_2$. To calculate these coefficients we first use
\begin{gather}
 \prod_{s=0}^{k-1} g(a\pm x_{N_1+k\ell'-s};t_1,1)=g(a\pm y_{\ell'};q^{1/2}t_1^{1/2},q^{-1/2}t_1^{1/2}), \label{gProofB1} \\
\prod_{s=0\atop s\ne N_1+k\ell-i}^{k-1}g(\epsilon x_i-x_{N_1+k\ell-s};t_1,1)g(\epsilon x_i+x_{N_1+k\ell-s};t_1,1)=g(0;q,t)g(2\epsilon y_\ell;qt_1^{-1},1), \label{gProofB2}
\end{gather}
where $i=N_1+k\ell-(1+\epsilon)(k-1)/2$.
Further, if $i=N_1+k\ell-(1+\epsilon)(k-1)/2$ then we have $\epsilon x_i=\epsilon(y_\ell-(k-1)\eta_1/2+(1+\epsilon)(k-1)\eta_1/2)=\epsilon y_\ell+(k-1)\eta_1/2$ and hence
\begin{align}
 g(\epsilon x_i+a;b,c)=g(\epsilon y_\ell+a;bq^{1/2}t_1^{-1/2},cq^{1/2}t_1^{-1/2}) \label{gProofB3}
\end{align}
for any $b,c\ne0$. Thus we obtain~\eqref{M1BnG}. \qed \\

If $\alpha_n\in S_0$ then we suppose that the numeration of the components $S_\ell$ is such that $\alpha_n\in S_m$. Then $R_\ell\cong A_{k-1}$ for $\ell=1,\ldots,m-1$ and $R_m\cong B_{\tilde k}$, where $k=|S_1|+1$ and $\tilde k=|S_m|$. Put $N_1=n-(m-1)k-\tilde k$, $N_2=m-1$, and let $S_1,\ldots,S_{m-1}$ have the form~\eqref{SellBn}, while
\begin{align}
 S_m&=\{e_{n-\tilde k+1}-e_{n-\tilde k+2},e_{n-\tilde k+2}-e_{n-\tilde k+3},\ldots,e_{n-1}-e_n,e_n\}. \label{SmBn}
\end{align}
The variables $x_1,\ldots,x_{N_1}$ together with~\eqref{ylBn} are coordinates on $\mathcal D_0$.

\begin{Prop}
 Let $t_1=e^{\eta_1/2}$ and $t_n=e^{\eta_n/2}$. In the case $\alpha_n\in S_0$ the restriction of the operator~\eqref{M1Bn} onto $\mathcal D_0$ takes the form
\begin{align}
 \overline{M_{-e_1}}= \sum_{i=1\atop\epsilon=\pm1}^{N_1} \wt A_i^\epsilon(x,y)
 \;\T^{-\epsilon\hbar}_{x_i}
 +\frac{q-q^{-1}}{t_1-t_1^{-1}}\sum_{\ell=1\atop\epsilon=\pm1}^{N_2} \wt B_\ell^\epsilon(x,y)
   \;\T^{-\epsilon\eta_1}_{y_\ell}  + \wt  C(x,y), \label{M1resBnk_m}
\end{align}
where
\begin{multline}
 \wt A_i^\epsilon(x,y)=g(\epsilon x_i;t_1t_n^{-1},q^{-1/2})g(\epsilon x_i;q^{1/2}t_1,1)
\prod_{j=1\atop j\ne i}^{N_1}g(\epsilon x_i-x_j;t_1,1)g(\epsilon x_i+x_j;t_1,1) \times \\
\prod_{\ell'=1}^{N_2}g(\epsilon x_i-y_{\ell'};q^{1/2}t_1^{1/2},q^{-1/2}t_1^{1/2})g(\epsilon x_i+y_{\ell'};q^{1/2}t_1^{1/2},q^{-1/2}t_1^{1/2}),
\end{multline}
\begin{multline}
 \wt B_\ell^\epsilon(x,y)=g(\epsilon y_\ell;q^{1/2}t_1^{1/2}t_n^{-1},t_1^{-1/2})
 g(\epsilon y_\ell;qt_1^{1/2},q^{1/2}t_1^{-1/2}) g(2\epsilon y_\ell;qt_1^{-1},1) \times \\
\prod_{j=1}^{N_1} g(\epsilon y_\ell-x_j;q^{1/2}t_1^{1/2},q^{1/2}t_1^{-1/2})g(\epsilon y_\ell+x_j;q^{1/2}t_1^{1/2},q^{1/2}t_1^{-1/2}) \times \\
 \prod_{\ell'=1\atop\ell'\ne\ell}^{N_2}g(\epsilon y_\ell-y_{\ell'};q,1)g(\epsilon y_\ell+y_{\ell'};q,1),
\end{multline}
\begin{multline}
 \wt C (x,y)=\frac{(q^{1/2}t_1t_n^{-1}-q^{-1/2}t_1^{-1}t_n)(q^{1/2}+q^{-1/2})}{t_1-t_1^{-1}} \prod_{j=1}^{N_1}g(-x_j;q^{1/2}t_1,q^{1/2}) g(x_j;q^{1/2}t_1,q^{1/2}) \times \\
 \prod_{\ell'=1}^{N_2} g(-y_{\ell'};qt_1^{1/2},t_1^{1/2})
g(y_{\ell'};qt_1^{1/2},t_1^{1/2}),
\end{multline}
and $g(z; a,b)$ is defined by~\eqref{gdef}.
\end{Prop}

\noindent{\bf Proof.} Note that $x_{n-s}=\eta_n+s\eta_1$ for all $s=0,\ldots,\tilde k-1$. So the functional coefficient at $\T^{-\epsilon\hbar}_{x_i}$ has non-zero restriction on $\mathcal D_0$ only if $1\le i\le N_1$ or $i=N_1+k\ell-(1+\epsilon)(k-1)/2$ for some $\ell=1,\ldots,N_2$ or $i-n+\tilde k=\epsilon=+1$. To calculate the coefficients for the first two cases we use the formulae~\eqref{gProofB1}, \eqref{gProofB2}, \eqref{gProofB3} and the formula $$g(a;t_n,1)\prod_{s=0}^{\tilde k-1}g(a-x_{n-s};t_1,1)g(a+x_{n-s};t_1,1)=g(a;t_1t_n^{-1};q^{1/2})g(a;q^{1/2}t_1,1).$$
The restriction of $\T^{-\hbar}_{x_{n-\tilde k+1}}$ yields the identity operator $1$. The corresponding coefficient is calculated using
\begin{gather*}
g(x_{n-\tilde k+1}+a;b,c)=g(a;bt_1^{\tilde k-1}t_n,ct_1^{\tilde k-1}t_n)=g(a;bq^{1/2},cq^{1/2}), \\
g(x_i;t_n,1)\prod\nolimits_{s=0}^{\tilde k-2}g(x_i-x_{n-s};t_1,1)g(x_i+x_{n-s};t_1,1)=g(0;q^{1/2}t_1t_n^{-1},t_1)g(0;q,q^{1/2}),
\end{gather*}
where $i=n-\tilde k+1$ and we used the condition $2(\tilde k-1)\eta_1+2\eta_n=\hbar$ in the form  $t_1^{\tilde k-1}t_n=\pm q^{1/2}$. \qed

\begin{Rem} \normalfont The operator~\eqref{M1resBnk_m} can be rewritten in the form
\begin{align}
 \overline{M_{-e_1}}= \sum_{i=1\atop\epsilon=\pm1}^{N_1} \wt A_i^\epsilon(x,y)
 \;\Big(\T^{-\epsilon\hbar}_{x_i}-1\Big)
 +\frac{q-q^{-1}}{t_1-t_1^{-1}}\sum_{\ell=1\atop\epsilon=\pm1}^{N_2} \wt B_\ell^\epsilon(x,y)
   \;\Big(\T^{-\epsilon\eta_1}_{y_\ell}-1\Big)  + m_{-e_1,t}\, , \label{M1resBnk_mConst}
\end{align}
where $m_{-e_1,t}=\frac{t_1^n-t_1^{-n}}{t_1-t_1^{-1}}(t_nt_1^{n-1}+t_n^{-1}t_1^{1-n})$.
\end{Rem}

\begin{Rem} \normalfont Rational degeneration of the operators~\eqref{M1BnG} and \eqref{M1resBnk_m} can be obtained by substitution $x_i\to x_i/\omega$, $y_i\to y_i/\omega$, $\eta_\alpha\to\eta_\alpha/\omega$, $\hbar\to\hbar/\omega$ and taking the limit $\omega\to\infty$. This gives particular case of the operator considered by Sergeev and Veselov in~\cite[formula~(94)]{SVBC}, see also \cite{Fbisp}.
\end{Rem}

\subsection{An example from root system ${F_4}$}
\label{sec46}

The root system $R=F_4$ is given in $V=\mathbb C^4$ by
\begin{multline}
 R=\{\epsilon e_i\mid1\le i\le4,\epsilon=\pm1\}\cup\{\epsilon_1e_i+\epsilon_2e_j\mid1\le i<j\le4,\epsilon_1,\epsilon_2=\pm1\}\cup \\
  \cup\Big\{\frac12\sum_{k=1}^4\nu_k e_k\mid\nu_1,\ldots,\nu_4=\pm1\Big\}. \label{rsF4}
\end{multline}
The simple roots are $\alpha_1=e_2-e_3$, $\alpha_2=e_3-e_4$, $\alpha_3=e_4$ and $\alpha_4=\frac12(e_1-e_2-e_3-e_4)$. The maximal root is $\theta=e_1+e_2$. The corresponding Weyl group $W$ has two orbits on $R$. There are no minuscule coweights in this case. The $W$-invariant function $\eta$ satisfies $\eta_1=\eta_2$ and $\eta_3=\eta_4$. We impose $t_0=t_1=t_2=e^{\eta_1/2}$ and $t_3=t_4=e^{\eta_3/2}$. The Macdonald operator~\eqref{RMop_th} takes the form
\begin{multline}
 M_{-\theta^\vee}=\sum_{1\le i<j\le4\atop\epsilon_1,\epsilon_2=\pm1}g(\epsilon_1x_i;t_3,1)g(\epsilon_2x_j;t_3,1)
 \prod_{k\ne i,j\atop\epsilon_3=\pm1}\big(g(\epsilon_1x_i+\epsilon_3x_k;t_1,1)g(\epsilon_2x_j+\epsilon_3x_k;t_1,1)\big)\times \\
\prod_{\nu_1,\ldots,\nu_4=\pm1\atop\nu_i=\epsilon_1,\,\nu_j=\epsilon_2}g\big(\frac12\sum_{k=1}^4\nu_kx_k;t_3,1\big)
g(\epsilon_1x_i+\epsilon_2x_j;t_1,1)g(\epsilon_1x_i+\epsilon_2x_j;q^{-1}t_1,q^{-1})\Big(\T^{-\epsilon_1\hbar}_{x_i}\T^{-\epsilon_2\hbar}_{x_j}-1\Big)+ \\
+m_{-\theta^\vee,t}, \label{MRopF4}
\end{multline}
where $g(z; a,b)$ is defined by~\eqref{gdef}.

Let us consider the simplest example $S_0=\{\alpha_3\}\cong A_1$. Then $\mathcal D_0=\{x\in V\mid x_4=\eta_3\}$ and the condition~\eqref{condInvGCN} becomes $2\eta_3=\hbar$. The restriction of the operator~\eqref{MRopF4} takes the form
\begin{multline}
 \overline{M_{-\theta^\vee}}=\sum_{1\le i<j\le3\atop\epsilon_1,\epsilon_2=\pm1}g(\epsilon_1x_i;q^{1/2},1)g(\epsilon_1x_i;t_1q^{1/2},q^{1/2}) g(\epsilon_1x_i;t_1q^{-1/2},q^{-1/2})g(\epsilon_2x_j;q^{1/2},1) \times \\
 g(\epsilon_2x_j;t_1q^{1/2},q^{1/2}) g(\epsilon_2x_j;t_1q^{-1/2},q^{-1/2})
 g(\epsilon_1x_i+\epsilon_2x_j;t_1,1)g(\epsilon_1x_i+\epsilon_2x_j;q^{-1}t_1,q^{-1}) \times \\
 \prod_{\epsilon_3=\pm1}\Big(g(\epsilon_1x_i+\epsilon_3x_k;t_1,1)g(\epsilon_2x_j+\epsilon_3x_k;t_1,1)g\big(\frac12(\epsilon_1 x_i+\epsilon_2 x_j+\epsilon_3 x_k);q^{3/4},q^{-1/4}\big) \Big) \times \\
\Big(\T^{-\epsilon_1\hbar}_{x_i}\T^{-\epsilon_2\hbar}_{x_j}-1\Big)+ \notag
\end{multline} \\[-10pt]
\begin{multline}
+(q^{1/2}+q^{-1/2})\sum_{1\le i\le3\atop\epsilon_1=\pm1}
g(\epsilon_1x_i;q^{1/2},1)
 \prod_{l\ne i,4\atop\epsilon_3=\pm1}\big(g(\epsilon_1x_i+\epsilon_3x_l;t_1,1)g(\epsilon_3x_l;q^{1/2}t_1,q^{1/2})\big) \times\\
\prod_{\nu_1,\ldots,\nu_3=\pm1\atop\nu_i=\epsilon_1}g\big(\frac12\sum_{l=1}^3\nu_lx_l;q^{3/4},q^{1/4}\big)\;
g(\epsilon_1x_i;t_1q^{1/2},q^{1/2})
g(\epsilon_1x_i;q^{-1/2}t_1,q^{-1/2})\Big(\T^{-\epsilon_1\hbar}_{x_i}-1\Big)+ \\
 +m_{-\theta^\vee,t}\,,
\label{MRopF4G}
\end{multline}
where in every term in the sum $k$ is such that $1\le k\le3$, $k\ne i,j$; and $t_1=e^{\eta_1/2}$, $q^{1/2}=e^{\hbar/4}$ are parameters.

\begin{Rem} \normalfont
 It is possible to obtain a commuting operator to~\eqref{MRopF4G} by restriction of the relevant Macdonald operators for the small weights~\cite{EmvD}. Similar one can get Hamiltonians and quantum integrals starting from the Macdonald operators for the small weights in the case of other exceptional root systems.
\end{Rem}

\section{The case $\pmb{(C^\vee_n,C_n)}$}
\label{sec5}

The affine root system of type $(C^\vee_n,C_n)$ is given by
\begin{multline} \label{sysCC}
\SS=\{[\epsilon_1e_i+\epsilon_2e_j,\hbar k]\mid1\le i<j\le n; \epsilon_1,\epsilon_2=\pm1;k\in\mathbb Z\}\cup \\
\cup\{\epsilon\cdot[e_i,\frac{\hbar k}2]\mid 1\le i\le n;\epsilon=\pm1,\pm2;k\in\mathbb Z\}\subset\wh V
\end{multline}
(see~\cite{Mac1,Mac2}). Let $\wt W$ be the group generated by the affine reflections $s_\alpha$ with $\alpha\in\SS$. The corresponding DAHA is defined in~\cite{Sahi} (see also~\cite{Noumi}). We recall the associated commuting difference operators below.

Let $R=C_n$, that is
\begin{align}
 R=\{\epsilon_1e_i+\epsilon_2e_j\mid1\le i<j\le n; \epsilon_1,\epsilon_2=\pm1\}\cup\{\epsilon 2 e_i\mid 1\le i\le n;\epsilon=\pm1\}, \label{RCdef} \\
 \Delta=\{\alpha_i=e_i-e_{i+1}\mid i=1,\ldots,n-1\}\cup\{\alpha_n=2e_n\}\subset R,
\end{align}
let $\wh R=\{[\alpha,\hbar k]\mid\alpha\in R;k\in\mathbb Z\}$ and $\wt\Delta=\Delta\cup\{\alpha_0=[-\theta,\hbar]=[-2e_1,\hbar]\}$. Then $\wt W = \{s_\alpha\mid\alpha\in\SS\}=\{s_\alpha\mid\alpha\in\wh R\}$ is the affine Weyl group for $C_n$.

Let $Q^\vee$, $P^\vee$, $P$ and $W$ be the lattices of coroots, of coweights, of weights and the Weyl group for the root system $R=C_n$ respectively. Let also $\wh P\subset\wh V$ be the subset defined in Subsection~\ref{sec21} for $R=C_n$.

\subsection[Noumi representation and Koornwinder operator]{Noumi representation of $\pmb{(C^\vee_n,C_n)}$ type DAHA and Koornwinder operator}
\label{sec51}

Let $t,u\colon\wh R\to\mathbb C\backslash\{0\}$ be such that $t_\alpha=t_{\tilde w\alpha}=t(\alpha)$, $u_{\alpha}=u_{\tilde w\alpha}=u(\alpha)$ for any $\tilde w\in\wt W$, $\alpha\in\wh R$. Put $u_{\pm e_i\pm e_j}=1$, $t_i=t_{\alpha_i}$, $u_i=u_{\alpha_i}$. So we have $t_1=t_2=\ldots=t_{n-1}$ and $u_1=u_2=\ldots=u_{n-1}=1$.

\begin{Def}{\normalfont\cite{Sahi}}
 Double Affine Hecke Algebra of type $(C^\vee_n,C_n)$ is the algebra $\HH_{t,u}$ generated by the elements $T_0,T_1,\ldots,T_n$, and by the set of pairwise commuting elements $\{X^\mu\}_{\mu\in\wh P}$ with the relations~\eqref{DAHAxx}--\eqref{DAHAttt}, \eqref{DAHAtx0} and {\normalfont
\begin{align}
 T_i^{-1}X^\mu T_i^{-1}&=X^{\mu-\alpha_i}+(u_i-u_i^{-1})X^{\mu-\alpha_i/2}\,T_i^{-1}, & &\text{if $(\mu,\alpha_i^\vee)=1$}, & &\mu\in\wh P,
\end{align}}
where $i=0,1,\ldots,n$.
\end{Def}

The algebra $\HH_{t,u}$ depends on six independent parameters $t_1$, $t_n$, $t_0$, $u_n$, $u_0$ and $\hbar$. Note that the elements $\pi_r$ are not included in this case. One can also define an `extended' version of $(C_n^\vee,C_n)$ type DAHA where affine Weyl group $\wt W=W\ltimes\tau(Q^\vee)$ is replaced with the extended affine Weyl group $\wh W=W\ltimes\tau(P^\vee)$. However the $\wh W$-invariance of the functions $t, u$ leads to less independent parameters in this algebra and subsequently in the commuting difference operators. Note also that if $u_0=u_n=1$ then the algebra $\HH_{t,u}$ can be interpreted as a `non-extended' version of DAHA of type $\wh{C_n}$.

Consider the action of the generators of $\HH_{t,u}$ on the space $\mathcal F=C^\infty(V)^{2\pi i Q^\vee}$ given by the Noumi representation \cite{Noumi}
\begin{align}
 X^\mu&=e^{\mu(x)}, & & T_i=t_i+t_i^{-1}\frac{(1-t_iu_ie^{-\alpha_i(x)/2})(1+t_iu_i^{-1}e^{-\alpha_i(x)/2})}{1-e^{-\alpha_i(x)}}\,(s_i-1), \label{ReprDAHACC}
\end{align}
where $\mu\in\wh P$, $i=0,1,\ldots,n$.

The elements $T_{\tilde w}\in\HH_{t,u}$ with $\tilde w\in \wt W$ and Cherednik-Dunkl operators $Y^\lambda\in\HH_{t,u}$ with $\lambda\in Q^\vee$ are defined by the formula~\eqref{Twdef}, where $\pi_r=1$, and
\begin{align} \label{YdefCC}
 Y^{\lambda-\nu}=T_{\tau(\lambda)}T_{\tau(\nu)}^{-1}, &&&\lambda,\nu\in P^\vee_+\cap Q^\vee.
\end{align}
The operators $Y^{\lambda}$ commute with each other and satisfy the relations~\cite{Noumi,Mac2}
\begin{align}
 T_iY^\lambda-Y^{s_i\lambda}T_i&=\frac{t_i-t_i^{-1}+(\bar t_i-\bar t_i^{-1})Y^{-\alpha_i^\vee}}{1-Y^{-2\alpha_i^\vee}}(Y^\lambda-Y^{s_i\lambda}), & &i=1,\ldots,n, \label{TiYCC}
\end{align}
where $\lambda\in Q^\vee$, $\bar t_i=t_i$ for $i=1,\ldots,n-1$ and $\bar t_n=t_0$.

 Let $\mathbb C[Y]$ be the commutative subalgebra of $\HH_{t,u}$ spanned by the Cherednik-Dunkl operators~\eqref{YdefCC}. The relation~\eqref{TiYCC} implies that the elements of the subalgebra $\mathbb C[Y]^W\subset\mathbb C[Y]$ preserve the subspace of $W$-invariant functions $\mathcal F$, that is $m_\lambda(Y)\mathcal F^W\subset\mathcal F^W$, where $\lambda\in Q^\vee$ and $m_\lambda(Y)$ is the operator defined by the formula~\eqref{flambdadef}. By restricting these operators to the subspace $\mathcal F^W$ we obtain the difference operators $M_\lambda$ given by formula~\eqref{Flamndadef} with~\eqref{VIAg}.

Thus we get a commutative algebra $\mathfrak R$ spanned by $M_\lambda$ with $\lambda\in Q^\vee$. It is generated by the algebraically independent elements $M_{\nu_i}=M_{-\nu_i}$, where $\nu_i=e_1+e_2+\ldots+e_i$ and $i=1,\ldots,n$. The first of the generators is {\it the Koornwinder operator}~\cite{Noumi,Koornwinder}
\begin{align} \label{MacKoor}
\begin{split}
 M_{-e_1}=&M_{-\theta^\vee}=\sum_{i=1\atop\epsilon=\pm1}^n(t_0t_n)^{-1}
\prod_{j=1\atop j\ne i}^n\big(g(\epsilon x_i-x_j;t_1, 1)g(\epsilon x_i+x_j;t_1, 1)\big)\times \\
&\frac{(1-t_0u_0q^{-1}e^{\epsilon x_i})(1+t_0u_0^{-1}q^{-1}e^{\epsilon x_i})(1-t_nu_ne^{\epsilon x_i})(1+t_nu_n^{-1}e^{\epsilon x_i})}
{(1-q^{-2}e^{2\epsilon x_i})(1-e^{2\epsilon x_i})}
\Big(\T^{-\epsilon\hbar}_{x_i}-1\Big)+\kappa_t\,,
\end{split}
\end{align}
where $g(z; a,b) = \dfrac{a e^z - a^{-1}}{b e^z - b^{-1}}$ and  $\kappa_t=\dfrac{t_1^n-t_1^{-n}}{t_1-t_1^{-1}}(t_0t_nt_1^{n-1}+t_0^{-1}t_n^{-1}t_1^{1-n})$.

\subsection{Invariant ideals and generalized Koornwinder operator}
\label{sec52}

Let $\eta\colon R\to\mathbb C$ be a $W$-invariant function. It is defined by independent parameters $\eta_1=\eta_{\alpha_1}$ and $\eta_n=\eta_{\alpha_n}$. Let $S_0\subset\Delta$ and let $\mathcal D_0$, $W_+$, $\mathcal D$ and $I_{\mathcal D}$ be defined in the same way as in Section~\ref{sec3} for the case $R=C_n$. Let $S_0=\bigsqcup_{\ell=1}^m S_\ell$ and $R_0=\bigsqcup_{\ell=1}^m R_\ell$ be the decompositions defined in Subsection~\ref{sec32}. We denote $k_\ell=|S_\ell|$.

There are two possible cases: $\alpha_n\notin S_0$ and $\alpha_n\in S_0$. In the former case the components $R_\ell\cong A_{k_\ell}$ for all $\ell=1,\ldots,m$. In the latter case one of the components is of type $C$. We suppose $\alpha_n\in S_m$ so that $R_\ell\cong A_{k_\ell}$ for all $\ell=1,\ldots,m-1$ and $R_m\cong C_{k_m}$.

\begin{Lem} \label{Lem_inv_s0DC}
 Assume that
\begin{align} \label{hellA}
 h^\ell_\eta=\hbar &&&\text{\normalfont whenever $\alpha_n\notin R_\ell$}.
\end{align}
Let $w\in W_+$ be such that $w\theta_\ell\ne\theta$ for all $\ell=1,\ldots,m$.
Then $s_0 H_\eta(wS_0)\subset{\cal D}$.
\end{Lem}

\noindent{\bf Proof.} If $\alpha_n\notin S_0$ then the assumptions of Lemma~\ref{Lem_inv_s0D} hold and, therefore, $s_0 H_\eta(wS_0)\subset{\cal D}$. If $\alpha_n\in S_0$ then the statement can be proved in the same way as Lemma~\ref{Lem_inv_s0D}. Indeed, in the proof of Lemma~\ref{Lem_inv_s0D} the condition $h^\ell_\eta=\hbar$ was used only if $(w\theta_\ell,\theta^\vee)=1$. Thus it is sufficient to establish that $(w\theta_m,\theta^\vee)=0$. Since $w S_m$ consists of positive roots only this set has the form
\begin{align}
 wS_m&=\{e_{i_1}-e_{i_2},e_{i_2}-e_{i_3},\ldots,e_{i_{k_m-1}}-e_{i_{k_m}},2e_{i_{k_m}}\},\label{wSmCC}
\end{align}
for some $1\le i_1<i_2<\ldots<i_{k_m}\le n$. Since $w\theta_m=2e_{i_1}$ the condition $w\theta_m\ne\theta=2e_1$ implies $i_1\ne1$ and, hence, $(w\theta_m,\theta^\vee)=(2e_{i_1},e_1)=0$. \qed

\begin{Th} \label{ThInvCC}
 Suppose the subset $S_0\subset\Delta$ and the function $\eta$ are such that the condition~\eqref{hellA} is satisfied. Suppose the parameters of DAHA $\HH_{t,u}$ are such that $t_1^2=e^{\eta_1}$ if $n>1$ and $\epsilon_nt_nu_n^{\epsilon_n}=e^{\eta_n/2}$ for some $\epsilon_n=\pm1$. Suppose also that if $\alpha_n\in S_0$ then
\begin{align} \label{tnt0unu0}
 \epsilon_0\epsilon_n t_nt_0u_n^{\epsilon_n}u_0^{\epsilon_0}t_1^{2(\tilde k-1)}q^{-1}=1
\end{align}
for some $\epsilon_0=\pm1$, where $\tilde k=|S_m|$ is the rank of the irreducible component of type $C$.
Then the ideal $I_{\cal D}\subset{\cal F}$ is invariant under the action of $\HH_{t,u}$ in the Noumi representation~\eqref{ReprDAHACC}.
\end{Th}

\noindent{\bf Proof.} Let $f\in I_{\mathcal D}$ and $x\in w\mathcal D_0$ for some $w\in W_+$.
The operators $T_1,\ldots,T_{n-1}$ have the same form as the operators $T_1,\ldots,T_{n-1}$ of Section~\ref{sec3} in the case $R=C_n$, and the proof of invariance under their action follows by Proposition~\ref{propInvTi}.

To prove the invariance under $T_n$ note that if $\alpha_n\notin wS_0$ then we have $f(x)=f(s_n x)=0$ and, hence, $T_nf(x)=0$ (see the proof of Proposition~\ref{propInvTi}). If $\alpha_n\in wS_0$ then $1-t_nu_ne^{-x_n}=0$ or $1+t_nu_n^{-1}e^{-x_n}=0$ hence $T_nf(x)=0$.

To prove the invariance under $T_0$ suppose first  that $w\theta_\ell=\theta$ for some $\ell$. This is possible only if $\alpha_n\in S_0$ and $\ell=m$ since $\alpha_n\in S_m$. In this case we obtain $-\alpha_0(x)=(w\theta_m,x)-\hbar=h^m_\eta-\eta_{\theta}-\hbar=2(\tilde k-1)\eta_1+\eta_n-\hbar$. Then the condition~\eqref{tnt0unu0} implies that $\epsilon_0t_0u_0^{\epsilon_0}e^{-\alpha_0(x)/2}=1$ for some $\epsilon_0=\pm1$. Now $T_0f(x)=0$ since $t_0u_0e^{-\alpha_0(x)/2}=1$ or $t_0u_0^{-1}e^{-\alpha_0(x)/2}=-1$. In the case when $w\theta_\ell\ne\theta$ for all $\ell=1,\ldots,m$ the plane $s_0(w\mathcal D_0)=s_0H_{\eta}(wS_0)$ is contained in $\mathcal D$ due to Lemma~\ref{Lem_inv_s0DC} and hence $f(s_0x)=0$. As in the proof of Proposition~\ref{propInvT0} it follows that $T_0f(x)=0$ for all $x\in w\mathcal D_0$. \qed

\begin{Rem} \normalfont
 The condition~\eqref{tnt0unu0} can be rewritten in the form $\epsilon_0t_0u_0^{\epsilon_0}=e^{\eta_0/2}$ where $\eta_0\in\mathbb C$ is such that $2(\tilde k-1)\eta_1+\eta_n+\eta_0=\hbar$. Thus the number $\hat h^m_\eta=2(\tilde k-1)\eta_1+\eta_n+\eta_0$ can be interpreted as a generalized Coxeter number for the affine root system $\wh R_m\subset\SS$ of type $(C_{\tilde k}^\vee,C_{\tilde k})$.
\end{Rem}

\begin{Rem} \normalfont
 Submodules in the polynomial representation of DAHA of type $(C_n^\vee,C_n)$ were investigated by Kasatani~\cite{KasCC}. In the case $\alpha_n\notin S_0$ the invariant ideals $I_{\cal D}$ correspond to the Kasatani's submodules $I_m^{(k,r)}$ with $r=2$ (see~\cite[Proposition~4.13]{KasCC}). In the case $\alpha_n\in S_0$, $m=1$ the ideals $I_{\cal D}$ correspond to the submodules $V$ with $r=2$ constructed in the paper~\cite[Subsection~4.8]{KasCC}, the conditions~\eqref{tnt0unu0} correspond to the conditions~(49), (50) of~\cite{KasCC} with $r=2$.
\end{Rem}

Theorem~\ref{ThInvCC} guarantees that the arguments of Subsection~\ref{sec41} remain valid in the $(C_n^\vee,C_n)$ case. Thus we obtain the commutative algebra $\overline{\mathfrak R}$ spanned by the operators $\overline{M_\lambda}$, $\lambda\in Q^\vee$, defined by the formula~\eqref{dbar}.

To establish the integrability we use the method of Subsection~\ref{sec42} where the lattice $P^\vee$ is replaced by $Q^\vee$. Note that $\nu_i=e_1+e_2+\ldots+e_i$, $i=1,\ldots,n$, form the basis of $Q^\vee$. As before we denote $J=\{j=1,\ldots,n \mid \alpha_j\notin S_0\}=\{j_1,\ldots,j_{\bar n}\}$, where $\bar n=\dim\mathcal D_0$.

\begin{Th}
  The algebra $\overline{\mathfrak R}$ contains $\bar n$ algebraically independent elements, namely, the operators $\overline{M_{-\nu_{j_1}}},\ldots,\overline{M_{-\nu_{j_{\bar n}}}}$.
\end{Th}

Indeed, the operators $M_\lambda$ in this case have the form~\eqref{MlambdaY} with
$$g^\lambda(x)=\prod_{\alpha\in\wh R_\lambda}\frac{(1-t_\alpha u_\alpha e^{\alpha(x)/2})(1+t_\alpha u_\alpha^{-1} e^{\alpha(x)/2})}{1-e^{\alpha(x)}}.$$
If the function $g^{-\nu_j}(x)$ is identically zero on $\mathcal D_0$ for some $j\in J$ then the function $1-t_\alpha u_\alpha e^{\alpha(x)/2}$ or $1+t_\alpha u_\alpha^{-1} e^{\alpha(x)/2}$ is identically zero on $\mathcal D_0$ for some $\alpha=[\alpha^{(1)},\alpha^{(0)}]\in\wh R_{-\nu_j}$. It is possible only if $\alpha$ is a linear combination of the elements of $S_0$ and hence $(\alpha^{(1)},\nu_j)=0$. From the definition of $\wh R_\lambda$ it follows that $(\alpha^{(1)},\nu_j)<0$, hence $g^{-\nu_j}\big|_{\mathcal D_0}\not\equiv0$, and the operators $M_{-\nu_j}$, $j\in J$, are algebraically independent (see the proof of Theorem~\ref{Th43} for details).

The algebra $\overline{\mathfrak R}$ contains {\it the generalized Koornwinder operator}. Let us introduce this operator explicitly. It acts on functions of two sets of variables $x_1,\ldots, x_{N_1}$, $y_1,\ldots, y_{N_2}$ and it depends on six independent parameters: two shift parameters $\hbar$, $\xi$ and the parameters $a,b,c,d$. It is defined as
\begin{align} \label{MacKoorGeneralized}
{\cal M}_{N_1,N_2}=
 \sum_{i=1\atop\epsilon=\pm1}^{N_1} A_i^\epsilon(x,y)
\Big(\T^{-\epsilon\hbar}_{x_i}-1\Big)
+\frac{1-q^{-2}}{1-s^{-2}}\sum_{\ell=1\atop\epsilon=\pm1}^{N_2} B_\ell^\epsilon(x,y) \Big(\T^{-\epsilon\xi}_{y_\ell}-1\Big),
\end{align}
where $q = e^{\hbar/2}$, $s=e^{\xi/2}$ and
\begin{multline}
A_i^\epsilon(x,y)=\prod_{j=1\atop j\ne i}^{N_1} g(\epsilon x_i-x_j; s,1)g(\epsilon x_i+x_j; s,1)\times \\
 \prod_{\ell'=1}^{N_2} g(\epsilon x_i-y_{\ell'};q^{1/2}s^{1/2},q^{-1/2}s^{1/2})g(\epsilon x_i+y_{\ell'};q^{1/2}s^{1/2},q^{-1/2}s^{1/2}) \times \\
\frac{(1-a e^{\epsilon x_i})(1-b e^{\epsilon x_i})(1-c e^{\epsilon x_i})(1-d e^{\epsilon x_i})}
{(1-q^{-2}e^{2\epsilon x_i})(1-e^{2\epsilon x_i})},
\end{multline}
\begin{multline}
B_\ell^\epsilon(x,y)=\prod_{j=1}^{N_1} g(\epsilon y_\ell-x_j;q^{1/2}s^{1/2},q^{1/2}s^{-1/2})g(\epsilon y_\ell+x_j;q^{1/2}s^{1/2},q^{1/2}s^{-1/2}) \times \\
 \prod_{\ell'=1\atop\ell'\ne\ell}^{N_2} g(\epsilon y_\ell-y_{\ell'};q,1)g(\epsilon y_\ell+y_{\ell'};q,1) \times \\
\frac{(1-a q s^{-1}e^{\epsilon y_\ell})(1-b q s^{-1}e^{\epsilon y_\ell})(1-c q s^{-1}e^{\epsilon y_\ell})(1-d q s^{-1}e^{\epsilon y_\ell})}
{(1-s^{-2}e^{2\epsilon y_\ell})(1-e^{2\epsilon y_\ell})}.
\end{multline}
For $N_2=0$ this operator coincides with the Koornwinder operator~\eqref{MacKoor} up to a constant and a factor (where $a,b,c,d,\xi$ are expressed through the parameters of DAHA $\HH_{t,u}$).

We note that this operator is invariant (up to a factor) under the involution
\begin{align}
 &N_1\leftrightarrow N_2, &&x\leftrightarrow y, && \xi\leftrightarrow\hbar, && s\leftrightarrow q, && a\mapsto a q s^{-1},
 && b\mapsto b q s^{-1},&& c\mapsto c q s^{-1},&& d\mapsto d q s^{-1}.
\end{align}

Now we explain how this operator arises by the restriction of the Koornwinder operator  $M_{-e_1}$.
Consider first the case $\alpha_n\notin S_0$. In this case we assume that all the subsets $S_\ell$ coincide with~\eqref{SellBn} and $t_1=e^{\eta_1/2}$. Then $N_1=n-mk$, $N_2=m$, $\hbar=k\eta_1$, $q=t_1^k$. We introduce orthogonal coordinates $x_1,\ldots, x_{N_1}$, $y_1,\ldots, y_{N_2}$ on the corresponding plane ${\cal D}_0$ by~\eqref{ylBn}.

\begin{Prop}\label{Koornwinder-generalized-1}
In the case $\alpha_n\notin S_0$
we have $\overline{M_{-e_1}}= (t_0 t_n)^{-1}{\cal M}_{N_1,N_2} + \kappa_t$, where ${\cal M}_{N_1,N_2}$ is the generalized Koornwinder operator~\eqref{MacKoorGeneralized} with $N_1=n-mk$, $N_2=m$ and the parameters
\begin{align}
s&=t_1, &\xi &= \eta_1,& a&= t_0u_0q^{-1},& b&= - t_0u_0^{-1}q^{-1},& c&= t_n u_n,& d &= -t_n u_n^{-1} \notag
\end{align}
satisfying the condition $\hbar/\eta_1 = k \in \mathbb{Z}_{\ge 2}$.
\end{Prop}

In the case $\alpha_n\in S_0$ the subsets $S_\ell$ with $\ell=1,\ldots,m-1$ coincide with~\eqref{SellBn} and
\begin{align}
 S_m&=\{e_{n-\tilde k+1}-e_{n-\tilde k+2},e_{n-\tilde k+2}-e_{n-\tilde k+3},\ldots,e_{n-1}-e_n,2e_n\}, \label{SmCC}
\end{align}
where $\tilde k=k_m=|S_m|$; $t_1=e^{\eta_1/2}$, $\epsilon_n t_nu_n^{\epsilon_n}=e^{\eta_n/2}$. Then
\begin{align}
 N_1=n-(m-1)k-\tilde k, &&&N_2=m-1, \label{N1N2CCC}
\end{align}
and the orthogonal coordinates $x_1,\ldots, x_{N_1}$, $y_1,\ldots, y_{N_2}$ on ${\cal D}_0$ are given by~\eqref{ylBn}

\begin{Prop}\label{Koornwinder-generalized-2}
In the case $\alpha_n\in S_0$ we have $\overline{M_{-e_1}}= (t_0 t_n t_1^{2\tilde k})^{-1}{\cal M}_{N_1,N_2} + \kappa_t$, where ${\cal M}_{N_1,N_2}$ is the generalized Koornwinder operator~\eqref{MacKoorGeneralized} with $N_1$, $N_2$ given by~\eqref{N1N2CCC} and the parameters
\begin{align}
&s=t_1, &&\xi = \eta_1, &&a= \epsilon_n t_1^2 t_n^{-1} u_n^{-\epsilon_n}, &&b= - \epsilon_n t_n u_n^{-\epsilon_n}, &&c= \epsilon_0 q t_1^2 t_0^{-1} u_0^{-\epsilon_0}, &&d = -q^{-1} \epsilon_0 t_0 u_0^{-\epsilon_0} \notag
\end{align}
where the condition~\eqref{tnt0unu0} holds for $\tilde k\in\mathbb Z_{\ge1}$, and if $m>1$ then $\hbar/\eta_1 = k \in \mathbb{Z}_{\ge 2}$.
\end{Prop}

The propositions can be established by direct calculations. In both cases the obtained operators $\overline{M_{-e_1}}$ depend on six independent parameters. Namely, in the case $\alpha_n\notin S_0$ we have five continuous parameters $\xi,a,b,c,d$ and one discrete parameter $k$ (where $\hbar=k\xi$), and in the case $\alpha_n\in S_0$ we have four continuous parameters $\xi,b,c,d$ and two discrete ones $k,\tilde k$ ($\hbar=k\xi$, $a=c^{-1}e^{(\tilde k+1)\xi}$).



\begin{Rem} \normalfont The operator~\eqref{MacKoorGeneralized} generalizes {\it the deformed rational Koornwinder operator}~\cite[formula~(94)]{SVBC}. Namely the latter can be obtained by the substitution $x_i\to x_i/\omega$, $y_i\to y_i/\omega$, $\xi\to\xi/\omega$, $\hbar\to\hbar/\omega$, $a=e^{a'/\omega}$, $b=-e^{b'/\omega}$, $c=e^{c'/\omega}$, $d=-e^{d'/\omega}$ to~\eqref{MacKoorGeneralized} and by taking the limit $\omega\to\infty$.
\end{Rem}


It follows from above that the operator ${\cal M}_{N_1,N_2}$ comes together with the complete algebra of commuting difference operators. We specify these operators explicitly and establish commutativity (without extra integrality restriction for a combination of parameters).
We start with the explicit formula for the higher Koornwinder operators obtained by van Diejen~\cite{vD96}:
\begin{align}
 H_r=&\sum_{J_1\subset\{1,\ldots,n\},\; |J_1|=r\atop\epsilon_i=\pm1,\; i\in J_1}\;
\prod_{i\in J_1}G(e^{\epsilon_i x_i})
\sum_{J_2\subset J_1}(-1)^{r-|J_2|}\; \mathcal U_1\,\mathcal U_2\;\prod_{i\in J_2}\T^{-\epsilon_i\hbar}_{x_i}, \label{Koornr}
\end{align}
where $r=1,\ldots,n$,
\begin{align}
 G(z)=\frac{(1-t_0u_0q^{-1}z)(1+t_0u_0^{-1}q^{-1}z)(1-t_nu_nz)(1+t_nu_n^{-1}z)}
{(1-q^{-2}z^2)(1-z^2)},
\end{align}
\begin{multline} \label{Ubg}
 \mathcal U_b=\prod_{i,j\in J_b\backslash J_{b+1}\atop i<j} g(\epsilon_i x_i+\epsilon_j x_j;t_1,1)g(\epsilon_i x_i+\epsilon_j x_j;q^{-1}t_1^{2b-3},q^{-1}) \times \\
 \prod_{i\in J_b\backslash J_{b+1}\atop j\notin J_b} g(\epsilon_i x_i+ x_j;t_1,1)g(\epsilon_i x_i- x_j;t_1,1),
\end{multline}
where $b=1,2$, $J_3=\emptyset$ and $q=e^{\hbar/2}$.
These operators commute with each other, and for $r=1$ we have $H_1=t_0t_n(M_{-e_1}-\kappa_t)$, where $M_{-e_1}$ is the Koornwinder operator~\eqref{MacKoor}.

Note that the operators have the form $H_r=g^r_{-\nu_r}(x)\tau(-\nu_r)+\sum_{\lambda>-\nu_r}g^r_\lambda(x)\tau(\lambda)$, where $\nu_r=e_1+\ldots+e_r$ and $g^r_\lambda(x)$ are meromorphic functions such that $g^r_{-\nu_r}(x)\not\equiv0$. Hence the operators~\eqref{Koornr} with $r=1,\ldots,n$ are algebraically independent (c.f. Subsection~\ref{sec22}).

The operators $H_r$ are $W$-invariant. Since for generic parameters the centralizer of $M_{-e_1}$ in the algebra of $W$-invariant difference operators (with suitable coefficients) coincides with the algebra $\mathfrak R$ (see~\cite[Theorem~3.15]{LS}) we obtain $H_r\in\mathfrak R$. Hence $H_r\in\mathfrak R$ for any values of the parameters, and the restrictions of the operators $H_r$ onto the corresponding plane $\mathcal D_0$ can be defined and they commute with each other.

Consider restriction of the operators~\eqref{Koornr} in the case $\alpha_n\notin S_0$. The structure of the functions~\eqref{Ubg} implies that if $\mathcal U_b$ is not zero on the plane $\mathcal D_0$ then $J_b$ have the form
\begin{align} \label{JbIb}
 J_b=I_b\cup\bigcup\nolimits_{\ell=1}^{N_2}\{N_1+k\ell-s\mid
 0\le s\le m^-_{\ell b}-1 \text{ or } k-m^+_{\ell b}\le s\le k-1\}, &&&b=1,2,
\end{align}
where $I_2\subset I_1\subset\{1,\ldots,N_1\}$ and $m^\pm_{\ell b}$ are non-negative integers such that $m^{\pm}_{\ell2}\le m^{\pm}_{\ell1}$, $m^-_{\ell 1}+m^+_{\ell 1}\le k$ and either $m^-_{\ell2}$ or $m^+_{\ell2}$ vanishes for each $\ell$. Moreover, the signs $\epsilon_i$ are given by
\begin{align} \label{epsm}
 \epsilon_{N_1+k\ell-s}=\begin{cases} -1, \text{ if $0\le s\le m^-_{\ell 1}-1$} \\ +1, \text{ if $k-m^+_{\ell 1}\le s\le k-1$} \end{cases}.
\end{align}
We recall the coordinates $x_1,\ldots,x_{N_1}$, $y_1,\ldots,y_{N_2}$ defined by~\eqref{ylBn} on the plane $\mathcal D_0$.

\begin{Prop} \label{Prophrrestr}
The restriction of the operator $H_r$ takes the form
\begin{multline}\label{hrrestr}
 \overline{H_r}=\sum_{\substack{I_1\subset\{1,\ldots,N_1\},\\ \epsilon_i=\pm1,\; i\in I_1 \\ m^{\pm}_{\ell1}\in\mathbb Z_{\ge0},\;1\le\ell\le N_2 \\ |I_1|+\sum_{\ell=1}^{N_2}(m^+_{\ell1}+m^-_{\ell1})=r}}
\prod_{i\in I_1}G(e^{\epsilon_i x_i}) \prod_{\ell=1\atop\epsilon=\pm1}^{N_2} \prod_{p=0}^{m_{\ell1}^\epsilon-1} G(qt_1^{-1-2p}e^{\epsilon y_\ell})\times \\
 \sum_{\substack{I_2\subset I_1,\;m^{\pm}_{\ell2}\in\mathbb Z_{\ge0} \\ m^\pm_{\ell2}\le m^\pm_{\ell 1},\; m_{\ell2}^+m_{\ell2}^-=0}}(-1)^{r-|I_2|-\sum_{\ell=1}^{N_2}(m^+_{\ell2}+m^-_{\ell2})}\;\; \overline{\mathcal U}_1\overline{\mathcal U}_2\;
\prod_{i\in I_2}\T^{-\epsilon_i\hbar}_{x_i}
\prod_{\ell=1}^{N_2}\T^{-(m_{\ell2}^+-m_{\ell2}^-)\eta_1}_{y_\ell},
\end{multline}
where
\begin{multline}
 \overline{\mathcal U}_b=\prod_{\ell=1}^{N_2}\prod_{p=m^+_{\ell,b+1}}^{m^+_{\ell b}-1}
g(0;qt_1^{-p-m^-_{\ell,b+1}},t_1^{-p+m^+_{\ell b}})
g(0;t_1^{-p-1-m^-_{\ell b}},t_1^{-p-1-m^-_{\ell,b+1}}) \times \\
\prod_{\ell=1}^{N_2}\prod_{p=m^-_{\ell,b+1}}^{m^-_{\ell b}-1}
g(0;qt_1^{-p-m^+_{\ell b}},t_1^{-p+m^-_{\ell b}})
 \prod_{i,j\in I_b\backslash I_{b+1}\atop i<j} g(\epsilon_i x_i+\epsilon_j x_j;t_1,1)g(\epsilon_i x_i+\epsilon_j x_j;q^{-1}t_1^{2b-3},q^{-1}) \times \\
 \prod_{{i\in I_b\backslash I_{b+1}\atop\epsilon=\pm1}\atop j\notin I_b} g(\epsilon_i x_i+\epsilon x_j;t_1,1)\;\;
\prod_{{\ell=1 \atop \epsilon,\epsilon'=\pm1}\atop j\notin I_b}^{N_2} g(\epsilon y_\ell+\epsilon'x_j;q^{\frac12}t_1^{\frac12-m^\epsilon_{\ell,b+1}},q^{\frac12}t_1^{\frac12-m^\epsilon_{\ell b}})   \times \notag
\end{multline} \\[-30pt]
\begin{multline}
\prod_{{\ell=1\atop\epsilon=\pm1}\atop i\in I_b\backslash I_{b+1}}^{N_2} g(\epsilon_i x_i+\epsilon y_\ell;q^{\frac12}t_1^{\frac12-m^\epsilon_{\ell,b+1}},q^{-\frac12}t_1^{\frac12+m^{-\epsilon}_{\ell b}})g(\epsilon_i x_i+\epsilon y_\ell;q^{-\frac12}t_1^{\frac{2b-3}2-(2-b) m^\epsilon_{\ell1}},q^{-\frac12}t_1^{\frac{2b-3}2-m^\epsilon_{\ell2}}) \times \\
\prod_{1\le\ell<\ell'\le N_2 \atop \epsilon,\epsilon'=\pm1}\;\prod_{p=m^\epsilon_{\ell,b+1}}^{m^\epsilon_{\ell b}-1}
\Big(g(\epsilon y_\ell+\epsilon' y_{\ell'};qt_1^{-p-m^{\epsilon'}_{\ell',b+1}},t_1^{-p+m^{-\epsilon'}_{\ell' b}})\times \\
g(\epsilon y_\ell+\epsilon' y_{\ell'};t_1^{-p-(2-b)(1+m^{\epsilon'}_{\ell'1})},t_1^{-p-(2-b)-m^{\epsilon'}_{\ell'2}})\Big)  \times \notag
\end{multline} \\[-30pt]
\begin{multline} \label{urrestr}
\prod_{1\le\ell'<\ell\le N_2 \atop \epsilon,\epsilon'=\pm1}\;\prod_{p=m^\epsilon_{\ell,b+1}}^{m^\epsilon_{\ell b}-1} g(\epsilon y_\ell+\epsilon'y_{\ell'};qt_1^{-p-m^{\epsilon'}_{\ell' b}},t_1^{-p+m^{-\epsilon'}_{\ell' b}}) \times \\
\prod_{\ell=1\atop\epsilon=\pm1}^{N_2}\prod_{p=m^\epsilon_{\ell,b+1}}^{m^\epsilon_{\ell b}-1}  g(2\epsilon y_\ell;qt_1^{-1-2p},t_1^{-p+m^{-\epsilon}_{\ell b}})
g(2\epsilon y_\ell;t_1^{-1-bp-(2-b)m^\epsilon_{\ell1}},t_1^{-2(2-b)-(3-b)p-(b-1) m^\epsilon_{\ell2}}),
\end{multline}
with $b=1,2$, $I_3=\emptyset$, $m^\pm_{\ell3}=0$ and $q=e^{\hbar/2}$, $t_1=e^{\eta_1/2}$.
\end{Prop}

Note that we do not impose the condition $m^-_{\ell1}+m^+_{\ell1}\le k$ in the sum in~\eqref{hrrestr} since the product $\overline{\mathcal U}_1\,\overline{\mathcal U}_2$ vanishes if $m^-_{\ell1}+m^+_{\ell1}>k$, so the operators $\overline{H_r}$ do not have explicit dependence on $k$. \\

\noindent{\bf Proof of Proposition~\ref{Prophrrestr}.} The structure of the operator~\eqref{hrrestr} follows from the formulae~\eqref{JbIb}, \eqref{epsm}.
Let $i=N_1+k\ell-s\in J_{\ell b}$ for some $\ell=1,\ldots, N_2$, $s=0,\ldots,k-1$, $b=1,2$. Then the formulae~\eqref{JbIb}, \eqref{epsm} imply that on $\mathcal D_0$ one has $\epsilon x_i=\epsilon y_\ell+\big(\frac{k-1}2-p\big)\eta_1$  for some $p=0,\ldots,m^\epsilon_{\ell b}-1$ with $\epsilon=\epsilon_i$. Hence the restriction of $\prod_{i\in J_1}G(e^{\epsilon_i x_i})$ takes the form
\begin{align}
\prod_{i\in I_1}G(e^{\epsilon_i x_i}) \prod_{\ell=1\atop\epsilon=\pm1}^{N_2} \prod_{p=0}^{m_{\ell1}^\epsilon-1} G(qt_1^{-1-2p}e^{\epsilon y_\ell}). \notag
\end{align}

Consider the restriction of the first product in the expression~\eqref{Ubg} which is taken over $i,j\in J_b\backslash J_{b+1}$. Its part corresponding to the case $i<j\le N_1$ is equal to
\begin{align}
\prod_{i,j\in I_b\backslash I_{b+1}\atop i<j} g(\epsilon_i x_i+\epsilon_j x_j;t_1,1)g(\epsilon_i x_i+\epsilon_j x_j;q^{-1}t_1^{2b-3},q^{-1}). \notag
\end{align}
For $i\le N_1<j$ we obtain
\begin{multline}
\prod_{i\in I_b\backslash I_{b+1}}\prod_{\ell=1\atop\epsilon=\pm1}^{N_2} \prod_{p=m^\epsilon_{\ell,b+1}}^{m^\epsilon_{\ell b}-1}g(\epsilon_i x_i+\epsilon y_\ell;q^{\frac12}t_1^{\frac12-p},q^{\frac12}t_1^{-\frac12-p})g(\epsilon_i x_i+\epsilon y_\ell;q^{-\frac12}t_1^{-\frac12-p+2b-3},q^{-\frac12}t_1^{-\frac12-p})= \\
=\prod_{{\ell=1\atop\epsilon=\pm1}\atop i\in I_b\backslash I_{b+1}}^{N_2} g(\epsilon_i x_i+\epsilon y_\ell;q^{\frac12}t_1^{\frac12-m^\epsilon_{\ell,b+1}},q^{\frac12}t_1^{\frac12-m^\epsilon_{\ell b}})g(\epsilon_i x_i+\epsilon y_\ell;q^{-\frac12}t_1^{\frac{2b-3}2-(2-b) m^\epsilon_{\ell1}},q^{-\frac12}t_1^{\frac{2b-3}2-m^\epsilon_{\ell2}}). \label{pr56gxy1}
\end{multline}
Here we used the formulae
\begin{align*}
 \prod_{p=m}^{M}g(z;at_1^{1-p},at_1^{-p})=g(z;at_1^{1-m},at_1^{-M}), &&&
 \prod_{p=m}^{M}g(z;at_1^{-1-p},at_1^{-p})=g(z;at_1^{-1-M},at_1^{-m}),
\end{align*}
where $m,M\in\mathbb Z$ are such that $m\le M+1$, $a,z\in\mathbb C$ (the second formula is used in the case $b=1$ only).
If $N_1<i<j$ then $i=N_1+k\ell-s$ and $j=N_1+k\ell'-s'$, where $1\le\ell,\ell'\le N_2$, $0\le s,s'\le k-1$. If $\ell\ne\ell'$ then we obtain
\begin{multline}
\prod_{1\le\ell<\ell'\le N_2 \atop \epsilon,\epsilon'=\pm1}\prod_{p=m^\epsilon_{\ell,b+1}}^{m^\epsilon_{\ell b}-1}\prod_{p'=m^{\epsilon'}_{\ell',b+1}}^{m^{\epsilon'}_{\ell' b}-1} g(\epsilon y_\ell+\epsilon' y_\ell;qt_1^{-p-p'},qt_1^{-1-p-p'})g(\epsilon y_\ell+\epsilon' y_\ell;t_1^{-1-p-p'+2b-3},t_1^{-1-p-p'})= \\
=\prod_{1\le\ell<\ell'\le N_2 \atop \epsilon,\epsilon'=\pm1}\prod_{p=m^\epsilon_{\ell,b+1}}^{m^\epsilon_{\ell b}-1}
\Big(g(\epsilon y_\ell+\epsilon' y_{\ell'};qt_1^{-p-m^{\epsilon'}_{\ell',b+1}},qt_1^{-p-m^{\epsilon'}_{\ell' b}})\times \\
g(\epsilon y_\ell+\epsilon' y_{\ell'};
t_1^{-p-(2-b)(1+m^{\epsilon'}_{\ell'1})},t_1^{-p-(2-b)-m^{\epsilon'}_{\ell'2}})\Big), \label{pr56gellell}
\end{multline}
while if $\ell=\ell'$ we derive the following two factors. The factor corresponding to $\epsilon_i=-\epsilon_j=+1$ is equal to
\begin{multline}
\prod_{\ell=1}^{N_2} \prod_{p=m^+_{\ell,b+1}}^{m^+_{\ell b}-1}\prod_{p'=m^-_{\ell,b+1}}^{m^-_{\ell b}-1} g(0;qt_1^{-p-p'},qt_1^{-1-p-p'})g(0;t_1^{-1-p-p'+2b-3},t_1^{-1-p-p'})= \\
\prod_{\ell=1}^{N_2}\prod_{p=m^+_{\ell,b+1}}^{m^+_{\ell b}-1}
g(0;qt_1^{-p-m^-_{\ell,b+1}},qt_1^{-p-m^-_{\ell b}})
g(0;t_1^{-p-1-m^-_{\ell b}},t_1^{-p-1-m^-_{\ell,b+1}}). \label{pr56g02s1}
\end{multline}
(In the case $b=2$ both sides of the equality \mref{pr56g02s1} equal 1).
The factor corresponding to $\epsilon_i=\epsilon_j=\epsilon=\pm1$ is equal to
\begin{multline}
\prod_{\ell=1\atop\epsilon=\pm1}^{N_2} \prod_{p=m^\epsilon_{\ell,b+1}}^{m^\epsilon_{\ell b}-1}\prod_{p'=p+1}^{m^\epsilon_{\ell b}-1} g(2\epsilon y_\ell;qt_1^{-p-p'},qt_1^{-1-p-p'})g(2\epsilon y_\ell;t_1^{-1-p-p'+2b-3},t_1^{-1-p-p'})= \\
\prod_{\ell=1\atop\epsilon=\pm1}^{N_2}\prod_{p=m^\epsilon_{\ell,b+1}}^{m^\epsilon_{\ell b}-1}  g(2\epsilon y_\ell;qt_1^{-1-2p},qt_1^{-p-m^\epsilon_{\ell b}})
g(2\epsilon y_\ell;t_1^{-1-bp-(2-b)m^\epsilon_{\ell1}},t_1^{-2(2-b)-(3-b)p-(b-1) m^\epsilon_{\ell2}}). \label{pr56g02s2}
\end{multline}

Now consider the restriction of the second product in the expression~\eqref{Ubg} which is taken over $i\in J_b\backslash J_{b+1}$, $j\notin J_b$. If $1\le i,j\le N_1$ we obtain
\begin{align}
 \prod_{i\in I_b\backslash I_{b+1}\atop j\notin I_b} g(\epsilon_i x_i+ x_j;t_1,1)g(\epsilon_i x_i- x_j;t_1,1). \notag
\end{align}
For $j\le N_1<i$ we obtain
\begin{align}
\prod_{\ell=1 \atop \epsilon,\epsilon'=\pm1}^{N_2}\prod_{j=1\atop j\notin I_b}^{N_1} g(\epsilon y_\ell+\epsilon'x_j;q^{\frac12}t_1^{\frac12-m^\epsilon_{\ell,b+1}},q^{\frac12}t_1^{\frac12-m^\epsilon_{\ell b}}). \notag
\end{align}
If $i\le N_1<j$ then $\epsilon x_j=\epsilon y_\ell+\big(\frac{k-1}2-p\big)\eta_1$ with $p=m^\epsilon_{\ell b},\ldots,k-1-m^{-\epsilon}_{\ell b}$ for any $\epsilon=\pm1$ (it follows from $j\notin J_b$), so we obtain
\begin{align}
 \prod_{{\ell=1\atop\epsilon=\pm1}\atop i\in I_b\backslash I_{b+1}}^{N_2} \prod_{p=m^\epsilon_{\ell b}}^{k-1-m^{-\epsilon}_{\ell b}}
 g(\epsilon_i x_i+\epsilon y_\ell;q^{\frac12}t_1^{\frac12-p},q^{\frac12}t_1^{-\frac12-p})
=\prod_{{\ell=1\atop\epsilon=\pm1}\atop i\in I_b\backslash I_{b+1}}^{N_2} g(\epsilon_i x_i+\epsilon y_\ell;q^{\frac12}t_1^{\frac12-m^\epsilon_{\ell b}},q^{-\frac12}t_1^{\frac12+m^{-\epsilon}_{\ell b}}). \label{pr56gxy2}
\end{align}
If $N_1<i<j$ then $i=N_1+k\ell-s$ and $j=N_1+k\ell'-s'$, where $1\le\ell,\ell'\le N_2$, $0\le s,s'\le k-1$. In the case $\ell\ne\ell'$ we obtain
\begin{multline}
\prod_{{\ell,\ell'=1\atop\ell\ne\ell'}\atop\epsilon,\epsilon'=\pm1}^{N_2} \prod_{p=m^\epsilon_{\ell,b+1}}^{m^\epsilon_{\ell b}-1}\prod_{p'=m^{\epsilon'}_{\ell' b}}^{k-1-m^{-\epsilon'}_{\ell' b}}
 g(\epsilon y_\ell+\epsilon' y_{\ell'};qt_1^{-p-p'},qt_1^{-1-p-p'}) = \\
=\prod_{{\ell,\ell'=1\atop\ell\ne\ell'}\atop\epsilon,\epsilon'=\pm1}^{N_2} \prod_{p=m^\epsilon_{\ell,b+1}}^{m^\epsilon_{\ell b}-1} g(\epsilon y_\ell+\epsilon'y_{\ell'};qt_1^{-p-m^{\epsilon'}_{\ell' b}},t_1^{-p+m^{-\epsilon'}_{\ell' b}}), \label{pr56gellellne}
\end{multline}
and for $\ell=\ell'$ we derive
\begin{multline}
\prod_{\ell=1\atop\epsilon,\epsilon'=\pm1}^{N_2} \prod_{p=m^\epsilon_{\ell,b+1}}^{m^\epsilon_{\ell b}-1}\prod_{p'=m^{\epsilon'}_{\ell b}}^{k-1-m^{-\epsilon'}_{\ell b}}
 g(\epsilon y_\ell+\epsilon' y_\ell;qt_1^{-p-p'},qt_1^{-1-p-p'}) = \\
\prod_{\ell=1\atop\epsilon=\pm1}^{N_2} \prod_{p=m^\epsilon_{\ell,b+1}}^{m^\epsilon_{\ell b}-1} g(0;qt_1^{-p-m^{-\epsilon}_{\ell b}},t_1^{-p+m^\epsilon_{\ell b}})\;
\prod_{\ell=1\atop\epsilon=\pm1}^{N_2} \prod_{p=m^\epsilon_{\ell,b+1}}^{m^\epsilon_{\ell b}-1} g(2\epsilon y_\ell;qt_1^{-p-m^\epsilon_{\ell b}},t_1^{-p+m^{-\epsilon}_{\ell b}}), \label{pr56g02d}
\end{multline}
(the first factor in the right hand side corresponds to $\epsilon'=-\epsilon$ and the second one corresponds to $\epsilon'=\epsilon$).
By combining the factors~\eqref{pr56g02s1}, \eqref{pr56g02s2} and~\eqref{pr56g02d} we obtain the beginning and ending of the formula~\eqref{urrestr}. Then combining~\eqref{pr56gellell} and the factors of~\eqref{pr56gellellne} with $\ell<\ell'$ we obtain the corresponding factors of~\eqref{urrestr}. By combining~\eqref{pr56gxy1} with \eqref{pr56gxy2} and taking into account all the remaining factors we obtain the formula~\eqref{urrestr}. \qed

For $r=1$ we obtain the generalized Koornwinder operator $\overline{H_1}={\cal M}_{N_1,N_2}$ with the parameters as in Proposition~\ref{Koornwinder-generalized-1}, so we get a collection of operators commuting with the operator~\eqref{MacKoorGeneralized}. Note that if $\hbar=\eta_1$ the operators $\overline{H_r}$ coincide with the specialization of $H_r$ for $n=N_1+N_2$. Thus we obtain the following theorem.

\begin{Th} \label{PropHOMOcommCC}
 For any values of the parameters $\hbar$, $\eta_1$, $t_0,t_n,u_0,u_n$ $(u_0,u_n\ne0)$ the operators $\overline{H_r}$ given by the formulae~\eqref{hrrestr},\eqref{urrestr}, $r\in\mathbb Z_{>0}$, pairwise commute. The operators $\overline{H_r}$ with $r=1,\ldots,N_1+N_2$ are algebraically independent at generic values of the parameters.
\end{Th}


\section{The case of DAHA of $\pmb{\wh{R}^\vee}$ type}
\label{sec62}

Let $\wh{R}^\vee=\{\alpha^\vee\mid\alpha\in\wh R\}\subset\wh V$, where $R$ is an irreducible reduced root system in $\mathbb R^n\subset V\cong\mathbb C^n$ with set of simple roots $\Delta = \{\a_1,\ldots,\a_n\}$. We have $\wh R^\vee\not\cong\wh R$ only for $R=B_n$, $C_n$, $F_4$ and $G_2$, but we will consider the general $R$.

Let $W$, $\wh W$, $\Pi$, ${\cal O}$ be the same as in Section~\ref{sec2} and let $t\colon\wh R^\vee\to\mathbb C\backslash\{0\}$ be a $\wh W$-invariant function. Put $t_i=t_{\alpha_i^\vee}=t(\a_i^\vee)$ and $\wh{P^\vee}=\{\mu=[\mu^{(1)},\mu^{(0)}]\mid\mu^{(1)}\in P^\vee,\,\mu^{(0)}\in\frac\hbar{\epsilon_0}\mathbb Z\}\subset \wh V$, where $\epsilon_0\in\mathbb Z$ is such that $(P^\vee,P^\vee)\subset\epsilon_0^{-1}\mathbb Z$. Then $P^\vee$ is a sublattice of $\wh{P^\vee}$. The associated DAHA is defined as follows.

\begin{Def}(\cite{CherDAHA}, \cite{Cher92DAHA})
 Double Affine Hecke Algebra of type $\wh R^\vee$ is the algebra $\HH^{\wh R^\vee}_t$ generated by the elements $T_0,T_1,\ldots,T_n$, by the set of pairwise commuting elements $\{X^\mu\}_{\mu\in\hat P}$ and by the group $\Pi=\{\pi_r\mid r\in\mathcal O\}$ with the relations~\eqref{DAHAxx}--\eqref{DAHApt}, \eqref{DAHApx} (with $\wh{P^\vee}$ instead of $\wh P$) and {\normalfont
\begin{align}
 T_i^{-1}X^\mu T_i^{-1}&=X^{\mu-\alpha_i^\vee}, & &\text{if $(\mu,\alpha_i)=1$},& &\mu\in\wh{P^\vee}; \label{DAHAtx1vee} \\
 T_iX^\mu&=X^\mu T_i, & &\text{if $(\mu,\alpha_i)=0$},& &\mu\in\wh{P^\vee}; \label{DAHAtx0vee}
\end{align}}
where $i=0,1,\ldots,n$.
\end{Def}

Consider the space $\check{\cal F}=C^\infty(V)^{2\pi i Q}$ of smooth functions $f$ such that $f(x+2\pi i\alpha_j)=f(x)$ $\forall j=1,\ldots,n$. The representation of DAHA $\HH^{\wh R^\vee}_t$ in this space is defined by the formulae
\begin{align}
 X^\mu&=e^{\mu(x)}, & \pi_r&=\pi_r, & T_i=\frac{t_i-t_i^{-1}}{1-e^{-\alpha^\vee_i(x)}}+\frac{t_i^{-1}-t_ie^{-\alpha^\vee_i(x)}}{1-e^{-\alpha^\vee_i(x)}}\,s_i, \label{ReprDAHAvee}
\end{align}
where $\mu\in\wh{P^\vee}$, $r\in{\cal O}$, $i=0,1,\ldots,n$, so
\begin{align} \label{Tifxvee}
 T_i f(x)=t_i f(x)+\big(t_i^{-1}-t_ie^{-\alpha^\vee_i(x)}\big)\tilde g(x)
\end{align}
where  $f\in\check{\mathcal F}$ and $\tilde g(x)=\frac{f(s_ix)-f(x)}{1-e^{-\alpha^\vee_i(x)}}\in\check{\cal F}$.

Let $m_\lambda(Y)$ be the elements~\eqref{flambdadef} of $\HH^{\wh R^\vee}_t$ in the representation~\eqref{ReprDAHAvee}, where $Y^\lambda$ are defined by the formulae~\eqref{Twdef}, \eqref{Ydef}. Their restrictions to $\check{\mathcal F}^W$ are difference operators $M_\lambda=\Res\big(m_\lambda(Y)\big)$ spanning the commutative algebra $\mathfrak R$. The algebra $\mathfrak R$ is generated by the algebraically independent operators $M_{-b_1},\ldots,M_{-b_n}$, where $b_i$ are the fundamental coweights for $R$. This is the Macdonald-Ruijsenaars system corresponding to the affine root system $\wh R^\vee$. In particular, it always contains the Macdonald operator
\begin{align}
 &M_{-\theta^\vee}=\sum_{\nu\in W\theta^\vee}\Big(\prod_{\alpha\in R\atop (\alpha,\nu)>0}\frac{t_{\alpha^\vee} e^{(\alpha^\vee,x)}-t_{\alpha^\vee}^{-1}}{e^{(\alpha^\vee,x)}-1}\Big)\;\frac{t_0e^{(\nu,x)-\hbar(\theta^\vee,\theta^\vee)/2}-t_0^{-1}}{e^{(\nu,x)-\hbar(\theta^\vee,\theta^\vee)/2}-1}\;\Big(\tau(\nu)-1\Big)+m_{-\theta^\vee,t}\,. \label{RMop_thvee}
\end{align}

Let $S_0\subset\Delta$ and let $\eta\colon\wh R^\vee\to\mathbb C$ be a $\wh W$-invariant function (it is sufficient to define it on the subsystem $R^\vee\subset\wh{R}^\vee$). Define the corresponding ideal $I_{\cal D}\subset{\cal F}$ as follows
\begin{align}
 W_+&=\{w\in W\mid wS_0\subset R_+\}\subset W, \\
 \mathcal D_0&=H\big((S_0^\vee)'\big)=\{x\in V \mid (\alpha^\vee,x)=\eta_{\alpha^\vee},\forall\alpha\in S_0\}, \label{D0defvee} \\[5pt]
 {\cal D}&=W_+\mathcal D_0, \qquad\qquad
 I_{\cal D}=\{f\in{\cal F}\mid f(x)=0, \forall x\in{\cal D}\}. \label{DandIdefvee}
\end{align}

Let $R_0=\bigsqcup_{\ell=1}^m R_\ell$ and $S_0=\bigsqcup_{\ell=1}^m S_\ell$ be the decompositions into irreducible components, where $R_0=\langle S_0\rangle\cap R$. Let $\theta_\ell\in R_\ell$ be the maximal root with respect to the basis $S_\ell$. Let $S_\ell=\{\alpha^\ell_j\mid j=1,\ldots,k_\ell\}$, $\theta_\ell=\sum_{j=1}^{k_\ell} n^\ell_j\alpha^\ell_j$ and $\eta^\ell_j=\eta_{(\alpha^\ell_j)^\vee}$. Define the corresponding generalized Coxeter numbers as
\begin{align} \label{genCoxcheck}
\check h^\ell_\eta=
\frac{(\theta_\ell,\theta_\ell)}2\eta_{\theta^\vee_\ell}+\sum_{j=1}^{k_\ell}\frac{(\alpha_j^\ell,\alpha_j^\ell)}2 n^\ell_j\eta^\ell_j.
\end{align}

\begin{Th} \label{ThInvvee}
 Suppose the subset $S_0\subset\Delta$ and the invariant functions $\eta,t$ satisfy $\check h^\ell_\eta=\hbar$ for all $\ell=1,\ldots,m$ and $t_{\alpha^\vee}^2=e^{\eta_{\alpha^\vee}}$ for all $\alpha\in R$. Then the ideal $I_{\cal D}\subset\check{\cal F}$ is invariant under the action of $\HH^{\wh R^\vee}_t$.
\end{Th}

\noindent{\bf Proof.} First we define the $W$-invariant function $\eta'\colon R\to\mathbb C$ by the formula $\eta'_{\alpha}=\frac{(\alpha,\alpha)}{2}\eta_{\alpha^\vee}$. Note that one can represent~\eqref{D0defvee} as
\begin{align}
 \mathcal D_0=H_{\eta'}(S_0)=\{x\in V \mid (\alpha,x)=\eta'_{\alpha},\forall\alpha\in S_0\}, \label{D0etaprime}
\end{align}
By virtue of~\eqref{D0etaprime} the ideal~\eqref{DandIdefvee} for the function $\eta\colon R^\vee\to\mathbb C$ coincides with the ideal~\eqref{DandIdef} for the function $\eta'$, moreover the conditions $\check h^\ell_\eta=\hbar$ are equivalent to $h^\ell_{\eta'}=\hbar$ where $h^\ell_{\eta'}$ is the generalized Coxeter number for the function $\eta'$. Hence the invariance under the action of $\pi_r$ follows from Proposition~\ref{propInvPi}.

Let $f\in I_{\mathcal D}$ and $x\in w\mathcal D_0$ for some $w\in W_+$.
Consider first the invariance under $T_i$ with $i=1,\ldots,n$. In the case $\alpha_i\in wS_0$ it follows from $f(x)=0$ and $(t^{-1}_i-t_ie^{-\alpha_i^\vee(x)})=0$. In the case $\alpha_i\notin wS_0$ we have $f(x)=f(s_ix)=0$ and the plane $w\mathcal D_0=H_{\eta'}(wS_0)$ does not lie in any plane $H([\alpha_i^\vee,2\pi i k])=H([\alpha_i,\pi i k(\alpha_i,\alpha_i)])$ due to Lemma~\ref{LemSumPhi} (see the proof of Proposition~\ref{propInvTi} for detail).

Consider now the invariance under $T_0$. If $w\theta_\ell=\theta$ for some $\ell$ then $\alpha^\vee_0(x)=\eta_{\theta^\vee}$ and $t_0^{-1}-t_0e^{-\alpha^\vee_0(x)}=0$; so that $T_0f(x)=0$. If $w\theta_\ell\ne\theta$ for all $\ell=1,\ldots,m$ then the plane $s_0\mathcal D_0=s_0H_{\eta'}(wS_0)$ contains in $\mathcal D$ due to Lemma~\ref{Lem_inv_s0D} and the condition $h^\ell_{\eta'}=\hbar$. Hence $f(x)=f(s_0x)=0$ and $T_0f(x)=0$ for $x\in w\mathcal D_0$ such that $e^{\alpha^\vee_0(x)}\ne1$. Thus we are left to prove that $w\mathcal D_0=H_{\eta'}(wS_0)$ does not lie in the affine hyperplane $H(\alpha^\vee_0+[0,2\pi i k])=H(\alpha_0+[0,\pi i k(\theta,\theta)])$. This is done in the same way as in the proof of Proposition~\ref{propInvT0}. \qed \\

Let $\overline{M_\lambda}$ be the restriction of the difference operator $M_\lambda$ to the plane $\mathcal D_0$ defined as in Subsection~\ref{sec41}. Repeating the reasoning of Subsections~\ref{sec41} and~\ref{sec42} we conclude that the operators $\overline{M_\lambda}$ span the commutative algebra $\overline{\mathfrak R}$ containing $\bar n=\dim\mathcal D_0$ algebraically independent elements. This is the {\it generalized Macdonald-Ruijsenaars system for the  $\wh R^\vee$ case}.

\section{Concluding remarks}

We considered submodules in the polynomial representation of DAHA given by vanishing conditions and derived the corresponding generalized Macdonald-Ruijsenaars operators and their quantum integrals acting in the quotient module. In some particular cases the eigenfunctions of these operators can be given by Baker-Akhiezer functions (\cite{Chalykhbisp}, \cite{ChalykhAdv}, \cite{Fbisp})
but more generally the understanding of eigenfunctions is a major challenge in the area. It follows from our construction that the restriction of a Weyl invariant eigenfunction of the higher Macdonald operators $M_\lambda$ onto the plane ${\cal D}_0$ is an eigenfunction of the corresponding operators $\overline{M_\lambda}$. However the question of finding the solutions is non-trivial as the specializations of the parameters so that the ideals $I_{\cal D}$ are invariant are such that the Macdonald polynomials do not exist in general at these specializations. We also note that more generally, the non-symmetric generalized eigenfunctions for the Cherednik-Dunkl operators  at the special parameters are studied in \cite{CherNonsemisimple}. Further comparison with \cite{CherNonsemisimple} may be important for extracting the information on the eigenfunctions of the operators $\overline{M_\lambda}$. At generic values of the parameters for the restricted operators the eigenfunctions might have different structure and description. Thus in type $A$ a complete set of eigenfunctions (so called super Macdonald polynomials) is constructed in \cite{SV} at generic $\eta, \hbar$, these polynomials are symmetric in each set of variables $x_i$, $y_j$ and span a subspace in the space of functions of $x_i, y_j$ with simple description.

Another important question is on the composition series of the polynomial representation. In the case of $R=A_N$ it is given in terms of the ideals defined by vanishing (``wheel") conditions (\cite{KasA}, \cite{Enomoto}). It would be interesting to see if the ideals $I_{\cal D}$ can be used to define the composition factors in some other cases. Further, the vanishing conditions of the submodules in the polynomial representation for the rational Cherednik algebras were helpful to establish unitarity of certain submodules \cite{FS}. It would be interesting to see if analogous approach can be applied for the case of DAHA.

\appendix

\section{Other affine root systems}
\label{sec6}

The Macdonald-Ruijsenaars systems and their generalized versions can be defined starting from a general affine root system. The cases of affine root systems not considered in the main part of the paper lead to particular cases of the operators $\overline{M_\lambda}$ from Section~\ref{sec5}. Nevertheless we present these cases for completeness.

Let $\SS\subset\wh V$ be an irreducible affine root system~\cite{Mac1,Mac2} and let $\wt W=\wt W(\SS)$ be the group generated by $\{s_\alpha\mid\alpha\in\SS\}$ -- {\it affine Weyl group for the affine root system $\SS$}. Let $\Delta_1=\{a_0,a_1,\ldots,a_n\}$ be a basis of $\SS$ and let
\begin{align}
 \bar a_i=
\begin{cases}
 a_i &  \text{if $2a_i\notin\SS$} \\
 2a_i & \text{if $2a_i\in\SS$}
\end{cases}, &&& \Delta_2=\{\bar a_0,\bar a_1,\ldots,\bar a_n\}.
\end{align}
The sets $\SS_1=\{\alpha\in\SS\mid\alpha/2\notin\SS\}$ and $\SS_2=\{\alpha\in\SS\mid2\alpha\notin \SS\}$ are reduced affine root subsystems in $\wh V$ and $\Delta_1$ and $\Delta_2$ are bases of $\SS_1$ and $\SS_2$ respectively. The group $\wt W$ is generated by the {\it affine simple reflections} $s_i=s_{a_i}=s_{\bar a_i}$.

Let $t\colon\SS_2\to\mathbb C\backslash\{0\}$ and $u\colon\SS_2\to\mathbb C\backslash\{0\}$ be $\wt W$-invariant functions such that $u_\a=u(\alpha)=1$ for all $\alpha\in\SS_1\cap\SS_2$. Let $t_i=t_{\bar a_i}$, $u_i=u_{\bar a_i}$, $i=0,1,\ldots,n$. The functions $t$ and $u$ are completely determined by these values and the restriction condition on $u$ can be reformulated as $u_i=1$ if $a_i=\bar a_i$.
Consider the operators
\begin{align}
  T_i=t_i+t_i^{-1}\frac{(1-t_iu_ie^{-\bar a_i(x)/2})(1+t_iu_i^{-1}e^{-\bar a_i(x)/2})}{1-e^{-\bar a_i(x)}}\,(s_i-1), &&& i=0,1,\ldots,n. \label{TiGeneral}
\end{align}
These operators preserve the space $\mathcal F\subset C^\infty(V)$ of smooth functions $f(x)$ on $V$ such that $f(x)=f\big(x+2\pi i(\bar a_j^{(1)})^\vee\big)$ for all $j=0,1,\ldots,n$.

\begin{Prop} \label{Gtt}
 The operators~\eqref{TiGeneral} satisfy the affine Hecke algebra relations~\eqref{DAHAtt} and \eqref{DAHAttt} where $m_{ij}$ is the order of the element $s_is_j\in\wt W$.
\end{Prop}

The Proposition follows from the established cases $\SS=(C_n^\vee,C_n)$ and $\SS=\wh{R}$ where $R$ is a reduced root system.

Now consider the affine root systems not covered in the previous sections. These are $\SS=(B_n, B_n^\vee)$, $\wh{BC_n}$, $(BC_n, C_n)$ and $(C_n^\vee,BC_n)$. In the above notations we have $\SS_1= \wh{B_n}$, $\SS_2=\wh{B_n}^\vee$ for $\SS=(B_n, B_n^\vee)$, $\SS_1=\SS_2=\wh{BC_n}$ for $\SS=\wh{BC_n}$; $\SS_1= \wh{BC_n}$, $\SS_2=\wh{C_n}$ for $\SS=(BC_n, C_n)$;
$\SS_1= \wh{C_n}^\vee$, $\SS_2=\wh{BC_n}$ for $\SS=(C_n^\vee, BC_n)$. We consider the standard realizations of these affine root systems, so in particular, $\SS\supset D_n=\{\pm e_i\pm e_j\mid1\le 1<j\le n\}$. Define the lattice $L_Y=\mathbb Z^n=P^\vee_{B_n}=Q^\vee_{C_n}=P^\vee_{D_n}$, where $P^\vee_{B_n}$, $Q^\vee_{C_n}$ and $P^\vee_{D_n}$ are the coweight lattice for $R=B_n$, the coroot lattice for $R=C_n$ and the coweight lattice for $R=D_n$ respectively. Let $H^{L_Y}_\SS$ be the {\it affine Hecke} algebra generated by the operators $T_{\hat w}=T_{\hat w}^\SS$ with $\hat w\in W\ltimes\tau(L_Y)$ defined by the formulae~\eqref{Twdef}, \eqref{TiGeneral}. So we consider `extended' affine Hecke algebra for $\SS=(B_n,B_n^\vee)$ and `non-extended' for other cases as the group $W\ltimes\tau(L_Y)$ is the extended affine Weyl group $\wh W=\Pi\ltimes\wt W(\SS)$ for $\SS=(B_n,B_n^\vee)$ and it is the usual affine Weyl group $\wt W=\wt W(\SS)$ for other $\SS$. We will also consider the `non-extended' affine Hecke algebras $H^{L_Y}_\SS$ for $\SS=\wh{C_n}$ and $\wh{C_n}^\vee$, where $L_Y$ is as above. The Cherednik-Dunkl operators $Y^\lambda=Y^\lambda_\SS$, $\lambda\in L_Y$, are defined by the formula $Y^{\mu-\nu}=T_{\tau(\mu)}T_{\tau(\nu)}^{-1}$, where $\mu,\nu\in P^\vee_+\cap L_Y$. These operators form the commutative subalgebra $\mathbb C[Y^{L_Y}_\SS] \subset H^{L_Y}_\SS$.

The elements $m_\lambda(Y)=\sum_{\nu\in W\lambda}Y^\nu$ preserve the space of invariants $\mathcal F^W$, their restrictions $M_\lambda=m_\lambda\big|_{{\cal F}^W}$ form commutative algebra $\mathfrak R^{L_Y}$ of difference operators.
The simplest operators in this algebra has the form
\begin{align}
 M_{-\nu}=\frac1{|W_{-\nu}|}\sum_{w\in W}g^{-\nu}(w^{-1}x)\Big(\tau(-w\nu)-1\Big)+ const, \label{RMop_brthetaGeneral}
\end{align}
where $\nu$ is a minuscule coweight $b_r$ or $\nu=\theta^\vee$, and
\begin{align}
 g^{-\nu}(x)=\prod_{\alpha\in\SS_2(-\nu)}\frac{(1-t_\alpha u_\alpha e^{\alpha(x)/2})(1+t_\alpha u_\alpha^{-1} e^{\alpha(x)/2})}{t_\alpha (1-e^{\alpha(x)})},
\end{align}
where $\SS_2(-\nu)=\SS_2^-\cap\big(\tau(-\nu)\SS_2^+\big)$, $\SS_2^+=-\SS_2^-=\{\alpha\in\SS_2\mid\alpha(x)>0 \text{ whenever $\bar a_i(x)>0$ } \forall i=0,1,\ldots,n\}$.

The following diagram explains how one gets the subalgebras $\mathbb C[Y^{L_Y}_\SS] \subset H^{L_Y}_\SS$ by specializing the parameters starting from the subalgebra generated by the Cherednik-Dunkl operators in the $(C_n^\vee,C_n)$ DAHA. We also include in the diagram the affine Hecke algebras parts of DAHA of $\wh{R}$ and $\wh{R}^\vee$ type where $R=B_n, D_n$.

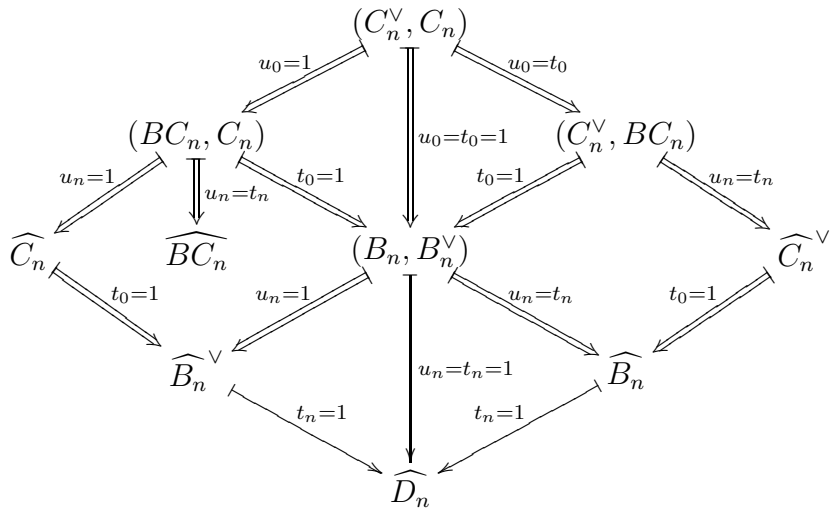
\begin{align}
\xymatrix{
&& \,(C_n^\vee,C_n)\,\ar@{|=>}[dl]_{u_0=1}\ar@{|=>}[dd]^{u_0=t_0=1}\ar@{|=>}[dr]^{u_0=t_0} \\
&(BC_n,C_n)\ar@{|=>}[dl]_{u_n=1}\ar@{|=>}[d]^{u_n=t_n}\ar@{|=>}[dr]^{t_0=1}&&(C_n^\vee,BC_n)\ar@{|=>}[dl]_{t_0=1}\ar@{|=>}[dr]^{u_n=t_n} \\
\wh{C_n}\ar@{|=>}[dr]^{t_0=1} &\wh{BC_n}& (B_n,B_n^\vee)\ar@{|=>}[dl]_{u_n=1}\ar@{|->}[dd]^{u_n=t_n=1}\ar@{|=>}[dr]^{u_n=t_n} & & \wh{C_n}^\vee\ar@{|=>}[dl]_{t_0=1} \\
& \wh{B_n}^\vee\ar@{|->}[dr]^{t_n=1} & &\wh{B_n}\ar@{|->}[dl]_{t_n=1} \\
&& \wh{D_n}
} \notag
\end{align}

The double arrows in the diagram mean the isomorphism of the corresponding Cherednik-Dunkl subalgebras under the specialization and the single arrows mean the embedding of these subalgebras. The parameters $t_0, t_n, u_0, u_n$ are the parameters of the operators~\eqref{TiGeneral} for $\SS=(C_n^\vee, C_n)$. The validity of these isomorphisms essentially follows from the following two lemmas.

\begin{Lem} \label{lemCCtoBB}
 If $t_0=u_0=1$ then the affine Hecke algebra $H^{Q^\vee}_{(C_n^\vee,C_n)}$ coincides with the affine Hecke algebra $H^{P^\vee}_{(B_n,B_n^\vee)}$. More exactly, under the specialization the generators correspond as
\begin{align}
 T_i^{(C_n^\vee, C_n)}&=T_i^{(B_n, B_n^\vee)}, & T_0^{(C_n^\vee, C_n)}&=\pi_1^{(B_n, B_n^\vee)},
\end{align}
where $i=1,\ldots,n$, and the Cherednik-Dunkl operators correspond as $Y^\lambda_{(C_n^\vee,C_n)}= Y^\lambda_{(B_n, B_n^\vee)}$, $\forall \lambda \in L_Y$.
\end{Lem}

\noindent{\bf Proof.} The affine Hecke algebra $H^{Q^\vee}_{(C_n^\vee,C_n)}$ is generated by $T_0^{(C_n^\vee,C_n)},T_1^{(C_n^\vee,C_n)},\ldots,T_n^{(C_n^\vee,C_n)}$ while the algebra $H^{P^\vee}_{(B_n,B_n^\vee)}$ is generated by $\pi_1^{(B_n,B_n^\vee)},T_1^{(B_n,B_n^\vee)},\ldots,T_n^{(B_n,B_n^\vee)}$. The operators $T_i^{(C_n^\vee,C_n)}$ and $T_i^{(B_n,B_n^\vee)}$ coincide for each $i=1,\ldots,n$ (at any parameters) so denote them by $T_i$. The operator $T_0^{(C_n^\vee,C_n)}$ at the specialization of parameters becomes $s_0^{(C_n^\vee,C_n)}=\tau(e_1)s_{e_1}=\pi_1^{(B_n,B_n^\vee)}$. Thus we are left to prove the formula $Y^\lambda_{(C_n^\vee,C_n)}=Y^\lambda_{(B_n,B_n^\vee)}$ at this specialization. It is enough to check this formula for $\lambda=e_i$, $i=1,\ldots,n$. From the definition one yields $Y^{e_1}_{(C_n^\vee,C_n)}=T_0^{(C_n^\vee,C_n)}T_{s_{e_1}}^{(C_n^\vee,C_n)}$ and $Y^{e_1}_{(B_n,B_n^\vee)}=\pi_1^{(B_n,B_n^\vee)}T_{s_{e_1}}^{(B_n,B_n^\vee)}$. Then for $i=2,\ldots,n$ the equality $Y^{e_i}_{(C_n^\vee,C_n)}=Y^{e_i}_{(B_n,B_n^\vee)}$ follows by induction from the formula $Y^{e_i}_{\SS}=T_{i-1}^{-1}Y^{e_{i-1}}_{\SS}T^{-1}_{i-1}$ valid for both $\SS=(C_n^\vee,C_n)$ and $\SS=(B_n,B_n^\vee)$. \qed

Consider now $\SS=\wh{D_n}$, $n\ge2$. The simple roots of $D_n$ are $\alpha_i=e_i-e_{i+1}$ for $i=1,\ldots,n-1$ and $\alpha_n=e_{n-1}+e_n$. Denote by $\pi_r^{D_n}$, $T_i^{D_n}$  the corresponding elements of the algebra $H^{P^\vee}_{\wh{D_n}}$ where $r\in{\cal O}=\{0,1,n-1,n\}$, $i=0,1,\ldots,n$.

Note that if $t_n=1$ then the affine Hecke algebra $H^{P^\vee}_{\wh{B_n}}$ in the representation~\eqref{ReprDAHA} coincides with the affine Hecke algebra $H^{P^\vee}_{\wh{B_n}^\vee}$ in the representation~\eqref{ReprDAHAvee}. Indeed,
the algebras $H^{P^\vee}_{\wh{B_n}}$, $H^{P^\vee}_{\wh{B_n}^\vee}$ have the only different generator $T_n$ which becomes $T_n=s_{e_n}$ at $t_n=1$. Since the corresponding operators $\pi_1,T_i$ of the algebras $H^{P^\vee}_{\wh{B_n}}$ and $H^{P^\vee}_{\wh{B_n}^\vee}$
 coincide at this specialization we denote them by $\pi_1^{B_n}$, $T_i^{B_n}$ where $i=0,\ldots,n$.

\begin{Lem} \label{LemYBnDn}
At the specialization $t_n=1$ one has
\begin{align}
  Y^\lambda_{\wh{B_n}}=   Y^\lambda_{\wh{B_n}^\vee} = Y^\lambda_{\wh{D_n}} \label{YBnYDn}
\end{align}
for all $\lambda\in L_Y$. That is the corresponding algebras are embedded at this specialization: $\mathbb C[Y^{P^\vee}_{\wh{B_n}}]=\mathbb C[Y^{P^\vee}_{\wh{B_n}^\vee}]\subset\mathbb C[Y^{P^\vee}_{\wh{D_n}}]$.
\end{Lem}

\noindent{\bf Proof.} Note first that the operators $T_i^{B_n}$ and $T_i^{D_n}$ coincide for $i=1,\ldots,n-1$ (at any parameters) and we denote them by $T_i$. Since $T_n^{B_n}=s_{e_n}$ (at the specialization $t_n=1$) we obtain
\begin{align}
 \pi_1^{D_n}=\tau(e_1)s_{e_1}s_{e_n}=\pi_1^{B_n}T_n^{B_n}, &&& T_n^{D_n}=s_nT_{n-1}s_n=(T_n^{B_n})^{-1}T_{n-1}T_n^{B_n}. \label{AHADnBnEmbed}
\end{align}
Comparing the expressions
\begin{align}
 Y^{e_1}_{B_n}&=\pi_1^{B_n}T_1\cdots T_{n-2}T_{n-1}T_n^{B_n}T_{n-1}\cdots T_1, \\
 Y^{e_1}_{D_n}&=\pi_1^{D_n}T_1\cdots T_{n-2}T_n^{D_n}T_{n-1}\cdots T_1
\end{align}
we derive $Y^{e_1}_{B_n}=Y^{e_1}_{D_n}$.
Due to the formula $Y^{e_i}_{R}=T_{i-1}^{-1}Y^{e_{i-1}}_{R}T^{-1}_{i-1}$ valid for both $R=B_n$ and $R=D_n$ we establish the equality~\eqref{YBnYDn} for all $\lambda=e_i$, where $i=1,\ldots,n$, and, hence, for all $\lambda\in L_Y$. \qed \\

It follows that the operators $M_\lambda$ for the affine root system $\SS$ coincide with the operators $M_\lambda$ from Section~\ref{sec5} at the corresponding specialization. In particular the specializations of the Koornwinder operator~\eqref{MacKoor} give the operators $M_{-e_1}$ determined by~\eqref{RMop_brthetaGeneral}. The specializations of the generalized Koornwinder operator~\eqref{MacKoorGeneralized} give the generalized Macdonald-Ruijsenaars operators $\overline{M_{-e_1}}$ for the affine root system $\SS$. Each specialization of the operators $\overline{M_\lambda}$ from Section~\ref{sec5} gives commutative algebra containing the operator $\overline{M_{-e_1}}$. These algebras were constructed in Sections~\ref{sec4},~\ref{sec62} for $\SS = \wh{B_n}, \wh{D_n}, \wh{B_n}^\vee$. For other cases the invariant ideals and these operators can also be defined using the specialization of the invariant ideals for $\SS=(C_n^\vee, C_n)$.

For example consider the case $\SS=(B_n, B_n^\vee)$, so $\SS_2=\wh{B_n}^\vee$.
Let $\eta\colon B_n^\vee\to\mathbb C$ be a $W$-invariant function and let $S_0$ be a subset of the set of simple roots of $B_n$.  Define the ideal $I_{\cal D}$ by the formulae~\eqref{D0defvee} and \eqref{DandIdefvee}. Let $S_0=\bigsqcup_{\ell=1}^m S_\ell$ be the decomposition to the irreducible components $S_1,\ldots,S_m\subset S_0$. Let $\check h_\eta^\ell$ be the corresponding generalized Coxeter number defined by~\eqref{genCoxcheck}.
Let $\check{\cal F}$ be the space defined in Subsection~\ref{sec62} for $R=B_n$. Then the following statement takes place.

\begin{Th} \label{ThInvBB}
  Suppose that $t_1^2=e^{\eta_1}$ if $n>1$, $\epsilon_nt_nu_n^{\epsilon_n}=e^{\eta_n/2}$ for some $\epsilon_n=\pm1$ and $\check h^\ell_\eta=\hbar$ for all $\ell=1,\ldots,m$. Then the ideal $I_{\cal D}\subset\check{\cal F}$ is invariant under the actions of the operators~\eqref{TiGeneral} and $\pi_1$. The corresponding generalized Macdonald-Ruijsenaars operator $\overline M_{-e_1}$ can be obtained by the specialization of the parameters in the generalized Koornwinder operator~\eqref{MacKoorGeneralized}.
\end{Th}

\end{document}